\documentclass[10pt]{article}

\usepackage[a4paper]{geometry}
\usepackage{indentfirst}
\usepackage{hyperref}
\usepackage{units}
\usepackage{amssymb, amsmath, amsthm}
\usepackage{graphicx}
\usepackage[labelsep=none]{caption}
\usepackage[labelformat=simple]{subcaption}

\newcommand{\reveq}{\stackrel{-1}{=}}
\newcommand{\revsim}{\stackrel{-1}{\sim}}
\newcommand{\sto}{{\to}} 
\newcommand{\tto}{\mathbin{\rightsquigarrow}} 
\newcommand{\eqdir}{\approx}
\newcommand{\opdir}{\stackrel{-1}{\eqdir}}
\newcommand{\shift}[1]{\text{\rm \textcircled{\ensuremath{#1}}}}
\newcommand{\mymod}[1]{\mkern 8mu (\mathrm{mod} \mkern 6mu #1)}

\newcommand{\llfloor}{\lfloor\!\lfloor}
\newcommand{\rrfloor}{\rfloor\!\rfloor}
\newcommand{\llceil}{\lceil\!\lceil}
\newcommand{\rrceil}{\rceil\!\rceil}
\DeclareSymbolFont{bbold}{U}{bbold}{m}{n}
\DeclareMathSymbol{\lopen}{\mathopen}{bbold}{`(}
\DeclareMathSymbol{\ropen}{\mathclose}{bbold}{`)}

\newcommand{\half}{\nicefrac{1}{2}}
\newcommand{\bfp}{\mathbf{p}}
\newcommand{\bfq}{\mathbf{q}}
\newcommand{\bfu}{\mathbf{u}}
\newcommand{\bfv}{\mathbf{v}}

\newcommand{\bfs}{\mathbf{s}}
\newcommand{\bft}{\mathbf{t}}
\newcommand{\olu}{\overline{\mathbf{u}}}

\newcommand{\ols}{\overline{\mathbf{s}}}
\newcommand{\olt}{\overline{\mathbf{t}}}

\newcommand{\frf}{\mathfrak{f}}
\newcommand{\bfphi}{\boldsymbol\varphi}
\newcommand{\bfalpha}{\boldsymbol\alpha}
\newcommand{\bfdelta}{\boldsymbol\delta}
\newcommand{\bfeta}{\boldsymbol\eta}
\newcommand{\bfmu}{\boldsymbol\mu}
\newcommand{\bfrho}{\boldsymbol\varrho}
\newcommand{\bftau}{\boldsymbol\tau}
\newcommand{\bfomega}{\boldsymbol\omega}

\newcommand{\pfloat}{\mathbb{P}_\text{\rm Float}}
\newcommand{\tmptile}{\Theta_\text{\rm Tile}}
\newcommand{\tmpcube}{\Theta_\text{\rm Cube}}
\newcommand{\tmploop}{\Theta_\text{\rm Loop}}
\newcommand{\tmphex}{\Theta_\text{\rm Hex}}
\newcommand{\fodd}{f_\text{\rm Odd}}
\newcommand{\feven}{f_\text{\rm Even}}
\newcommand{\godd}{g_\text{\rm Odd}}
\newcommand{\geven}{g_\text{\rm Even}}
\newcommand{\gso}{g^\star_\text{\rm Odd}}
\newcommand{\gse}{g^\star_\text{\rm Even}}
\newcommand{\gsa}{g^\star_\text{\rm Any}}
\newcommand{\psiall}{\Psi_\text{\rm All}}
\newcommand{\psimain}{\Psi_\text{\rm Main}}
\newcommand{\psisymm}{\Psi_\text{\rm Symm}}

\DeclareMathOperator{\word}{word}
\DeclareMathOperator{\wind}{wind}
\DeclareMathOperator{\area}{area}
\DeclareMathOperator{\myneg}{neg}
\DeclareMathOperator{\mystar}{star}

\hypersetup{bookmarksnumbered, pdfstartview=FitH, pdfpagemode=UseNone}

\newcommand{\mynewtheorem}[2]{\newtheorem{#1}{\indent #2}}
\mynewtheorem{conjecture}{Conjecture}
\mynewtheorem{lemma}{Lemma}
\mynewtheorem{proposition}{Proposition}
\mynewtheorem{theorem}{Theorem}
\newenvironment{myproof}[1][Proof]{\begin{proof}[\indent #1]}{\end{proof}}

\begin{document}

\title{\textbf{Double Tiles}}
\author{Nikolai Beluhov}
\date{}

\maketitle

\begin{abstract} Which polygons admit two (or more) distinct lattice tilings of the plane? We call such polygons double tiles. It is well-known that a lattice tiling is always combinatorially isomorphic either to a grid of squares or to a grid of regular hexagons. We focus on the special case of the double tile problem where both tilings are in the square class. For this special case, we give an explicit description of all double tiles. We establish the result for polyominoes first; then, with little additional effort, we extend the proof to general polygons. Central to the description is a certain finite set of transformations which we apply iteratively to a base shape in order to obtain one family of ``fractal-like'' polyominoes. The double tiles are then given by these polyominoes together with particular ``deformations'' of them. \end{abstract}

\tableofcontents

\section{Introduction} \label{intro}

\subsection{Overview} \label{over}

A polygon $T$ \emph{tiles} the plane when it is possible to partition the plane into congruent copies of $T$. Each copy is called a \emph{tile}, and the partitioning is called a \emph{tiling}.

The simplest kind of tilings are the ones where, roughly speaking, ``each tile looks the same''. Formally, for each pair of tiles $T'$ and $T''$, we want some translation to map $T'$ onto $T''$, while also preserving the tiling as a whole. Clearly, all such translations, taken over all pairs of tiles, must form a lattice. So tilings of this kind are known as \emph{lattice tilings}.

For example, both the square and the regular hexagon admit lattice tilings of the plane. In fact, the square admits multiple lattice tilings as well as multiple non-lattice tilings. On the other hand, the regular hexagon's tiling is unique.

Here is a natural question: Which polygons admit two (or more) different lattice tilings of the plane? What can we say about them? This is a generalisation of an earlier question posed by Mebane in \cite{M}; see below. We call such polygons \emph{double tiles}.

This question appears to be difficult in full generality. We will focus on one natural special case, as follows: It can be shown that, in a lattice tiling, either each tile has $4$ neighbours and the tiling as a whole is ``combinatorially isomorphic'' to a grid of squares; or else each tile has $6$ neighbours and the tiling as a whole is combinatorially isomorphic to a grid of regular hexagons. Our special case is the one where the polygon admits two different $4$-neighbour lattice tilings. Or, in other words, two lattice tilings both of which belong to the square-grid combinatorial-isomorphism class.

\begin{figure}[ht] \centering \includegraphics{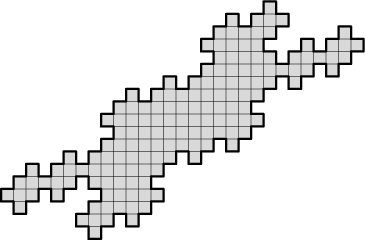} \caption{} \label{t} \end{figure}

\begin{figure}[ht] \null \hfill \begin{subfigure}[c]{175pt} \centering \includegraphics{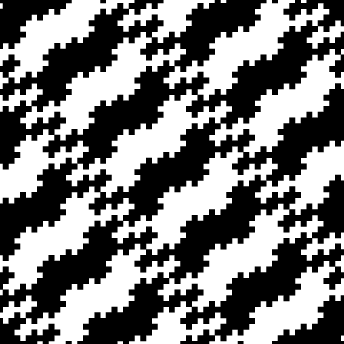} \caption{} \end{subfigure} \hfill \begin{subfigure}[c]{175pt} \centering \includegraphics{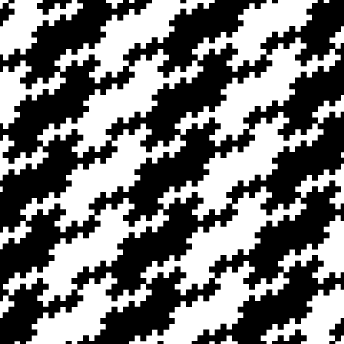} \caption{} \end{subfigure} \hfill \null \caption{} \label{tilings} \end{figure}

For example, consider the polygon in Figure \ref{t}. Two different $4$-neighbour lattice tilings of the plane based on it are shown in Figure \ref{tilings}.

Our main result is an explicit description of all $4$-neighbour double tiles. We go on now to an informal sketch of this description.

One class of polygons which will be particularly important to us is that of polyominoes. A \emph{polyomino} is a polygon formed out of unit squares joined up by their sides. For example, the polygon in Figure \ref{t} is a polyomino made out of $153$ unit squares.

We consider the problem first in the settings of polyominoes. So, to begin with, we are going to describe all polyominoes which are $4$-neighbour double tiles. The reason why we single out this class of polygons in particular will become clear soon.

Central to the description is one finite set of transformations over polyominoes. Roughly speaking, each one of these transformations takes in a polyomino and produces as output a larger, more complicated polyomino. When several such transformations are applied in a series (so that the output of each is fed into the next), usually the result is a ``fractal-like'' shape with a densely zigzagging boundary.

\begin{figure}[ht] \null \hfill \begin{subfigure}[c]{50pt} \centering \includegraphics{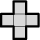} \caption{} \label{step-a} \end{subfigure} \hfill \begin{subfigure}[c]{50pt} \centering \includegraphics{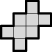} \caption{} \label{step-b} \end{subfigure} \hfill \begin{subfigure}[c]{50pt} \centering \includegraphics{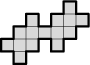} \caption{} \label{step-c} \end{subfigure} \hfill \null\\ \vspace{\baselineskip}\\ \null \hfill \begin{subfigure}[c]{135pt} \centering \includegraphics{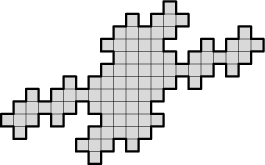} \caption{} \label{step-d} \end{subfigure} \hfill \null \caption{} \label{steps} \end{figure}

We always begin such a series with the same base shape, namely the \emph{Greek cross}. This is the polyomino formed by joining four unit squares to the sides of a fifth, central one; it is shown in Figure \ref{step-a}. By way of an example, there exists a series of four transformations which turns the Greek cross into the polyomino of Figure \ref{t}; the intermediate steps are shown in Figures \ref{step-b}--\subref{step-d}.

Eventually, we will see that every polyomino $\Phi$ obtained from the Greek cross by means of our set of transformations is a $4$-neighbour double tile.

However, this is not the complete list of solutions. There is one more part to the description. It hinges on the notion of a ``deformation'', which we now proceed to explain.

\begin{figure}[ht] \null \hfill \begin{subfigure}[c]{75pt} \centering \includegraphics{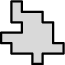} \caption{} \label{def-a} \end{subfigure} \hfill \begin{subfigure}[c]{135pt} \centering \includegraphics{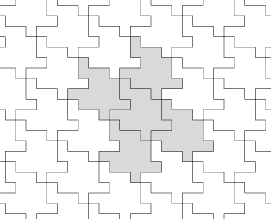} \caption{} \label{def-b} \end{subfigure} \hfill \begin{subfigure}[c]{75pt} \centering \includegraphics{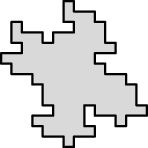} \caption{} \label{def-c} \end{subfigure} \hfill \null \caption{} \label{deff} \end{figure} 

Choose any polyomino $T$ which admits a $4$-neighbour lattice tiling. One example is shown in Figure \ref{def-a}, and the tiling is in Figure \ref{def-b}. Since each such tiling is combinatorially isomorphic to a unit-square grid, we can pick several tiles within the tiling so that they are joined together in the exact same manner as the cells of $\Phi$. Figure \ref{def-b} shows this for the Greek cross; the five chosen tiles are highlighted. The union of the chosen tiles will be a new polyomino which we call a \emph{deformation} of $\Phi$; for our $T$ and the Greek cross, this is shown in Figure \ref{def-c}.

It is straightforward to see that, if $\Phi$ is a $4$-neighbour double tile, then so are all of its deformations, too. Thus, to the special family of polyominoes specified above, we must add all of their deformations as well.

Our first main theorem states that this, now, is the complete list -- in the polyomino setting. The extension from polyominoes to general polygons is only going to require a small amount of additional effort.

Consider the following simple observation: In the definition of a deformation, there is no need to assume that $T$ is a polyomino. We can instead let $T$ be any polygon which admits a $4$-neighbour lattice tiling, and the definition will still make sense. Hence, out of our special family of polyominoes, by means of polygonal deformations we can obtain a lot of polygonal $4$-neighbour double tiles as well.

Our second main theorem states that the latter is, in fact, the complete list of all $4$-neighbour double tiles. So, in particular, the polyominoes form a class of special significance to the $4$-neighbour double tile problem. In some sense, the ``fundamental'' solutions are polyominoes; and all of the polygonal solutions can be derived from a certain family of such fundamental solutions.

Our approach is based on encoding polygons as words over finite alphabets. For example, to encode a polyomino, we go around its boundary once, and each time we make a step to the right, left, up, or down, we write down the letter $\mathsf{R}$, $\mathsf{L}$, $\mathsf{U}$, or $\mathsf{D}$, respectively.

This sort of relationship between words and polygons is a familiar idea in topology and group theory; \cite{S} goes into both the mathematics and the history. The same insight has also been fruitfully applied to purely combinatorial questions about tilings; see, e.g., \cite{CL}.

\subsection{Structure of the Paper} \label{struct}

The rest of the paper is structured as follows:

In Section \ref{init}, we go into deeper detail regarding our encoding of polygons with words. Usually (e.g., in \cite{S} and \cite{CL}), the boundaries of polygons are interpreted as group products, and the reasoning revolves around the group-theoretic properties of these products. By contrast, we interpret the boundaries of polygons as literal words, and we care about the ``syntactic'' properties of these words.

Section \ref{chain} applies the general tools developed in Section \ref{init} to the problem of $4$-neighbour double tiles in particular. We reduce this problem to the study of certain systems of equations over words.

However, the reduction comes at a cost. A word could describe the boundary of a polygon just as easily as some tangled-up self-intersecting trajectory. In the word setting, the double-tile conditions are easy to express -- but, on the flip side, words alone are not particularly sensitive to self-intersections.

So, in the transition from geometry to syntax, the ``balance of difficulties'' is turned on its head. The analysis of the double-tile conditions will be, for the most part, neat and tidy. But it will be trickier to filter out the ``double-tile-like'' objects which are self-intersecting trajectories instead of genuine polygons.

In Section \ref{lift}, we introduce the finite set of transformations referred to above. The definitions are entirely in terms of syntactic operations over words. However, because of the connection between words and polygons, of course the same definitions can be understood in geometric terms, too.

This is also the section where we are first able to state our main result. We begin with the special case when the double tile is a polyomino, and we only extend the result's statement to general polygons further on in the paper. This choice reflects the logical structure of the proof; almost all of the meaningful work is done in the polyomino setting.

The result splits naturally into a negative part and a positive part. The negative part says that shapes outside of our description cannot possibly work, and the positive part says that all shapes included in the description do in fact work. (Here, as above, the notion of a shape ``working'' is less about satisfying the double-tile conditions, and more about avoiding self-intersections.) The two parts are very different, conceptually, and we handle them separately.

The proof of the negative part is in Section \ref{neg}. It is mostly syntactic, and it revolves around manipulations with words.

Section \ref{prop} is a brief interlude where we work out some of the most basic properties of the shapes included in the description. It provides the foundation required for Sections \ref{pos} and \ref{prime}.

Following this, the proof of the positive part is in Section \ref{pos}. It is mostly geometric, and it revolves around manipulations with figures in the plane.

By this point, we will have resolved the polyomino case. In Section \ref{poly}, we state and prove the generalisation to polygons. The main idea is to show that $4$-neighbour double tiles are ``alphabetisable''. This means that, even if a $4$-neighbour double tile is not a polyomino, we can still encode it as a word over some finite alphabet -- possibly much larger than the alphabet $\{\mathsf{R}, \mathsf{L}, \mathsf{U}, \mathsf{D}\}$ which suffices for polyominoes. Then our analysis of the polyomino setting will carry over immediately.

As we remarked in the beginning, this work developed out of an earlier question due to Mebane. It is about the special case when the double tile is a polyomino and its two tilings are patterned after the two lattices formed by the complex multiples of a given Gaussian prime and its conjugate. We touch upon the subject briefly in Section \ref{prime}.

Finally, Section \ref{further} lists some promising directions for further research.

\subsection{Survey of Earlier Literature} \label{lit}

Following the submission of the original manuscript, a referee pointed out the existence of relevant literature which the author was previously unaware of. This subsection has been added in to address the matter, and the rest of the text has been left unaltered.

There are some overlaps between the present work and the earlier \cite{BN}--\cite{BMGL}. However, they involve solely the basic notions and tools -- not the main results. With a small number of exceptions, these overlaps are confined to Sections \ref{word}--\ref{rl}. (Our main results and their proofs occupy Sections \ref{smt}--\ref{adc} instead.) The details will be spelled out below.

Note that \cite{BN}--\cite{BMGL} are concerned exclusively with polyominoes, and not with general polygons. Here, this will be understood implicitly and we will not remark on it every time. E.g., when we refer to a ``double tile'', we will really mean a double tile in the polyomino setting; when we cite a reference for some result, we will really mean the polyomino variant of that result; and so on.

The rest of this subsection makes use of notation and terminology which have not been defined yet. For this reason, the reader is advised to postpone reading it at least until the end of Section~\ref{rl}.

The $6$-neighbour analogue of Lemma \ref{abl} at the beginning of Section \ref{six} can be found in \cite{BN}. Clearly, the original Lemma \ref{abl} is an immediate corollary. Building on these results, \cite{BP} contains an equivalent form of Lemma \ref{nil}. The proofs differ from ours.

The notion of a double tile is introduced in \cite{P}, with a focus on $4$-neighbour double tiles in particular. The same work includes a table of all $4$-neighbour double tiles of perimeter at most~$32$. Two conjectures formulated in \cite{P} have played a central role in the subsequent development of the subject. Both of them will be discussed shortly.

Two infinite families of $4$-neighbour double tiles are exhibited in \cite{BMBGL}. It is straightforward to see that the Christoffel tiles are $\{f_1, f_2\}$-descendants of the Greek cross, and the Fibonacci tiles are its $\{\gse\}$-descendants. Of course, the definitions of \cite{BMBGL} are not phrased in terms of descents; instead, they are based on the theory of Christoffel and Fibonacci words.

The first one of the two aforementioned conjectures of \cite{P} is that a polyomino cannot admit three distinct $4$-neighbour lattice tilings. Or, equivalently, a $4$-neighbour double tile must necessarily admit exactly two $4$-neighbour lattice tilings. (Note that questions of this form are not considered in the present work.) This was confirmed in \cite{BMBL}.

Before we can talk about the second conjecture, we must set up some terminology. The work \cite{P} introduces the notion of a deformation, and it defines a polyomino to be \emph{prime} when it cannot be obtained by means of a deformation in a nontrivial way. The conjecture is that all prime $4$-neighbour double tiles are centrally symmetric. This was confirmed in \cite{BMGL}. It is also a direct corollary of our Lemma \ref{symml}.

We proceed now to a detailed discussion of the overlaps between the tools developed in \cite{BMGL}, and our own machinery. The analysis of $4$-neighbour double tiles in \cite{BMGL} begins, as does ours, with the system of equations (\textbf{I}) over the words of $\mathcal{A}_\square$.

Next come three transformations called $\textsc{trim}$, $\textsc{extend}$, and $\textsc{swap}$. The first two have the same effect as our $f$-reductions and $f$-lifts. The last one has the same effect as our $g$-transforms, up to cyclic shift, in all instances where the $g$-transforms are well-defined. More generally, if $U$ and $V$ are related by a $\godd$-transform, then up to cyclic shift (and disregarding some degenerate cases) $\textsc{swap}$ also swaps $f_1^i \circ f_3^j(U)$ with $f_1^i \circ f_3^j(V)$; and similarly for $f_2$, $f_4$, and~$\geven$. The definitions in \cite{BMGL} are, once again, different from ours.

This is followed by a proof that every $4$-neighbour double tile can be reduced to a root using $\textsc{trim}$ and $\textsc{swap}$. Or, conversely, every $4$-neighbour double tile can be obtained from some root by means of $\textsc{extend}$ and $\textsc{swap}$. These results roughly match the material in our Section \ref{rl}, with one augmentation: Initially, we also allow the possibility of reduction to a loop, to be ruled out in our Section \ref{loop}.

In addition, \cite{BMGL} exhibits one algorithm for the generation of all $4$-neighbour double tiles. It operates by first employing $\textsc{extend}$ and $\textsc{swap}$ to generate all solutions to (\textbf{I}) over the polyomino alphabet $\mathcal{A}_\square$, and then separately it filters out the solutions which contain self-intersections. (Hence, it completely elides the ``self-intersectional'' part of the problem. By contrast, the bulk of the present work is devoted to sorting out the self-intersections.)

This algorithm is not too efficient, as the authors of \cite{BMGL} observe themselves. The issue is that the solutions to (\textbf{I}) which are free of self-intersections form a vanishingly small fraction of the total. A better approach would account for the self-intersections ``organically'', so that only the solutions to (\textbf{I}) which do describe polyominoes are generated. It is straightforward to rephrase our Theorem \ref{mt} as a generation algorithm which would indeed be more satisfactory in this sense.

\section{Words and Polylines} \label{init} 

Here, we collect some initial observations regarding words, polylines, and tilings.

Sections \ref{word}--\ref{line} introduce terminology and notation. Then Section \ref{val} details the way we use words to encode polylines. Section \ref{tile} is a quick review of some of the basics of lattice tilings. Finally, in Sections \ref{abab}--\ref{abc} we prove a number of lemmas which will be important to us later on.

\subsection{Words} \label{word}

Let $\mathcal{A}$ be an alphabet; i.e., a set of letters. Throughout, we assume that each alphabet $\mathcal{A}$ we consider comes equipped with some involution over $\mathcal{A}$, called \emph{reversal}. For each letter $w$ of $\mathcal{A}$, we denote its reverse by $w^{-1}$.

There is a natural extension of reversal to all words over $\mathcal{A}$. Explicitly, if $W = w_1w_2 \ldots w_k$ is a word over $\mathcal{A}$, then we define its reverse $W^{-1}$ to be $w_k^{-1}w_{k - 1}^{-1} \ldots w_1^{-1}$. Given two words $U$ and $V$ over $\mathcal{A}$, we write $U \reveq V$ when $U = V^{-1}$; or, equivalently, when $U^{-1} = V$.

The point of reversals will become clear in Section \ref{val}. The gist is that we are going to encode oriented polylines as words, and then reversing the word will correspond to reversing the orientation of the polyline.

One alphabet of particular importance to us will be \[\mathcal{A}_\square = \{\mathsf{R}, \mathsf{L}, \mathsf{U}, \mathsf{D}\}.\]

Its letters stand for ``right'', ``left'', ``up'', and ``down''. This is the alphabet we will use to encode polyominoes. We define reversal over $\mathcal{A}_\square$ in the natural way, by $\mathsf{R} \reveq \mathsf{L}$ and $\mathsf{U} \reveq \mathsf{D}$.

We write $\varepsilon$ for the empty word. We also write $|W|$ for the length of the word $W$.

Let $\operatorname{\#}(w, W)$ denote the number of occurrences of the letter $w$ in the word $W$. We say that two words $U$ and $V$ are \emph{anagrams} when $\operatorname{\#}(w, U) = \operatorname{\#}(w, V)$ for all letters $w$ of $\mathcal{A}$.

We also say that a word $W$ is \emph{balanced} when $\operatorname{\#}(w, W) = \operatorname{\#}(w^{-1}, W)$ for all $w$. Or, equivalently, $W$ is balanced if and only if $W$ and $W^{-1}$ are anagrams. A \emph{combinatorial self-intersection} in $W$ is a balanced subword of $W$ distinct from $\varepsilon$ and $W$; except that, if $W = ww^{-1}$ for some letter $w$ of $\mathcal{A}$, then we also consider $ww^{-1}$ to be a combinatorial self-intersection in $W$. The point of this definition will become clear in Section \ref{val}.

Suppose that $W = w_1w_2 \ldots w_k$. We define the $i$-th \emph{cyclic shift} of $W$, denoted $W^\shift{i}$, to be the word $w_{i + 1}w_{i + 2} \ldots w_kw_1w_2 \ldots w_i$. I.e., the first $i$ letters are moved to the end. When $i \not \in \{0, 1, \ldots, k\}$, we take its remainder modulo $k$; e.g., $\mathsf{R}\mathsf{U}\mathsf{L}\mathsf{D}^\shift{9} = \mathsf{U}\mathsf{L}\mathsf{D}\mathsf{R}$. Given a positive integer $\ell$, we write also $W^\ell$ for the concatenation of $\ell$ copies of $W$.

We denote the space of all words over $\mathcal{A}$ by $\word(\mathcal{A})$.

Let $\mathcal{B}$ be an alphabet with reversal, too -- possibly coinciding with $\mathcal{A}$. We define a \emph{substitution} to be a mapping $\zeta : \word(\mathcal{A}) \to \word(\mathcal{B})$ which agrees with reversal and concatenation, in the sense that $\zeta(W^{-1}) = \zeta(W)^{-1}$ and $\zeta(UV) = \zeta(U)\zeta(V)$ for all words $U$, $V$, $W$ over~$\mathcal{A}$.

Clearly, a substitution $\zeta$ is determined completely by its effect over the letters of $\mathcal{A}$. In fact, just one letter out of each reverse pair is sufficient. Let $v_1$, $v_2$, $\ldots$, $v_\ell$ be some letters of $\mathcal{A}$ such that $\zeta(w) = w$ for all letters of $\mathcal{A}$ outside of $\{v_1, v_1^{-1}, v_2, v_2^{-1}, \ldots, v_\ell, v_\ell^{-1}\}$; and let also $\zeta(v_i) = V_i$. Then instead of $\zeta(W)$ we may write $W[v_1 \to V_1, v_2 \to V_2, \ldots, v_\ell \to V_\ell]$. This notation signifies that $\zeta(W)$ is the word obtained from $W$ by substituting each occurrence of $v_i$ with $V_i$ and each occurrence of $v_i^{-1}$ with $V_i^{-1}$, simultaneously for all $i$.

\subsection{Translations} \label{vec}

The \emph{negation} of a vector $\bfu$ is the vector $-\bfu$.

Let $A$ be some geometric object, such as a point, segment, or polygon. Let also $\bfu$ be some vector. Then we write $A + \bfu$ for the translation copy of $A$ obtained when we translate $A$ by $\bfu$.

Let $A'$ and $A''$ be two geometric objects. We write $A' \sim A''$ when $A'$ and $A''$ are translation copies of one another.

Suppose that $A' \sim A''$. The \emph{displacement} from $A'$ to $A''$, denoted $A' \tto A''$, is the unique vector $\bfu$ such that $A' + \bfu = A''$.

Given a sequence of coefficients $\bfalpha = (\alpha_1, \alpha_2, \ldots, \alpha_k)$ and a sequence of vectors $\olu = (\bfu_1, \bfu_2,\allowbreak \ldots, \bfu_k)$ of the same length, we use the shorthand notation $\bfalpha \cdot \olu$ for the linear combination $\alpha_1\bfu_1 + \alpha_2\bfu_2 + \cdots + \alpha_k\bfu_k$.

We write $\mathcal{L}(\bfu, \bfv)$ for the lattice generated by the vectors $\bfu$ and $\bfv$.

Let $\bfu = (x_1, y_1)$ and $\bfv = (x_2, y_2)$ be two vectors. We write $\bfu \times \bfv$ for the oriented area of the parallelogram spanned by $\bfu$ and $\bfv$, in this order. Numerically, it equals $x_1y_2 - y_1x_2$. Notice that $\bfu$ and $\bfv$ are linearly dependant if and only if $\bfu \times \bfv = 0$. Furthermore, the counterclockwise angle $\theta$ from $\bfu$ to $\bfv$ satisfies $0^\circ < \theta < 180^\circ$ if and only if $\bfu \times \bfv > 0$.

\subsection{Polylines} \label{line}

Let $P_1$, $P_2$, $\ldots$, $P_k$ be a sequence of points. The \emph{polyline} determined by this sequence is the oriented curve which begins at $P_1$; goes along a straight-line segment to $P_2$; continues along a straight-line segment to $P_3$; and so on, until it ends at $P_k$. We denote it by $P_1 \sto P_2 \sto \cdots \sto P_k$.

Notice that it is possible for many distinct sequences of points to determine the same polyline. For example, both of $(0, 0) \sto (1, 0)$ and $(0, 0) \sto (\half, 0) \sto (1, 0)$ work out to the same oriented curve, and hence also to the same polyline.

We write $\varepsilon$ for an empty polyline; i.e., for any polyline of zero length. We also write $|S|$ for the length of the polyline $S$.

Let $S$ be the polyline $P_1 \sto P_2 \sto \cdots \sto P_k$. When we reverse the orientation of $S$, we obtain a new polyline which we call the \emph{reverse} of $S$ and which we denote by $S^{-1}$. One sequence of points which determines $S^{-1}$ would be $P_k \sto P_{k - 1} \sto \cdots \sto P_1$.

Let $S$ and $T$ be two polylines. We write $S \reveq T$ when $S = T^{-1}$; or, equivalently, when $S^{-1} = T$. We also write $S \revsim T$ when $S \sim T^{-1}$; or, equivalently, when $S^{-1} \sim T$.

We write $\lopen S \ropen$ for the oriented curve obtained from $S$ by omitting both of its endpoints. We also write $\lopen S)$ when only the initial point is omitted, and $(S \ropen$ when we omit only the final point.

We define $\sigma(S)$ to be the vector which points from the initial point of $S$ to the final point of~$S$; i.e., $\sigma(S) = P_1 \tto P_k$. We call this vector the \emph{span} of $S$. Clearly, $\sigma(S^{-1}) = -\sigma(S)$.

When $P_1 = P_k$, we say that $S$ is \emph{closed}. Or, equivalently, $S$ is closed if and only if $\sigma(S) = \mathbf{0}$.

Suppose that $S$ is closed. Then we write $\area(S)$ for the oriented area enclosed by $S$; notice that this is well-defined even if $S$ contains self-intersections. Given a point $P$ not on $S$, we also write $\wind(P, S)$ for the winding number of $S$ around $P$. Of course, $\area(S^{-1}) = -\area(S)$ and $\wind(P, S^{-1}) = -\wind(P, S)$. Furthermore, if $S$ is the boundary of a polygon, then $S$ is oriented counterclockwise if and only if $\area(S) > 0$.

Suppose that $S$ is nonempty, closed, and free of self-intersections. Then we write $[S]$ for the closed polygon determined by $S$. We also write $\langle S \rangle$ for the open polygon determined by $S$; i.e., for the interior of $[S]$. So $S$ and $\langle S \rangle$ are disjoint, with $S \cup \langle S \rangle = [S]$.

Suppose that $S'$ and $S''$ are two polylines such that the final point of $S'$ coincides with the initial point of $S''$. The \emph{concatenation} of $S'$ and $S''$, denoted $S'S''$, is the polyline obtained by joining up $S'$ and $S''$ at that point. I.e., if $S' = P_1 \sto P_2 \sto \cdots \sto P_k$ and $S'' = P_k \sto P_{k + 1} \sto \cdots \sto P_\ell$, then $S'S'' = P_1 \sto P_2 \sto \cdots \sto P_\ell$. Clearly, $\sigma(S'S'') = \sigma(S') + \sigma(S'')$.

We say that $S$ is a \emph{sub-polyline} of $T$ when there exist two polylines $S'$ and $S''$, possibly empty, such that $T = S'SS''$.

Sometimes, we use set-theoretic notation with polylines; e.g., we might say that $S' \cap S'' \subseteq T$ for some polylines $S'$, $S''$, $T$. Then we disregard orientation and we interpret each polyline as the set of its points. So ``$S' \cap S'' \subseteq T$'' would mean that every point which lies on both of $S'$ and $S''$ also lies on $T$. Notice that ``$S$ is a sub-polyline of $T$'' and ``$S \subseteq T$'' are two very different statements. For example, $TT^{-1} \subseteq T$ is true of all polylines $T$, but $TT^{-1}$ is not a sub-polyline of $T$ when $T$ is nonempty.

\subsection{Valuations} \label{val}

Let $\mathbb{P}$ be the set of all polylines in the plane. Let also $\pfloat$ be the set that we obtain when we factor $\mathbb{P}$ by translation. So, from the point of view of $\pfloat$, two polylines are ``the same'' if and only if they are translation copies of one another.

Formally, $\pfloat$ is the set of all equivalence classes of the $\sim$-relation within $\mathbb{P}$. For convenience, though, sometimes we will talk about individual polylines as if they belong to $\pfloat$; in such cases, the individual polylines will represent the equivalence classes which contain them.

We call the elements of $\pfloat$ \emph{floating polylines}. When we want to emphasise that we are viewing some polyline as an element of $\mathbb{P}$ rather than as an element of $\pfloat$, we call it a \emph{literal polyline}.

Since reversal agrees with translation, it is well-defined over $\pfloat$. So is concatenation, too. Indeed, given two floating polylines, we can always find two literal representatives of them such that the final point of the first one coincides with the initial point of the second one. The concatenation of these representatives will then be a literal representative of the concatenation of the two original floating polylines.

Let $\mathcal{A}$ be an alphabet with reversal. A \emph{valuation} over $\mathcal{A}$ is a mapping $\psi : \mathcal{A} \to \pfloat$ such that, first, $\psi(w)$ is nonempty for all letters $w$ of $\mathcal{A}$; and, second, $\psi$ respects reversal, in the sense that $\psi(w^{-1}) = \psi(w)^{-1}$ for all letters $w$ of $\mathcal{A}$.

We can naturally extend $\psi$ to the space of all words over $\mathcal{A}$. Specifically, if $W = w_1w_2 \ldots w_k$ is a word over $\mathcal{A}$, then we define $\psi(W)$ to be the concatenation of $\psi(w_1)$, $\psi(w_2)$, $\ldots$, $\psi(w_k)$, in this order.

(Notice that, in some sense, we are not giving up any expressive power by excluding $\varepsilon$ from the range of $\psi$. Indeed, suppose that $\psi_\varepsilon(W_\varepsilon) = S$ for some ``weak valuation'' $\psi_\varepsilon$ which maps some letters of $W_\varepsilon$ onto $\varepsilon$. Delete all such letters from $W_\varepsilon$ to obtain $W$, and also alter the values of $\psi_\varepsilon$ at these letters arbitrarily so as to obtain a valuation $\psi$. Then $\psi(W)$ will still work out to the same polyline $S$.)

Given a word $W$ and a polyline $S$, we say that $S$ \emph{realises} the \emph{template} $W$ when there exists some valuation $\psi$ over the alphabet of $W$ such that $\psi(W) = S$. For example, it is straightforward to check that the nonempty polylines $S$ satisfying $S = S^{-1}$ are exactly the ones which realise the template~$\mathsf{A}\mathsf{A}^{-1}$.

Clearly, every realisation of a balanced template is a closed polyline. Conversely, if every realisation of a given template is a closed polyline, then that template must be balanced.

Furthermore, if a template contains a combinatorial self-intersection, then every realisation of it will contain a geometric self-intersection. (Here it is crucial that $\varepsilon$ has been excluded from the range of $\psi$.) The converse turns out to be false; see Sections \ref{abc} and \ref{ptmp}.

A cyclic shift of a balanced template corresponds to retracing the associated closed polyline from a different initial point but with the same orientation. Notice that, if the template $W$ is balanced and $W^\shift{i}$ is any cyclic shift of it, then $W$ contains a combinatorial self-intersection if and only if $W^\shift{i}$ does. I.e., for a balanced template, whether it contains a combinatorial self-intersection or not does not depend on our choice of a ``combinatorial initial point''.

Of particular importance to us will be the ``standard'' valuation $\psi_\square$ over $\mathcal{A}_\square$ which, to $\mathsf{R}$, $\mathsf{L}$, $\mathsf{U}$, and $\mathsf{D}$, assigns unit segments pointing to the right, left, up, and down, respectively. Explicitly, \begin{align*} \psi_\square(\mathsf{R}) &= (0, 0) \sto (1, 0) & \psi_\square(\mathsf{L}) &= (0, 0) \sto (-1, 0)\\ \psi_\square(\mathsf{U}) &= (0, 0) \sto (0, 1) & \psi_\square(\mathsf{D}) &= (0, 0) \sto (0, -1). \end{align*}

We can now write down the boundary of any polyomino as a word over $\mathcal{A}_\square$.

Notice that, in the setting of polyominoes, for the alphabet $\mathcal{A}_\square$ and its standard valuation $\psi_\square$, combinatorial and geometric self-intersections are equivalent. I.e., a word $W$ over $\mathcal{A}_\square$ contains a combinatorial self-intersection if and only if the polyline $\psi_\square(W)$ associated with it contains a geometric self-intersection. Clearly, this is very much not the case for most other alphabets and valuations over them.

In a context where $\mathcal{A}$ and $\psi$ are understood, we might sometimes disregard the distinction between the words of $\mathcal{A}$ and the floating polylines which they represent. E.g., we might talk about ``$\sigma(W)$'' for some word $W$ over $\mathcal{A}$ when what we really mean is the span of its associated polyline~$\psi(W)$.

\subsection{Tilings} \label{tile}

Here, we briefly review some of the most basic properties of lattice tilings. For an in-depth discussion of tilings in general, readers are referred to \cite{GS}.

Let $S$ be a nonempty closed polyline free of self-intersections, and let $[S]$ be the polygon enclosed by it. Let also $L$ be a lattice, and consider the translation copies $[S] + \bfu$ as $\bfu$ ranges over~$L$. Roughly speaking, $[S]$ tiles the plane with $L$ when these translation copies form a partitioning of the plane.

This definition makes sense with the intuitive, geometric notion of a partitioning. However, it does not fare so well with the formal, set-theoretic notion. Then both of $[S]$ and the plane are interpreted as sets of points; and the translation copies of the point set $[S]$ do not quite form a partitioning of the set of all points in the plane.

So, for a formal definition, we must make one amendment: We say that $[S]$ \emph{tiles} the plane with $L$ when the translation copies of $[S]$ over $L$ form a partitioning of the plane; except that the boundaries of the copies are allowed to overlap. Or, in other words, for each point $P$ in the plane, either $P$ lies in the interior of exactly one copy; or else it lies simultaneously on the boundaries of two or more copies but not in any of the interiors.

Suppose that $[S]$ does tile the plane with $L$, and let $\mathfrak{T}$ be the resulting tiling. We call the individual copies of $[S]$ in $\mathfrak{T}$ \emph{tiles}. We say that two tiles are \emph{neighbours} when the intersection of their boundaries contains a nonempty polyline.

It is not too difficult to check that, in a lattice tiling, the common boundary of two neighbours is always contiguous. Here is a quick informal sketch of one argument: Suppose, for the sake of contradiction, that the common boundary of $[S]$ and $[S] + \bfu$ is not contiguous. Then the union of these two tiles must ``envelop'' some other tile. Consider the projections of all three tiles onto some line perpendicular to $\bfu$. The projections of $[S]$ and $[S] + \bfu$ will coincide; whereas the projection of the enveloped tile will be congruent to these two. From here, it is fairly straightforward to reach a contradiction.

Notice that, in \cite{GS}, this property is assumed by fiat as a part of the definition of a ``normal'' tiling. Indeed, for non-lattice tilings it is not necessarily true; see, e.g., Figure 3.2.4 of \cite{GS}.

We pause now to address one minor technical point regarding the tilings of polyominoes.

It is tempting to say that, in a polyomino tiling, the common boundary of two neighbours is always expressible over $\mathcal{A}_\square$. But, if we view polyominoes purely as geometric objects, then this becomes false. Consider, for example, any tiling of the unit square which uses a lattice of the form $\mathcal{L}((1, 0), (\alpha, 1))$ with $\alpha \not \in \mathbb{Z}$.

Thus we introduce the following convention: Whenever we are working with polyomino tilings in a combinatorial context, we assume all tiles to ``inhabit'' one and the same implicit unit-square grid. Or, in other words, implicitly we fix a grid first, and then we assume all tiles to be composed out of cells within it. The tiling itself we view as a partitioning of the set of all cells of the grid. The common boundary of each pair of neighbours will then lie along the grid lines, and so it will indeed be expressible over $\mathcal{A}_\square$.

With this, we return to lattice tilings in general.

There is a natural way to interpret $\mathfrak{T}$ as a planar graph drawing $\mathfrak{G}$, as follows: First, the vertices of $\mathfrak{G}$ are those points where three or more tiles meet. Second, the edges of $\mathfrak{G}$ are the common boundaries of all pairs of neighbours. The tiles of $\mathfrak{T}$ will then correspond to the faces~of~$\mathfrak{G}$.

Two tilings $\mathfrak{T}'$ and $\mathfrak{T}''$ are \emph{combinatorially isomorphic} when there exists an isomorphism between the corresponding planar graph drawings $\mathfrak{G}'$ and $\mathfrak{G}''$. Notice that, since the isomorphism is between the planar drawings rather than between the underlying graphs, it must biject not only vertices with vertices and edges with edges, but also faces with faces.

It is well-known that every lattice tiling is combinatorially isomorphic either to the grid of unit squares or to the grid of unit regular hexagons. See, e.g., Section 6 of \cite{GS}. The treatment there is in fact much more general, covering not just lattice tilings but also arbitrary ``isohedral'' tilings; and, additionally, it takes into account a lot more information about the tilings than just their combinatorial isomorphism classes. The results are summarised in Table 6.2.1 of \cite{GS}; but it is only rows $1$ and $41$ of this table that are relevant to our purposes.

Let $D$ be the set of all vectors which map a tile to one of its neighbours. Then $D$ must be centrally symmetric. Indeed, suppose that $\bfu \in D$ and $[S]$ is neighbours with $[S] + \bfu$; then also $[S] + \bfu$ is neighbours with $[S]$, and so $-\bfu \in D$ as well.

For a lattice tiling in the square-grid combinatorial-isomorphism class, each tile has $4$ neighbours and $|D| = 4$. Similarly, for a lattice tiling in the hexagon-grid combinatorial-isomorphism class, each tile has $6$ neighbours and $|D| = 6$.

The more detailed treatment in \cite{GS} shows that, furthermore, if the neighbours of $[S]$ occur around it in the cyclic order $[S] + \bfu_1$, $[S] + \bfu_2$, $\ldots$, $[S] + \bfu_{2k}$, with $k \in \{2, 3\}$, then $\bfu_i + \bfu_{i + k} = \mathbf{0}$ for all $1 \le i \le k$. I.e., if the displacements from a tile to two of its neighbours are negations of one another, then these two neighbours must occur at ``opposite'' positions.

From now on, we will focus mostly on $4$-neighbour lattice tilings -- except for a brief excursion into $6$-neighbour lattice tilings in Section \ref{six}.

\subsection{The Template \texorpdfstring{$\tmptile$}{ThetaTile}} \label{abab}

In the setting of the previous sub-section, suppose that $\mathfrak{T}$ is a $4$-neighbour lattice tiling with $D = \{-\bfu, \bfu, -\bfv, \bfv\}$.

Consider the common boundary of $[S]$ with its neighbour $[S] - \bfv$. Since it is contiguous, it must be a sub-polyline of $S$. Denote this sub-polyline by $A_+$. Define also $A_-$, $B_+$, $B_-$ similarly, for the neighbours $[S] + \bfv$, $[S] + \bfu$, $[S] - \bfu$ of $[S]$, respectively. Then, without loss of generality, $S = A_+B_+A_-B_-$.

Translating the ordered pair $([S], [S] - \bfv)$ by $\bfv$, we obtain the ordered pair $([S] + \bfv, [S])$. Since the common boundary of the two polygons in the first pair is $A_+$, and for the second pair it is $A_-$, it follows that $A_- \revsim A_+$. Similarly, $B_- \revsim B_+$ as well.

\begin{figure}[ht] \centering \includegraphics{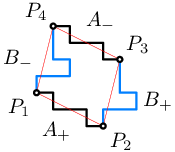} \caption{} \label{ababf} \end{figure} 

For example, Figure \ref{ababf} shows this for the tile of Figure \ref{def-a} and the tiling of Figure \ref{def-b}. (The parallelogram outlined in red will illustrate the proof of Lemma \ref{area} below.)

Recall the template terminology of Section \ref{val}. Since $A_- \revsim A_+$ and $B_- \revsim B_+$, we conclude that $S$ realises the template \[\tmptile = \mathsf{A}\mathsf{B}\mathsf{A}^{-1}\mathsf{B}^{-1}.\]

Or, to summarise, if $[S]$ admits a $4$-neighbour lattice tiling of the plane, then its boundary must realise a certain template. We proceed now to prove the converse as well. First, though, we establish one auxiliary result. It will be of independent usefulness to us also in Section \ref{orient}.

\begin{lemma} \label{area} Suppose that the polyline $S$ realises the template $\tmptile$ with $\bfu = \sigma(\mathsf{A})$ and $\bfv = \sigma(\mathsf{B})$. Then $\area(S) = \bfu \times \bfv$. \end{lemma}

\begin{myproof} Let the letters $\mathsf{A}$ and $\mathsf{B}$ of $\tmptile$ evaluate to $A$ and $B$. Let also $S = A_+B_+A_-B_-$ as above; so the literal polylines $A_+$, $B_+$, $A_-$, $B_-$ are representatives of the floating polylines $A$, $B$, $A^{-1}$, $B^{-1}$, respectively.

Denote the initial points of $A_+$, $B_+$, $A_-$, $B_-$ by $P_1$, $P_2$, $P_3$, $P_4$. Then $P_1 \tto P_2 = P_4 \tto P_3 = \bfu$ and $P_1 \tto P_4 = P_2 \tto P_3 = \bfv$; hence, $P_1$, $P_2$, $P_3$, $P_4$ are the vertices of a parallelogram. Let $T = P_1 \sto P_2 \sto P_3 \sto P_4 \sto P_1$ be this parallelogram's boundary. (See Figure \ref{ababf}.) Let also $A^\text{Close}_+$ be the closed polyline obtained from $A_+$ by joining its final point back to its initial point, using the straight-line segment $P_2 \sto P_1$; and define $B^\text{Close}_+$, $A^\text{Close}_-$, $B^\text{Close}_-$ similarly.

Then \[\area(S) = \area(T) + \area(A^\text{Close}_-) + \area(B^\text{Close}_-) + \area(A^\text{Close}_+) + \area(B^\text{Close}_+).\]

However, $\area(A^\text{Close}_-) + \area(A^\text{Close}_+) = 0$ because, modulo a suitable choice of initial points, $A^\text{Close}_- \revsim A^\text{Close}_+$; and similarly also for $B^\text{Close}_-$ and $B^\text{Close}_+$. Finally, $\area(T) = \bfu \times \bfv$. \end{myproof}

We continue on to the advertised converse.

\begin{lemma} \label{abl} A given polygon admits a $4$-neighbour lattice tiling of the plane if and only if its boundary realises the template $\tmptile$. (Up to a suitable choice of initial point.) \end{lemma} 

We have stated Lemma \ref{abl} as a necessary and sufficient condition. However, the ``only if'' direction has already been dealt with, and so in the proof we will consider the ``if'' direction only.

Notice that, since $\tmptile$ is balanced, every polyline which realises it must be closed. Hence, if a polyline is free of self-intersections and it realises $\tmptile$, then it must be the boundary of a polygon.

Let, in the setting of Lemma \ref{abl}, once again $\bfu = \sigma(\mathsf{A})$ and $\bfv = \sigma(\mathsf{B})$. From the proof, it will become clear that (as expected) we get a tiling with the lattice $\mathcal{L}(\bfu, \bfv)$; for this tiling, $D = \{-\bfu, \bfu, -\bfv, \bfv\}$; and the common boundaries of each tile with its four neighbours are the polylines associated with $\mathsf{A}$, $\mathsf{B}$, $\mathsf{A}^{-1}$, $\mathsf{B}^{-1}$. Whenever we invoke Lemma \ref{abl}, we will assume all of these ``supplements'' to be implied as well.

\begin{myproof} Let $S$ be a closed polyline free of self-intersections which realises $\tmptile$. Define also $A$ and $B$ as in the proof of Lemma \ref{area}.

Since $S$ is free of self-intersections, both of $\bfu$ and $\bfv$ are nonzero. So we can assume, modulo a suitable affine transformation, that $\bfu = (1, 0)$ and $\bfv = (0, 1)$. By Lemma \ref{area}, it follows that $\area(S) = 1$ and $S$ is oriented counterclockwise.

For all $(x, y) \in \mathbb{Z}^2$, let $A_{x, y}$ be an instance of $A$ with initial point $(x, y)$ and final point $(x + 1, y)$. Similarly, let $B_{x, y}$ be an instance of $B$ with initial point $(x, y)$ and final point $(x, y + 1)$. Let also $S_{x, y}$ be the closed polyline $A_{x, y}B_{x + 1, y}A_{x, y + 1}^{-1}B_{x, y}^{-1}$. Then $S_{x, y}$ is a translation copy of $S$.

We claim that the polygons enclosed by these translation copies form a tiling of the plane. Hence, our goal now becomes to show that the interiors $\langle S_{x, y} \rangle$ are pairwise disjoint as $(x, y)$ ranges over $\mathbb{Z}^2$.

Let $P$ be any point in the plane which does not lie on one of our translation copies of $S$. Since $S$ is oriented counterclockwise, the winding number $\wind(P, S_{x, y})$ will equal $1$ when $P \in \langle S_{x, y} \rangle$, and $0$ otherwise. So the total number of interiors that $P$ belongs to must be \[\sum_{(x, y) \in \mathbb{Z}^2} \wind(P, S_{x, y});\] notice that this sum is well-defined because only finitely many of its summands are nonzero.

For all $(x, y) \in \mathbb{Z}^2$, let $E_{x, y}$ be the closed polyline $(x, y) \sto (x + 1, y) \sto (x + 1, y + 1) \sto (x, y + 1) \sto (x, y)$. Let also $A^\text{Close}_{x, y}$ be the closed polyline obtained from $A_{x, y}$ by joining its final point back to its initial point, using the straight-line segment $(x + 1, y) \sto (x, y)$; and define $B^\text{Close}_{x, y}$ similarly.

For convenience, let us assume that $P$ does not lie on any polyline of the form $E_{x, y}$, either; or, equivalently, that both coordinates of $P$ are non-integers. (This restriction will turn out not to matter too much for our argument.) For all such $P$, the winding number $\wind(P, S_{x, y})$ will equal \[\wind(P, E_{x, y}) + \wind(P, A^\text{Close}_{x, y}) + \wind(P, B^\text{Close}_{x + 1, y}) - \wind(P, A^\text{Close}_{x, y + 1}) - \wind(P, B^\text{Close}_{x, y}).\]

Once again, for each term in this expression, there are only finitely many integer points $(x, y)$ where it is nonzero. Hence, we can safely take the sum over all $(x, y) \in \mathbb{Z}^2$. Telescoping, we arrive~at \[\sum_{(x, y) \in \mathbb{Z}^2} \wind(P, S_{x, y}) = \sum_{(x, y) \in \mathbb{Z}^2} \wind(P, E_{x, y}) = 1;\] where the final step follows because the unit squares $[E_{x, y}]$ form a tiling of the plane when $(x, y)$ ranges over $\mathbb{Z}^2$, and so $P$ must lie in the interior of exactly one of them.

We conclude that the interiors of our translation copies of $[S]$ can only intersect at points $P$ such that $P$ lies on some polyline of the form $E_{x, y}$ or $S_{x, y}$. Since each one of these interiors is an open subset of the plane, it follows that in fact they cannot intersect at all. \end{myproof}

\subsection{Deformations} \label{def}

Let $\Phi$ be a polyomino. Then we can interpret the boundary of $\Phi$ as a polyline $F$ described by some word $W$ over the alphabet $\mathcal{A}_\square$.

Fix some polyline $S$ which is free of self-intersections and which realises the template $\tmptile$. Define also $A$ and $B$ as in the previous sub-section.

Consider now a new, non-standard valuation $\psi_\diamond$ for $\mathcal{A}_\square$, given by \begin{align*} \psi_\diamond(\mathsf{R}) &= A & \psi_\diamond(\mathsf{U}) &= B\\ \psi_\diamond(\mathsf{L}) &= A^{-1} & \psi_\diamond(\mathsf{D}) &= B^{-1}. \end{align*}

Denote $F_\diamond = \psi_\diamond(W)$. Then we call $F_\diamond$ a \emph{deformation} of $F$.

\begin{lemma} \label{defl} In the setting described above, $F_\diamond$ is closed and free of self-intersections. \end{lemma} 

\begin{myproof} It is clear that $F_\diamond$ is closed, and so we turn to showing that it is free of self-intersections.

Let $\mathfrak{G}_\square$ be the planar graph drawing whose vertices are all integer points and whose edges are all unit segments between these points. Since $\Phi$ is a polyomino, we can assume without loss of generality that $F$ is a simple cycle in $\mathfrak{G}_\square$.

Let $\mathfrak{T}$ be the tiling of the plane with translation copies of $[S]$ constructed in the proof of Lemma \ref{abl}. Let also $\mathfrak{G}$ be the planar graph drawing associated with $\mathfrak{T}$, as in Section \ref{tile}.

Clearly, the planar graph drawings $\mathfrak{G}_\square$ and $\mathfrak{G}$ are isomorphic. Since $F$ is a simple cycle in $\mathfrak{G}_\square$, it follows that also $F_\diamond$ is a simple cycle in $\mathfrak{G}$. (Or, strictly speaking, it follows that some literal-polyline simple cycle in $\mathfrak{G}$ is a representative of the floating polyline $F_\diamond$.) This guarantees that $F_\diamond$ is free of self-intersections, as desired. \end{myproof}

Since $F_\diamond$ is nonempty, closed, and free of self-intersections, it is the boundary of a polygon. Let $\Phi_\diamond = [F_\diamond]$; we call this polygon a \emph{deformation} of $\Phi$.

From the proof of Lemma \ref{defl}, we get that $\Phi_\diamond$ is composed of several tiles of $\mathfrak{T}$ joined together in the exact same way as the cells of $\Phi$. Or, in other words, we can obtain $\Phi_\diamond$ from $\Phi$ by replacing each cell of $\Phi$ with a translation copy of the building block $[S]$.

Deformations will be important for the statement of our main result in Section \ref{smt}.

When both of $A$ and $B$ are expressible over $\mathcal{A}_\square$ (with its standard valuation $\psi_\square$), the deformation $\Phi_\diamond$ will be a polyomino as well. In this special case, we can also think of deformations in terms of substitutions. Indeed, suppose that $A = \psi_\square(W_\mathsf{A})$ and $B = \psi_\square(W_\mathsf{B})$. Then $F_\diamond$ will be described by the word $W[\mathsf{R} \to W_\mathsf{A}, \mathsf{U} \to W_\mathsf{B}]$ over $\mathcal{A}_\square$. We return to this observation in Section~\ref{root}.

\subsection{The Template \texorpdfstring{$\tmpcube$}{ThetaCube}} \label{abc}

Suppose that we are given a template $W$, and we want to determine whether $W$ can be realised without self-intersections. When $W$ contains a combinatorial self-intersection, the answer is trivially no, as we observed in Section \ref{val}. However, there are also some templates which fail for subtler reasons.

One such template will be particularly important to us later on. Let \[\tmpcube = \mathsf{A}\mathsf{B}\mathsf{C}\mathsf{B}^{-1}\mathsf{A}^{-1}\mathsf{B}\mathsf{C}^{-1}\mathsf{B}^{-1}.\]

(The choice of subscript is because, if we realise this template in three dimensions by evaluating $\mathsf{A}$, $\mathsf{B}$, and $\mathsf{C}$ with three mutually perpendicular unit segments, we get a Hamiltonian cycle of the unit cube. Since the cube admits an essentially unique Hamiltonian cycle, this observation can also serve as a mnemonic for the template.)

\begin{figure}[ht] \centering \includegraphics{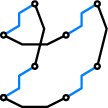} \caption{} \label{cube} \end{figure}

One example is shown in Figure \ref{cube}. The parts which correspond to the letter $\mathsf{B}$ and its reverse are coloured in blue; and the other parts of the realisation are all coloured in black.

\begin{lemma} \label{abcl} Suppose that the polyline $S$ realises the template $\tmpcube$. Then $S$ must contain a self-intersection. \end{lemma} 

\begin{myproof} Suppose, for the sake of contradiction, that $S$ is free of self-intersections.

Since \[\tmpcube = \tmptile[\mathsf{A} \to \mathsf{A}, \mathsf{B} \to \mathsf{B}\mathsf{C}\mathsf{B}^{-1}],\] it follows by Lemma \ref{abl} that $[S]$ tiles the plane using the lattice generated by $\bfu' = \sigma(\mathsf{A})$ and $\bfv' = \sigma(\mathsf{B}\mathsf{C}\mathsf{B}^{-1})$.

Similarly, since \[\tmpcube^\shift{7} = \tmptile[\mathsf{A} \to \mathsf{B}^{-1}\mathsf{A}\mathsf{B}, \mathsf{B} \to \mathsf{C}],\] it follows by Lemma \ref{abl} once again that $[S]$ also tiles the plane using the lattice generated by $\bfu'' = \sigma(\mathsf{B}^{-1}\mathsf{A}\mathsf{B})$ and $\bfv'' = \sigma(\mathsf{C})$.

But $\bfu' = \sigma(\mathsf{A}) = \sigma(\mathsf{B}^{-1}\mathsf{A}\mathsf{B}) = \bfu''$, and similarly $\bfv' = \bfv''$. We conclude that our two tilings are in fact one and the same. However, in the first tiling, the common boundary of $[S]$ with its neighbour $[S] - \bfv'$ is described by $\mathsf{A}$; whereas, in the second tiling, the same common boundary is described by $\mathsf{B}^{-1}\mathsf{A}\mathsf{B}$. Since the former is a strict subword of the latter, these two words cannot evaluate to the same floating polyline, and we arrive at a contradiction. \end{myproof}

Notice that all templates obtained from $\tmpcube$ by means of one or both of the substitutions $[\mathsf{A} \to \varepsilon]$ and $[\mathsf{C} \to \varepsilon]$ contain combinatorial self-intersections. Hence, the result will continue to hold even if we allow ``weak valuations'' where one or both of $\mathsf{A}$ and $\mathsf{C}$ map onto $\varepsilon$; but we must preserve the requirement that $\mathsf{B}$ evaluates to a nonempty polyline. We call this the ``strong form'' of Lemma \ref{abcl}. It is in fact this slight strengthening, rather than the original form, that we are going to use in Sections \ref{ni} and \ref{loop}.

We pick up the topic of ``paradoxical templates'' once again in Section \ref{ptmp}.

\section{Double Chains} \label{chain}

Here, we study the structure of $4$-neighbour double tiles.

Section \ref{bc} introduces the notion of a chain -- essentially, a closed polyline subdivided into parts. Then, in Section \ref{dcp}, we derive certain systems of equations over polylines which characterise the chains associated with $4$-neighbour double tiles. Finally, in Section \ref{dcw}, we translate these systems of equations from the language of polylines into the language of words.

\subsection{Chains} \label{bc} 

A \emph{cyclic sequence} is a finite sequence ``with wraparound''. Or, more formally, when we refer to some finite sequence $U_1$, $U_2$, $\ldots$, $U_k$ as cyclic, what we mean is that the expression $U_{k + 1}$ has been implicitly defined to mean $U_1$; the expression $U_0$ has been implicitly defined to mean $U_k$; and, in general, for all integers $i$, the expression $U_i$ has been implicitly defined to mean~$U_{i \bmod k}$.

A \emph{chain} is a cyclic sequence of polylines such that the final point of each polyline coincides with the initial point of the next one. We write $U_1 : U_2 : \cdots : U_k$ for the chain $U$ formed out of the polylines $U_1$, $U_2$, $\ldots$, $U_k$, in this order. So, for all $i$, the final point of $U_i$ coincides with the initial point of $U_{i + 1}$. Each one of $U_1$, $U_2$, $\ldots$, $U_k$ is called a \emph{part} of the chain.

Let $P_i$ be the initial point of $U_i$. Then the final point of $U_i$ will be $P_{i + 1}$. We call $P_1$, $P_2$, $\ldots$, $P_k$ the \emph{division points} of the chain.

Clearly, we can concatenate the parts of $U$ so as to obtain the closed polyline $U_1U_2 \cdots U_k$. Conversely, if we begin with a closed polyline and we subdivide it into parts, we get a chain. We write $|U|$ for the length of $U$; i.e., for $|U_1| + |U_2| + \cdots + |U_k|$. We also write $[U]$ as an abbreviation for $[U_1U_2 \cdots U_k]$.

The \emph{reverse} of $U$, denoted $U^{-1}$, is the chain $U_k^{-1} : U_{k - 1}^{-1} : \cdots : U_1^{-1}$. The $i$-th \emph{cyclic shift} of $U$, denoted $U^\shift{i}$, is the chain $U_{i + 1} : U_{i + 2} : \cdots : U_k : U_1 : U_2 : \cdots : U_i$.

Let $\iota$ be an isometry in the plane; e.g., a reflection or a rotation. We say that $\iota$ is a \emph{symmetry} of $U$ when it preserves $U$ modulo reversals and cyclic shifts. For example, if $\iota$ is a $90^\circ$ rotation around the origin which maps $U_i$ onto $U_{i + 2}$ for all $i$, then it will be a symmetry of $U$.

A \emph{self-intersection} in $U$ is a self-intersection in the concatenation $U_1U_2 \cdots U_k$ of $U$. Correspondingly, we say that a chain is free of self-intersections when its concatenation is free of self-intersections.

Suppose that $i \le j$. Then we write $U_{i : j}$ for the concatenation $U_iU_{i + 1} \cdots U_j$. We call each such concatenation a \emph{slice} of $U$. Notice that, by virtue of $U$ being a cyclic sequence, the expression $U_{i : j}$ continues to make sense even when $i$ and $j$ are not elements of the interval $\{1, 2, \ldots, k\}$. For example, if $k = 3$, then $U_{2 : 7} = U_2U_3U_1U_2U_3U_1$.

For convenience, we allow the use of the notation $U_{i : j}$ also when $1 \le j < i \le k$. (But not for any other pairs of indices with $i > j$.) We understand this to mean $U_iU_{i + 1} \cdots U_kU_1U_2 \cdots U_j$; or, equivalently, $U_{i : (j + k)}$. The intuition is that we begin at $U_i$, we wrap around at the end, and we stop at $U_j$. For example, if $k = 5$, then $U_{4 : 2} = U_4U_5U_1U_2$.

We say that a finite sequence of floating polylines is a \emph{chain} when some representatives of its elements form a chain of literal polylines. Or, equivalently, when the spans of its elements sum to $\mathbf{0}$.

Our main tool for the study of chains will be the encoding of their parts into words. So now we define a cyclic sequence of words over the alphabet with reversal $\mathcal{A}$ to be a \emph{chain} when its concatenation is balanced.

Let $W = W_1 : W_2 : \cdots : W_k$ be one such word chain. The definitions of $|W|$, $W^{-1}$, and $W^\shift{i}$ carry over from the polyline setting without change. So does the shorthand notation $W_{i : j}$. (Though notice that the $i$-th cyclic shift of $W$ is very different from the $i$-th cyclic shift of the concatenation $W_1W_2 \cdots W_k$ of $W$. The former is a word chain; the latter is a word; and the concatenation of $W^\shift{i}$ will differ from $(W_1W_2 \cdots W_k)^\shift{i}$ in general.)

A \emph{subword} of $W$ is any subword of the concatenation $W_1W_2 \cdots W_k$ of $W$. A \emph{combinatorial self-intersection} in $W$ is defined similarly. We apply substitutions to word chains part-wise; i.e., if $\zeta$ is a substitution, then we define $\zeta(W) = \zeta(W_1) : \zeta(W_2) : \cdots : \zeta(W_k)$.

Let $\psi$ be a valuation over $\mathcal{A}$. Then the $\psi$-images of the parts of $W$ will form a chain of floating polylines. Hence, we define $\psi(W) = \psi(W_1) : \psi(W_2) : \cdots : \psi(W_k)$. We may also say that a given word chain ``evaluates'' to a given chain of literal polylines when, strictly speaking, it evaluates to the corresponding chain of floating polylines.

We can now interpret word chains as templates, and the definition of a chain of a polylines realising such a template is fully analogous to the definition of a word template being realised by a polyline.

\subsection{Double Chains of Polylines} \label{dcp}

Let $\Phi$ be a $4$-neighbour double tile, and fix two distinct $4$-neighbour lattice tilings $\mathfrak{T}'$ and $\mathfrak{T}''$ of~$\Phi$. By Lemma \ref{abl}, from $\mathfrak{T}'$ we get four division points $P_1$, $P_2$, $P_3$, $P_4$ on the boundary of $\Phi$ which partition it into a polyline chain realising the chain template $\mathsf{A} : \mathsf{B} : \mathsf{A}^{-1} : \mathsf{B}^{-1}$. Define also the division points $Q_1$, $Q_2$, $Q_3$, $Q_4$ similarly for $\mathfrak{T}''$.

Notice that $P_1$ and $P_3$ are ``arc-length opposites'' on the boundary of $\Phi$, in the sense that they partition it into two polylines of the same length. The same observation is true also of the pairs $\{P_2, P_4\}$, $\{Q_1, Q_3\}$, and $\{Q_2, Q_4\}$. Hence, we can assume without loss of generality that these eight points occur on the boundary of $\Phi$ either in the cyclic order \[P_1, Q_1, P_2, Q_2, P_3, Q_3, P_4, Q_4 \tag{\textbf{i}}\] or in the cyclic order \[P_1, P_2, Q_1, Q_2, P_3, P_4, Q_3, Q_4\makebox[0pt][l]{.} \tag{\textbf{ni}}\]

We refer to the first one of these as the \emph{interleaved case}, and to the second one as the \emph{non-interleaved case}.

Let $F = F_1 : F_2 : \cdots : F_8$ be the polyline chain that we obtain when we partition the boundary of $\Phi$ by means of the ordered octuple of division points (\textbf{i}) in the interleaved case, and by means of (\textbf{ni}) in the non-interleaved case.

Consider the interleaved case first. Then $F_1$ begins at $P_1$ and ends at $Q_1$; $F_2$ begins at $Q_1$ and ends at $P_2$; $F_3$ begins at $P_2$ and ends at $Q_2$; and so on. By the definition of the division quadruple $P_1$, $P_2$, $P_3$, $P_4$, it follows that the chain $F_1F_2 : F_3F_4 : F_5F_6 : F_7F_8$ realises the chain template $\mathsf{A} : \mathsf{B} : \mathsf{A}^{-1} : \mathsf{B}^{-1}$. Similarly, the definition of the division quadruple $Q_1$, $Q_2$, $Q_3$, $Q_4$ tells us that the chain $F_2F_3 : F_4F_5 : F_6F_7 : F_8F_1$ realises the same chain template as well.

We conclude that, in the interleaved case, $F$ must satisfy the constraints \begin{align*} F_1F_2 &\sim F_6^{-1}F_5^{-1}\\ F_2F_3 &\sim F_7^{-1}F_6^{-1}\\ F_3F_4 &\sim F_8^{-1}F_7^{-1}\\ F_4F_5 &\sim F_1^{-1}F_8^{-1}. \end{align*}

Or, in the shorthand notation of the previous sub-section, \[F_{i : (i + 1)} \revsim F_{(i + 4) : (i + 5)} \qquad \text{for all $i$}.\]

\begin{figure}[ht] \centering \includegraphics{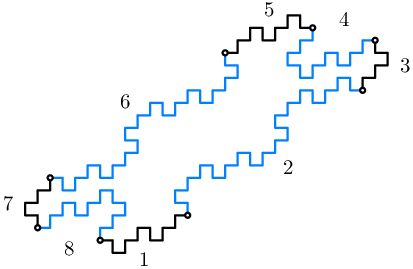} \caption{} \label{cf} \end{figure} 

We define an \emph{interleaved double chain} to be any eight-part chain which satisfies these constraints. For example, the chain which corresponds to the polyomino in Figure \ref{t} and its two tilings in Figure \ref{tilings} is of this form; it is shown in Figure \ref{cf}.

For the non-interleaved case, a similar analysis yields the constraints \begin{align*} F_1 &\sim F_5^{-1}\\ F_2F_3F_4 &\sim F_8^{-1}F_7^{-1}F_6^{-1}\\ F_3 &\sim F_7^{-1}\\ F_4F_5F_6 &\sim F_2^{-1}F_1^{-1}F_8^{-1}. \end{align*}

Or, using the shorthand notation of the previous sub-section, \[\begin{aligned} F_i &\revsim F_{i + 4}\\ F_{(i - 1) : (i + 1)} &\revsim F_{(i + 3) : (i + 5)} \end{aligned} \qquad \text{for all odd $i$}.\]

Just as in the interleaved case, we define a \emph{non-interleaved double chain} to be any eight-part chain which satisfies these constraints. Finally, we also define a \emph{double chain} to be either an interleaved double chain or a non-interleaved double chain.

\subsection{Double Chains of Words} \label{dcw}

Suppose, in the notation of the previous sub-section, that $\Phi$ is a polyomino. Suppose also that $\mathfrak{T}'$ and $\mathfrak{T}''$ respect the convention of Section \ref{tile} regarding polyomino tilings considered in a combinatorial context. Then the polyline chain $F$ will be described by some word chain over~$\mathcal{A}_\square$. Furthermore, the parts of this word chain will satisfy the exact same constraints as the parts of $F$, except that each instance of ``$\sim$'' will be replaced with an ``$=$''. We arrive at the following definitions:

Let $\mathcal{A}$ be an alphabet with reversal, and let $U = U_1 : U_2 : \cdots : U_8$ be a word chain over $\mathcal{A}$. We say that $U$ is an \emph{interleaved double chain} when it satisfies the system of equations \[\begin{aligned} U_1U_2 &= U_6^{-1}U_5^{-1}\\ U_2U_3 &= U_7^{-1}U_6^{-1}\\ U_3U_4 &= U_8^{-1}U_7^{-1}\\ U_4U_5 &= U_1^{-1}U_8^{-1}. \end{aligned} \tag{\textbf{I}}\]

Or, in the shorthand notation of Section \ref{bc}, \[U_{i : (i + 1)} \reveq U_{(i + 4) : (i + 5)} \qquad \text{for all $i$}.\]

We also say that $U$ is a \emph{non-interleaved double chain} when it satisfies the system of equations \[\begin{aligned} U_1 &= U_5^{-1}\\ U_2U_3U_4 &= U_8^{-1}U_7^{-1}U_6^{-1}\\ U_3 &= U_7^{-1}\\ U_4U_5U_6 &= U_2^{-1}U_1^{-1}U_8^{-1}. \end{aligned} \tag{\textbf{NI}}\]

Or, in concise form, \[\begin{aligned} U_i &\reveq U_{i + 4}\\ U_{(i - 1) : (i + 1)} &\reveq U_{(i + 3) : (i + 5)} \end{aligned} \qquad \text{for all odd $i$}.\]

Finally, we say that $U$ is a \emph{double chain} when it is either an interleaved double chain or a non-interleaved double chain.

The system (\textbf{I}) implies that $|U_i| + |U_{i + 1}| = |U_{i + 4}| + |U_{i + 5}|$ for all $i$. From these equations, it follows easily that also $|U_i| = |U_{i + 4}|$ for all $i$. By a similar argument, the same conclusion holds for the system (\textbf{NI}) as well. Let, then, \[a_i = |U_i| = |U_{i + 4}| \qquad \text{for all $i$}.\]

We define the \emph{type} of a double chain $U$ to be the ordered quadruple of nonnegative integers $(a_1, a_2, a_3, a_4)$. Clearly, $|U| = 2(a_1 + a_2 + a_3 + a_4)$ both in the interleaved and in the non-interleaved case.

For polyominoes, there is a kind of equivalence between solutions to the $4$-neighbour double tile problem, on the one hand, and double chains over $\mathcal{A}_\square$, on the other. Explicitly:

\begin{lemma} \label{dc} A word over $\mathcal{A}_\square$ describes the boundary of a $4$-neighbour double tile if and only if it is free of combinatorial self-intersections and, up to a suitable cyclic shift, it can be partitioned into a double chain $U$ such that the type $(a_1, a_2, a_3, a_4)$ of $U$ satisfies $a_1 + a_3 \neq 0$ and $a_2 + a_4 \neq 0$. \end{lemma} 

\begin{myproof} By Lemma \ref{abl} and the analysis of Section \ref{dcp}. The conditions on the type parameters of $U$ ensure that the two tilings we get out of Lemma \ref{abl} are in fact distinct. \end{myproof}

Thus the polyomino case of our main problem reduces to the study of those double chains over $\mathcal{A}_\square$ which are free of combinatorial self-intersections.

This is exactly what we will be doing over the course of Sections \ref{lift} and \ref{neg}. Though we will be working in somewhat greater generality; we will consider double chains over arbitrary alphabets with reversal, not just $\mathcal{A}_\square$. The reason for this is that we will in fact require the greater generality in Section \ref{poly}, where we are going to extend our results from polyominoes to arbitrary polygons.

For Section \ref{lift}, we will focus on the interleaved double chains. (Eventually, we will see that they are the only ones which actually yield polyominoes. Conversely, the non-interleaved double chains will turn out to be impossible to realise without geometric self-intersections. This is far from obvious; we prove it in Section \ref{ni}.) We will develop a tiny ``theory'' of them, showing that they exhibit a certain kind of ``fractal-like'' structure. Then, in Section \ref{neg}, we will apply this theory to the proof of our main theorem.

\section{Reductions and Lifts} \label{lift}

In this section, we introduce two families of transformations over word chains which will be crucial to our understanding of $4$-neighbour double tiles.

The first family, of the $f$-transforms, is introduced in Section \ref{ft}. Section \ref{iw} is a brief interlude which expands our set of tools for manipulating words. Thus equipped, in Section \ref{gt} we introduce also the second family, of the $g$-transforms. Then, in Section \ref{rl}, we apply these transformations to the study of the word chains which describe $4$-neighbour double tiles. Finally, in Section \ref{smt}, we state and discuss our main theorem in the special case of polyominoes.

\subsection{The \texorpdfstring{$f$}{f}-Transforms} \label{ft}

Let $V = V_1 : V_2 : \cdots : V_8$ be an interleaved double chain of type $(b_1, b_2, b_3, b_4)$ over the alphabet with reversal $\mathcal{A}$.

Suppose that there exists some index $i$ with $b_i \ge b_{i - 1} + b_{i + 1}$. (Here, we are viewing the type of $V$ as a cyclic sequence.) For the sake of clarity, below we will focus on the case when $i = 1$; the other three cases are analogous. So now $b_1 \ge b_2 + b_4$.

Since $b_1 \ge b_2$, the equation $V_1V_2 = V_6^{-1}V_5^{-1}$ of (\textbf{I}) tells us that $V_6^{-1}$ is a prefix of $V_1$. Similarly, the equation $V_8V_1 = V_5^{-1}V_4^{-1}$ tells us that $V_4^{-1}$ is a suffix of $V_1$. Furthermore, since $b_1 \ge b_2 + b_4$, these two subwords of $V_1$ are disjoint.

It follows that $V_1 = V_6^{-1}W_1V_4^{-1}$ for some word $W_1$ over $\mathcal{A}$. Similarly, also $V_5 = V_2^{-1}W_5V_8^{-1}$ for some word $W_5$ over $\mathcal{A}$.

Let us substitute the right-hand sides of these identities into the two equations of (\textbf{I}) which involve $V_1$ and $V_5$. Simplifying, we arrive at \begin{align*} W_1V_4^{-1} &= V_8W_5^{-1}\\ V_6^{-1}W_1 &= W_5^{-1}V_2. \end{align*}

Therefore, the chain \[W_5^{-1} : V_2 : V_3 : V_4 : W_1^{-1} : V_6 : V_7 : V_8\] satisfies the equations of (\textbf{I}), and so it is an interleaved double chain over $\mathcal{A}$.

Denote this chain by $U = U_1 : U_2 : \cdots : U_8$. For the type $(a_1, a_2, a_3, a_4)$ of $U$, we get that $a_1 = b_1 - (b_2 + b_4)$ and $a_i = b_i$ when $i \in \{2, 3, 4\}$.

We call $U$ the $f_1$-\emph{reduction} of $V$, and we write $U = f_1^{-1}(V)$. We also say that $V$ is the $f_1$-\emph{lift} of $U$, and we write $V = f_1(U)$.

Conversely, if we are given $U$ and we want to find its $f_1$-lift $V$, then the type parameters of $V$ will be given by $b_1 = a_4 + a_1 + a_2$ and $b_i = a_i$ for all $i \in \{2, 3, 4\}$; while the explicit form of $V$ itself, in terms of $U$, will be \[U_6^{-1}U_5^{-1}U_4^{-1} : U_2 : U_3 : U_4 : U_2^{-1}U_1^{-1}U_8^{-1} : U_6 : U_7 : U_8.\]

It is straightforward to check that, if $U$ is an interleaved double chain and $V$ is defined as above, then $V$ will be an interleaved double chain as well.

\begin{figure}[ht] \null \hfill \begin{subfigure}[c]{100pt} \centering \includegraphics{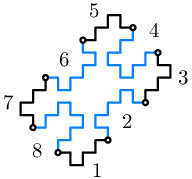} \caption{} \label{lift-f1-a} \end{subfigure} \hfill \begin{subfigure}[c]{130pt} \centering \includegraphics{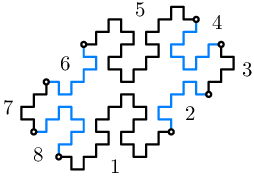} \caption{} \label{lift-f1-b} \end{subfigure} \hfill \null \caption{} \label{lift-f1} \end{figure}

For example, one possible $U$ with $\mathcal{A} = \mathcal{A}_\square$ is shown in Figure \ref{lift-f1-a}, and its corresponding $V$ is in Figure \ref{lift-f1-b}. (Or, strictly speaking, what is shown in these figures are the polyline chains which correspond to the word chains $U$ and $V$ under the standard valuation $\psi_\square$ for $\mathcal{A}_\square$.)

Notice that any given $U$ admits an $f_1$-lift. Indeed, the construction rules for an $f_1$-lift involve only reversal and concatenation, and so we can apply them to any interleaved double chain whatsoever. By contrast, only some $V$ admit an $f_1$-reduction. Specifically, a given $V$ admits an $f_1$-reduction if and only if its type parameters satisfy the constraint $b_1 \ge b_2 + b_4$.

One way to think about the issue is that the construction rules for an $f_1$-reduction involve not only reversal and concatenation, but also the deletion of certain prefixes and suffixes. These deletions can only be carried out if the words we are deleting from are long enough to accommodate them. This point of view will be important in Sections \ref{iw} and \ref{gt}.

The definitions of an $f_i$-reduction and an $f_i$-lift are analogous for all $i = 2$, $3$, $4$.

Since these transformations are going to be crucial for our work, let us write out the general case in detail. Suppose that we are given an interleaved double chain $U$, and we want to construct its $f_i$-lift $V$. The type parameters of $V$ will then be given by $b_i = a_{i - 1} + a_i + a_{i + 1}$ and $b_j = a_j$ for all $j \neq i$; whereas, using shorthand, the parts of $V$ will be \[\begin{aligned} V_j &= U_j & &\text{for all $j \not \equiv i \mymod 4$}\\ V_j &= U_{(j + 3) : (j + 5)}^{-1} & &\text{for all $j \equiv i \mymod 4$}. \end{aligned}\]

Collectively, we refer to $f_1$, $f_2$, $f_3$, $f_4$ as the four $f$-\emph{lifts} or $f$-\emph{transforms}.

Notice that $f_1$ and $f_3$ commute. Indeed, \[f_1(f_3(U)) = f_3(f_1(U)) = U_{4 : 6}^{-1} : U_2 : U_{6 : 8}^{-1} : U_4 : U_{8 : 2}^{-1} : U_6 : U_{2 : 4}^{-1} : U_8.\]

\begin{figure}[ht] \centering \includegraphics{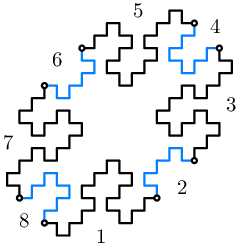} \caption{} \label{lift-fodd} \end{figure}

Similarly, $f_2$ and $f_4$ commute as well. We use the term $f$-\emph{commutativity} to refer collectively to $\{f_1, f_3\}$-commutativity and $\{f_2, f_4\}$-commutativity. We also introduce the notations \begin{align*} \fodd &= f_1 \circ f_3 = f_3 \circ f_1\\ \feven &= f_2 \circ f_4 = f_4 \circ f_2 \end{align*} for the compositions of these pairs. For example, Figure \ref{lift-fodd} shows the $\fodd$-lift of the interleaved double chain depicted in Figure \ref{lift-f1-a}.

\subsection{Imaginary Words} \label{iw}

Let $\mathcal{A}$ be an alphabet with reversal and let $\mathfrak{A}$ be the free group generated by the letters of~$\mathcal{A}$. We call the elements of $\mathfrak{A}$ \emph{imaginary words} over $\mathcal{A}$, and we call the group operation of $\mathfrak{A}$ \emph{imaginary concatenation}.

(The choice of term is by analogy with the historical development of complex numbers. Soon after the discovery of a general method for the solution of cubic equations, it was observed that in some instances the intermediate calculations involve ``imaginary quantities'' even though the end result works out to a real number anyway. This is similar to how we will be using imaginary words in the next sub-section.)

We identify each word over $\mathcal{A}$ with its corresponding group product in $\mathfrak{A}$. So, to us, each word over $\mathcal{A}$ is also an imaginary word over $\mathcal{A}$. In the context of imaginary words, sometimes we refer to the words of $\mathcal{A}$ as ``ordinary'' words.

Let $w$ be a letter of $\mathcal{A}$. We write $\widehat{w}$ for the group inverse of $w$ in $\mathfrak{A}$. (We cannot use $w^{-1}$ because this notation has already been reserved for the reverse of $w$ in $\mathcal{A}$.) Similarly, if $W$ is an imaginary word, we write $\widehat{W}$ or $\myneg(W)$ for the group inverse of $W$ in~$\mathfrak{A}$. We call this the \emph{negation} of $W$.

We extend reversal to imaginary words by $\widehat{w}^{-1} = \widehat{w^{-1}}$ for all letters $w$ of $\mathcal{A}$ and $(W_1W_2)^{-1} = W_2^{-1}W_1^{-1}$ for all imaginary words $W_1$ and $W_2$.

For example, say that $\mathcal{A} = \{\mathsf{A}, \mathsf{B}, \mathsf{C}, \mathsf{D}\}$ with $\mathsf{A} \reveq \mathsf{B}$ and $\mathsf{C} \reveq \mathsf{D}$. Then one imaginary word over $\mathcal{A}$ would be $\mathsf{A}\mathsf{B}\widehat{\mathsf{C}}$. Its reverse works out to $\widehat{\mathsf{D}}\mathsf{A}\mathsf{B}$, and its negation becomes $\mathsf{C}\widehat{\mathsf{B}}\widehat{\mathsf{A}}$.

We define a \emph{negative word} to be the negation of an ordinary word. Or, equivalently, a negative word is any group product, in $\mathfrak{A}$, of the negations of the letters of $\mathcal{A}$. In this context, sometimes we call the words of $\mathcal{A}$ \emph{positive words}. (So an ``ordinary'' word is the same thing as a ``positive'' word.)

The main point of imaginary concatenation is that it allows us to introduce negative words. Then we can write down concatenations between positive and negative words which cause some parts of the positive words to be erased. For example, consider the ordinary word $\mathsf{A}\mathsf{B}\mathsf{C}\mathsf{D}$ and the negative word $\widehat{\mathsf{D}}\widehat{\mathsf{C}}$. Their concatenation, in this order, will be the ordinary word $\mathsf{A}\mathsf{B}$.

More generally, suppose that $W_1$ and $W_2$ are two ordinary words. If $W_2$ is a prefix of $W_1$, then the imaginary concatenation $\widehat{W_2}W_1$ erases this prefix from the beginning of $W_1$. Similarly, if $W_2$ is a suffix of $W_1$, then the imaginary concatenation $W_1\widehat{W_2}$ erases this suffix from the end of $W_1$. In both of these cases, the result of the imaginary concatenation will be an ordinary word once again.

We can also attach some geometric intuitions to these notions. Consider, for example, the alphabet $\mathcal{A}_\square$. We can think of the word $\mathsf{R}\mathsf{U}$ as the instruction ``draw a unit segment pointing to the right, followed by a unit segment pointing up''. Then its negation, $\widehat{\mathsf{U}}\widehat{\mathsf{R}}$, will be the instruction to erase two such segments.

Of course, when we concatenate drawing and erasing instructions, sometimes we might find ourselves having to erase a polyline which is not there. This is why we cannot get away with positive and negative words only, and instead we must introduce the full space of imaginary words.

Notice that there is a big difference between, say, $\mathsf{R}\mathsf{R}^{-1}$ and $\mathsf{R}\widehat{\mathsf{R}}$. Geometrically, the former signifies a closed polyline of length two units which doubles back; whereas the latter signifies an empty polyline.

It will be convenient to define the length of an imaginary word by $|\widehat{w}| = -1$ for all letters $w$ of $\mathcal{A}$ and $|W_1W_2| = |W_1| + |W_2|$ for all imaginary words $W_1$ and $W_2$.

By analogy with chains of ordinary words, we consider chains of imaginary words, too. We call them \emph{pseudo-chains}. The definitions and notations carry over in an obvious manner. We define an \emph{interleaved double pseudo-chain} to be one which satisfies the system of equations (\textbf{I}), and the definition of a \emph{non-interleaved double pseudo-chain} is similar. Notice that, for a pseudo-chain, the equations of (\textbf{I}) and (\textbf{NI}) must all be reinterpreted in terms of imaginary reversal and imaginary concatenation.

\subsection{The \texorpdfstring{$g$}{g}-Transforms} \label{gt}

Just as before, let $V = V_1 : V_2 : \cdots : V_8$ be an interleaved double chain of type $(b_1, b_2, b_3, b_4)$ over the alphabet with reversal $\mathcal{A}$.

Recall the definition of $W_1$ from Section \ref{ft}. Provided that $b_1 \ge b_2 + b_4$, we obtain $W_1$ from $V_1$ by erasing a prefix and a suffix identical to $V_6^{-1}$ and $V_4^{-1}$, respectively. Now, equipped with the notion of imaginary concatenation, we can rewrite the same definition more concisely as \[W_1 = \widehat{V_6^{-1}}V_1\widehat{V_4^{-1}}.\]

However, notice that this new definition makes sense for all $V$; we do not need to impose any constraints on the chain's type parameters. On the flip side, for the new definition we cannot guarantee that it will always yield an ordinary word.

Let us investigate this question a bit more closely. We already know that, if $b_1 \ge b_2 + b_4$, then $W_1$ works out to an ordinary word. We claim that, furthermore, if $b_1 \le b_2 + b_4$, then $W_1$ becomes a negative word.

It will be more convenient for the proof to work with the negation of $W_1$ instead. So we denote $Z_1 = \widehat{W_1} = V_4^{-1}\widehat{V_1}V_6^{-1}$, and now we want to show that $Z_1$ is an ordinary word.

Since $b_1 \le b_2 + b_4$, there exist two nonnegative integers $b'$ and $b''$ such that $b' \le b_2$, $b'' \le b_4$, and $b_1 = b' + b''$. Let, then, $V_1 = V'V''$ with $|V'| = b'$ and $|V''| = b''$. Since $b' \le b_2$ and $V'$ is a prefix of $V_1$, from the equation $V_1V_2 = V_6^{-1}V_5^{-1}$ of (\textbf{I}) we get that $V'$ is a prefix of $V_6^{-1}$, too. Similarly, since $b'' \le b_4$ and $V''$ is a suffix of $V_1$, from the equation $V_8V_1 = V_5^{-1}V_4^{-1}$ of (\textbf{I}) we get that $V''$ is a suffix of $V_4^{-1}$ as well.

Hence, both of $V_4^{-1}\widehat{V''}$ and $\widehat{V'}V_6^{-1}$ are ordinary words. But their concatenation is $V_4^{-1}\widehat{V''}\allowbreak\widehat{V'}V_6^{-1} = V_4^{-1}\widehat{V_1}V_6^{-1} = Z_1$, and so $Z_1$ must be an ordinary word, too, as desired.

Similarly, for all $i$ we define \[Z_i = V_{i + 3}^{-1}\widehat{V_i}V_{i + 5}^{-1}.\]

By an analogous argument, $Z_i$ will be an ordinary word whenever $b_i \le b_{i - 1} + b_{i + 1}$. (Conversely, it will be a negative word whenever $b_i \ge b_{i - 1} + b_{i + 1}$. Either way, $|Z_i| = b_{i - 1} + b_{i + 1} - b_i$.)

In Section \ref{ft}, we saw that if $b_1 \ge b_2 + b_4$ then $W_1$ and $W_5$ satisfy $W_1V_4^{-1} = V_8W_5^{-1}$. In general, when $V$ is arbitrary and no conditions have been imposed on its type parameters, the same calculations yield $\widehat{Z_1}V_4^{-1} = V_8\widehat{Z_5}^{-1}$. (Though notice that, this time around, the expressions on both sides use imaginary concatenation.) Or, rearranging so as to get rid of the negations, $Z_1V_8 = V_4^{-1}Z_5^{-1}$.

Continuing to reason analogously, we obtain the identities \[\begin{aligned} Z_{i + 1}V_i &= V_{i + 4}^{-1}Z_{i + 5}^{-1}\\ V_{i + 1}Z_i &= Z_{i + 4}^{-1}V_{i + 5}^{-1} \end{aligned} \qquad \text{for all $i$}.\]

Suppose now that $b_i \le b_{i - 1} + b_{i + 1}$ for all odd $i$. Then also $Z_i$ becomes an ordinary word for all odd~$i$. Therefore, by the identities listed above, we get that \[Z_1^{-1} : V_2^{-1} : Z_3^{-1} : V_4^{-1} : Z_5^{-1} : V_6^{-1} : Z_7^{-1} : V_8^{-1}\] is an interleaved double chain over $\mathcal{A}$.

Denote this chain by $U = U_1 : U_2 : \cdots : U_8$, and let its type be $(a_1, a_2, a_3, a_4)$. So the type parameters of $U$ are given by $a_i = b_{i - 1} + b_{i + 1} - b_i$ for all odd $i$ and $a_i = b_i$ for all even $i$; whereas the parts of $U$ are constructed by means of \[\begin{aligned} U_i &= V_{i + 5}\widehat{V_i^{-1}}V_{i + 3} & &\text{for all odd $i$}\\ U_i &= V_i^{-1} & &\text{for all even $i$}. \end{aligned}\]

We call $U$ the $\godd$-\emph{transform} of $V$, and we write $U = \godd(V)$. (Shortly, we will define a $\geven$-transform, too.)

We claim that the $\godd$-transform is its own inverse. Or, in other words, $U = \godd(V)$ if and only if $V = \godd(U)$. (This is why in the previous paragraph we had ``$U$ is the $\godd$-transform of~$V$'' and ``$U = \godd(V)$'' where, by analogy with Section \ref{ft}, one might otherwise have expected to see ``$U$ is the $\godd$-reduction of $V$'' and ``$U = \godd^{-1}(V)$''.) Notice the stark contrast with the $f$-transforms, for which ``going up'' and ``going down'' are two very different transformations.

\begin{figure}[ht] \null \hfill \begin{subfigure}[c]{100pt} \centering \includegraphics{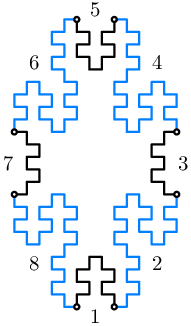} \caption{} \label{lift-godd-a} \end{subfigure} \hfill \begin{subfigure}[c]{130pt} \centering \includegraphics{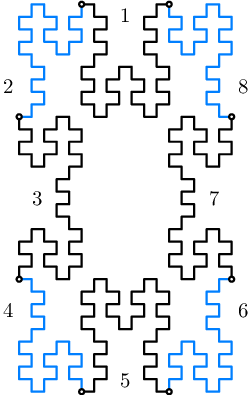} \caption{} \label{lift-godd-b} \end{subfigure} \hfill \null \caption{} \label{lift-godd} \end{figure}

By way of an example, Figure \ref{lift-godd} shows two interleaved double chains $U$ and $V$ with $\mathcal{A} = \mathcal{A}_\square$ each one of which is the $\godd$-transform of the other.

For the proof, observe first of all that the type parameters of $U$ also satisfy the constraints that we imposed on the type parameters of $V$; i.e., $a_i \le a_{i - 1} + a_{i + 1}$ for all odd $i$. This means that $\godd(U)$ is well-defined.

We must now verify the identities \[\begin{aligned} V_i &= U_{i + 5}\widehat{U_i^{-1}}U_{i + 3} & &\text{for all odd $i$}\\ V_i &= U_i^{-1} & &\text{for all even $i$}. \end{aligned}\]

For each even $i$, the verification is trivial. Since the odd $i$ are all analogous to one another, we are only going to carry out the calculations with $i = 1$. Then \begin{align*} U_6\widehat{U_1^{-1}}U_4 &= V_6^{-1}\myneg(V_6\widehat{V_1^{-1}}V_4)^{-1}V_4^{-1}\\ &= V_6^{-1}\widehat{V_6^{-1}}V_1\widehat{V_4^{-1}}V_4^{-1} = V_1, \end{align*} as needed. This confirms our claim that the $\godd$-transform is its own inverse.

Even so, we are still going to define $\godd$-reductions and $\godd$-lifts. We say that $U$ is a $\godd$-\emph{reduction} of $V$ when $|U| \le |V|$. This occurs if and only if $b_1 + b_3 \ge b_2 + b_4$. Conversely, we say that $U$ is a $\godd$-\emph{lift} of $V$ when $|U| \ge |V|$. The necessary and sufficient condition now becomes $b_1 + b_3 \le b_2 + b_4$.

(So, in particular, if $b_1 + b_3 = b_2 + b_4$, then the same $\godd$-transform will be classified simultaneously as a reduction and a lift. A similar scenario is possible also for the $f$-transforms, in the degenerate case where $b_1 = b_3 = 0$ or $b_2 = b_4 = 0$. We will revisit these exceptional cases in the next sub-section.)

The opposite-parity analogue of the $\godd$-transform is the $\geven$-transform. An interleaved double chain $U$ of type $(a_1, a_2, a_3, a_4)$ admits a $\geven$-transform if and only if its type parameters satisfy the constraints $a_i \le a_{i - 1} + a_{i + 1}$ for all even $i$. The formulas for the type parameters of $\geven(U)$, as well as the construction rules which produce the parts of $\geven(U)$, are all fully analogous to the ones for $\godd$. The only difference is that the words ``even'' and ``odd'' are swapped everywhere. Once again, $\geven$ is its own inverse, in the sense that $V = \geven(U)$ if and only if $U = \geven(V)$. The definitions of a $\geven$-reduction and a $\geven$-lift are both analogous to the corresponding definitions for $\godd$, too.

Collectively, we refer to $\godd$ and $\geven$ as the two $g$-\emph{transforms}.

\subsection{Roots and Loops} \label{rl}

Let $U = U_1 : U_2 : \cdots : U_8$ be an interleaved double chain of type $(a_1, a_2, a_3, a_4)$ over the alphabet with reversal $\mathcal{A}$.

An $f_i$-reduction, if admissible, will decrease $|U|$ by $2(a_{i - 1} + a_{i + 1})$. So the reduced chain will almost always be shorter. The only exception is the degenerate case where $a_{i - 1} = a_{i + 1} = 0$; in it, both of $U$ and $|U|$ will remain the same. We call such $f$-reductions \emph{false}.

Suppose now that $U$ does not admit any $f$-reductions -- neither true nor false ones. Then $a_i < a_{i - 1} + a_{i + 1}$ for all $i$. But these constraints imply that $U$ must admit both a $\godd$-transform and a $\geven$-transform.

If $a_1 + a_3 \ge a_2 + a_4$, then the $\godd$-transform will be a reduction and the $\geven$-transform will be a lift. Conversely, if $a_1 + a_3 \le a_2 + a_4$, then the $\geven$-transform will be a reduction and the $\godd$-transform will be a lift. Either way, $U$ necessarily admits at least one $g$-reduction.

When $a_1 + a_3 > a_2 + a_4$, the $\godd$-reduction of $U$ decreases its length by $4(a_1 - a_2 + a_3 - a_4)$, and so we get a shorter reduced chain. Similarly, when $a_1 + a_3 < a_2 + a_4$, the $\geven$-reduction of $U$ produces a shorter reduced chain instead.

This only leaves the exceptional case when $a_1 + a_3 = a_2 + a_4$. In it, both $g$-transforms qualify as reductions but neither one of them actually decreases the length of the chain. We call $g$-reductions of this form \emph{false} as well. Notice that a false $g$-reduction will, in general, alter the chain even though it will not change its length -- in contrast to false $f$-reductions, which do not affect the chain at all.

So, to summarise, each interleaved double chain admits at least one reduction. However, it might happen that all admissible reductions are false, and none of them actually decrease the chain's length.

Beginning with $U$, let us reduce our interleaved double chain iteratively, using the following algorithm on each step:

First, we check whether we are in one of the exceptional cases $a_1 = a_3 = 0$, $a_2 = a_4 = 0$, or $a_1 + a_3 = a_2 + a_4$. If yes, then we stop.

Else, we check whether any $f$-reductions are possible. If yes, then we choose one of them arbitrarily. It will necessarily be a true $f$-reduction; we apply it immediately, and we get a shorter interleaved double chain.

Otherwise, if $U$ does not admit any $f$-reductions, then it must admit both $g$-transforms. Furthermore, exactly one of them will be a true $g$-reduction. We apply it immediately, and once again our interleaved double chain becomes shorter.

This completes the description of the algorithm.

Since the length of a chain is a nonnegative integer, it cannot decrease indefinitely. Hence, sooner or later we are going to encounter one of our stopping conditions.

Suppose first that we stop at a point where either $a_1 = a_3 = 0$ or $a_2 = a_4 = 0$. We call an interleaved double chain $U$ of this form a \emph{root}. Suppose, for concreteness, that $a_2 = a_4 = 0$; the opposite case is analogous. Then the structure of $U$ must be \[X : \varepsilon : Y : \varepsilon : X^{-1} : \varepsilon : Y^{-1} : \varepsilon,\] where $X$ and $Y$ are two words over $\mathcal{A}$. Conversely, every eight-part chain of this form is indeed an interleaved double chain over $\mathcal{A}$.

Geometrically, in the setting of Section \ref{dcp}, a root would correspond to a scenario where the ``two'' tilings of $\Phi$ are in fact one and the same, and the two quadruples of division points $\{P_1, P_2, P_3, P_4\}$ and $\{Q_1, Q_2, Q_3, Q_4\}$ coincide.

Suppose, now, that we stop instead at a point where $a_1 + a_3 = a_2 + a_4$. We call an interleaved double chain $U$ of this form a \emph{loop}.

The reasoning behind this choice of term is as follows: Observe, to begin with, that a false $\godd$-reduction swaps $a_1$ and $a_3$; whereas a false $\geven$-reduction swaps $a_2$ and $a_4$. Either way, the transformed chain will still be a loop. This means that any series of $\godd$-transforms and $\geven$-transforms is admissible from a loop; but, on the flip side, any such series will consist of false $g$-reductions only, and the interleaved double chains that we obtain along the way will all be loops as well. Furthermore, these chains will all be of the same length; and so, if we keep going, eventually we will encounter a repeat. Or, in other words, the process will ``fall into a loop''.

Thus, in summary, we have established that every interleaved double chain can be reduced, iteratively, either to a root or to a loop.

Suppose that $U$ eventually reduces to $B$, with $B$ being either a root or a loop. Conversely, some sequence of $f$-lifts and $g$-lifts $\bfphi = (\varphi_1, \varphi_2, \ldots, \varphi_k)$ must transform $B$ back into $U$; i.e., \[U = \varphi_k(\varphi_{k - 1}( \cdots \varphi_1(B) \cdots )) = \varphi_k \circ \varphi_{k - 1} \circ \cdots \circ \varphi_1(B).\]

Then we say that $B$ is a \emph{base} for $U$, and that $\bfphi$ is a \emph{descent} for $U$ from $B$. (Though we are going to expand the definition of a ``descent'' in the next sub-section. Strictly speaking, here we are using this term as shorthand for ``descent from some base using the set of transformations $\{f_1, f_2, f_3, f_4, \geven, \godd\}$''.) For convenience, given a sequence of transformations $\bfphi = (\varphi_1, \varphi_2, \ldots, \varphi_k)$, from now on we will write just $U = \bfphi(B)$ as an abbreviation for $U = \varphi_k \circ \varphi_{k - 1} \circ \cdots \circ \varphi_1(B)$.

\subsection{Statement of the Main Theorem} \label{smt}

The \emph{Greek cross} is the polyomino obtained when four unit squares are attached to the sides of a fifth, central one, as shown in Figure \ref{step-a}. It admits exactly two tilings of the plane. Both of them are $4$-neighbour lattice tilings; they use the knight-move lattices $\mathcal{L}((2, 1), (-1, 2))$ and $\mathcal{L}((1, 2), (-2, 1))$. (The two tilings are mirror images of one another. However, they are distinct upon translation, and so they count as distinct for our purposes.)

\begin{figure}[ht] \centering \includegraphics{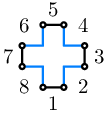} \caption{} \label{cross} \end{figure}

One double chain over $\mathcal{A}_\square$ which describes the Greek cross as a $4$-neighbour double tile is \[\maltese = \mathsf{R} : \mathsf{U}\mathsf{R} : \mathsf{U} : \mathsf{L}\mathsf{U} : \mathsf{L} : \mathsf{D}\mathsf{L} : \mathsf{D} : \mathsf{R}\mathsf{D};\] it is interleaved and of type $(1, 2, 1, 2)$. Figure \ref{cross} shows the corresponding double chain of polylines. For convenience, from now on we will use the term ``Greek cross'' both for the polyomino and for the double chain $\maltese$ over $\mathcal{A}_\square$.

We define the transformations $\gso$ and $\gse$ by \begin{align*} \gso &= \godd \circ \feven\\ \gse &= \geven \circ \fodd. \end{align*}

Collectively, we call $\gso$ and $\gse$ the two $g^\star$-\emph{transforms}.

Let $U = U_1 : U_2 : \cdots : U_8$ be an interleaved double chain of type $(a_1, a_2, a_3, a_4)$ over the alphabet with reversal $\mathcal{A}$. Recall that the four $f$-lifts of $U$ are always well-defined, but its $g$-transforms might not be. In this respect, $\gso$ and $\gse$ resemble the $f$-lifts. Both of them are total, in the sense that they are well-defined for all $U$.

Let us verify this claim for $\gso$; the case of $\gse$ is analogous. Since $\feven = f_2 \circ f_4$, the $\feven$-lift of $U$ is certainly well-defined. Denote its type by $(b_1, b_2, b_3, b_4)$. Then $b_i = a_i$ for all odd $i$ and $b_i = a_{i - 1} + a_i + a_{i + 1}$ for all even $i$. Hence, $b_i \le b_{i - 1} + b_{i + 1}$ for all odd $i$, and so the $\godd$-transform of $\feven(U)$ is indeed well-defined, too. Furthermore, in this instance the $\godd$-transform will in fact be a $\godd$-lift, by virtue of $b_1 + b_3 \le b_2 + b_4$.

Since the two $g^\star$-transforms are going to be crucial for the statement of our main result, let us write out their construction rules explicitly. Say that $V = \gso(U) = V_1 : V_2 : \cdots : V_8$. Then the parts of $V$ are given by \[\begin{aligned} V_i &= U_{(i - 2) : (i + 2)}^{-1} & &\text{for all odd $i$}\\ V_i &= U_{(i + 3) : (i + 5)} & &\text{for all even $i$}. \end{aligned}\]

The construction rules for $\gse$ are analogous. The only difference is that the words ``odd'' and ``even'' are swapped.

Or, alternatively, unpacking the shorthand, we find that $\gso(U)$ works out to \begin{gather*} U_3^{-1}U_2^{-1}U_1^{-1}U_8^{-1}U_7^{-1} \; : \; U_5U_6U_7 \; : \; U_5^{-1}U_4^{-1}U_3^{-1}U_2^{-1}U_1^{-1} \; : \; U_7U_8U_1 \; : \; {}\\ {} \; : \; U_7^{-1}U_6^{-1}U_5^{-1}U_4^{-1}U_3^{-1} \; : \; U_1U_2U_3 \; : \; U_1^{-1}U_8^{-1}U_7^{-1}U_6^{-1}U_5^{-1} \; : \; U_3U_4U_5. \end{gather*}

Once again, the expanded form of $\gse(U)$ is similar.

It is clear now that the construction rules for the two $g^\star$-transforms can be expressed entirely in terms of reversal and concatenation. This is true of the four $f$-transforms as well. By contrast, the construction rules for the two $g$-transforms require also negation and imaginary concatenation.

Let $\Psi$ be any set of transformations chosen out of the $f$-transforms, the $g$-transforms, and their compositions. We define a $\Psi$-\emph{descendant} of $U$ to be any interleaved double chain obtained from $U$ by means of some series of applications of the transformations of $\Psi$.

Or, more formally: We say that $V$ is a $\Psi$-descendant of $U$ when there exists some sequence $\bfphi = (\varphi_1, \varphi_2, \ldots, \varphi_k)$ of transformations, all chosen out of $\Psi$, with $V = \bfphi(U)$. (Of course, this implies that, in particular, $\varphi_i \circ \varphi_{i - 1} \circ \cdots \circ \varphi_1(U)$ is well-defined for all $1 \le i \le k$.) We also say that $\bfphi$ is a $\Psi$-\emph{descent} for $V$ from $U$.

By way of an example, if \[\psiall = \{f_1, f_2, f_3, f_4, \godd, \geven\},\] then Section \ref{rl} tells us that every interleaved double chain over $\mathcal{A}$ is a $\psiall$-descendant either of some root or of some loop.

We are ready now to state our main result, in its version for polyominoes:

\begin{theorem} \label{mt} Every polyomino which admits two different $4$-neighbour lattice tilings of the plane is a deformation of some polyomino described by an $\{f_1, f_2, f_3, f_4, \gso, \gse\}$-descendant of the Greek cross. Conversely, each such descendant of the Greek cross describes a polyomino, and every deformation of such a polyomino admits two different $4$-neighbour lattice tilings of the plane. \end{theorem}

For convenience, let \[\psimain = \{f_1, f_2, f_3, f_4, \gso, \gse\},\] and define a \emph{descendant chain} to be any $\psimain$-descendant of the Greek cross. Notice that all six transformations of $\psimain$ are total, and so we can apply them to the Greek cross any number of times and in any order. Or, in other words, for every finite sequence $\bfphi$ over $\psimain$, the expression $\bfphi(\maltese)$ is well-defined and it works out to a descendant chain.

\begin{figure}[p] \null \hfill
\begin{subfigure}[c]{35pt} \centering \includegraphics{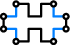} \caption{$f_1$} \label{chain-f1} \end{subfigure} \hfill
\begin{subfigure}[c]{30pt} \centering \includegraphics{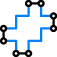} \caption{$f_2$} \label{chain-f2} \end{subfigure} \hfill
\begin{subfigure}[c]{50pt} \centering \includegraphics{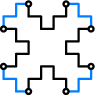} \caption{$\gso$} \label{chain-go} \end{subfigure} \hfill
\begin{subfigure}[c]{55pt} \centering \includegraphics{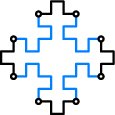} \caption{$\gse$} \label{chain-ge} \end{subfigure} \hfill \null\\ \vspace{\baselineskip}\\ \null \hfill
\begin{subfigure}[c]{50pt} \centering \includegraphics{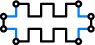} \caption{$f_1$, $f_1$} \label{chain-f1-f1} \end{subfigure} \hfill
\begin{subfigure}[c]{55pt} \centering \includegraphics{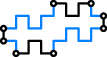} \caption{$f_1$, $f_2$} \label{chain-f1-f2} \end{subfigure} \hfill
\begin{subfigure}[c]{55pt} \centering \includegraphics{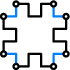} \caption{$f_1$, $f_3$} \label{chain-f1-f3} \end{subfigure} \hfill
\begin{subfigure}[c]{85pt} \centering \includegraphics{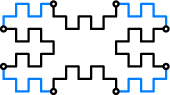} \caption{$f_1$, $\gso$} \label{chain-f1-go} \end{subfigure} \hfill
\begin{subfigure}[c]{95pt} \centering \includegraphics{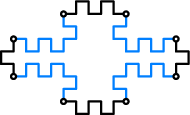} \caption{$f_1$, $\gse$} \label{chain-f1-ge} \end{subfigure} \hfill \null\\ \vspace{\baselineskip}\\ \null \hfill
\begin{subfigure}[c]{50pt} \centering \includegraphics{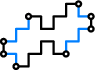} \caption{$f_2$, $f_1$} \label{chain-f2-f1} \end{subfigure} \hfill
\begin{subfigure}[c]{40pt} \centering \includegraphics{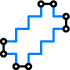} \caption{$f_2$, $f_2$} \label{chain-f2-f2} \end{subfigure} \hfill
\begin{subfigure}[c]{35pt} \centering \includegraphics{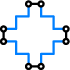} \caption{$f_2$, $f_4$} \label{chain-f2-f4} \end{subfigure} \hfill
\begin{subfigure}[c]{65pt} \centering \includegraphics{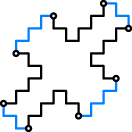} \caption{$f_2$, $\gso$} \label{chain-f2-go} \end{subfigure} \hfill
\begin{subfigure}[c]{75pt} \centering \includegraphics{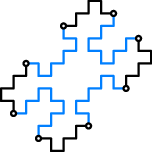} \caption{$f_2$, $\gse$} \label{chain-f2-ge} \end{subfigure} \hfill \null\\ \vspace{\baselineskip}\\ \null \hfill
\begin{subfigure}[c]{75pt} \centering \includegraphics{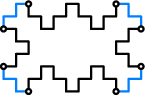} \caption{$\gso$, $f_1$} \label{chain-go-f1} \end{subfigure} \hfill
\begin{subfigure}[c]{65pt} \centering \includegraphics{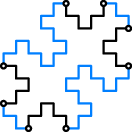} \caption{$\gso$, $f_2$} \label{chain-go-f2} \end{subfigure} \hfill
\begin{subfigure}[c]{80pt} \centering \includegraphics{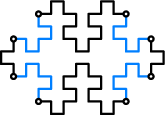} \caption{$\gse$, $f_1$} \label{chain-ge-f1} \end{subfigure} \hfill
\begin{subfigure}[c]{75pt} \centering \includegraphics{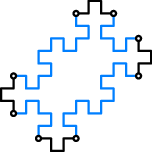} \caption{$\gse$, $f_2$} \label{chain-ge-f2} \end{subfigure} \hfill \null\\ \vspace{\baselineskip}\\ \null \hfill
\begin{subfigure}[c]{130pt} \centering \includegraphics{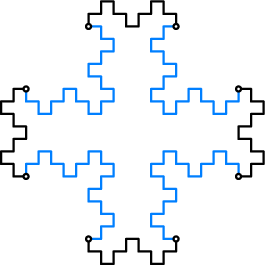} \caption{$\gso$, $\gse$} \label{chain-go-ge} \end{subfigure} \hfill
\begin{subfigure}[c]{115pt} \centering \includegraphics{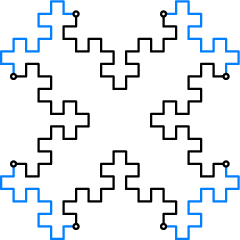} \caption{$\gse$, $\gso$} \label{chain-ge-go} \end{subfigure} \hfill \null
\caption{} \label{tree} \end{figure}

\begin{figure}[t!] \ContinuedFloat \null \hfill
\begin{subfigure}[c]{110pt} \centering \includegraphics{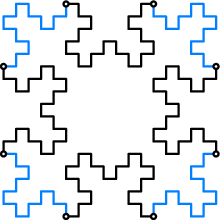} \caption{$\gso$, $\gso$} \label{chain-go-go} \end{subfigure} \hfill
\begin{subfigure}[c]{140pt} \centering \includegraphics{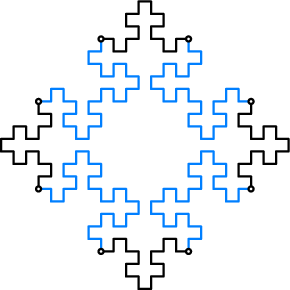} \caption{$\gse$, $\gse$} \label{chain-ge-ge} \end{subfigure} \hfill \null
\caption{{ }(continued)} \end{figure}

For example, Figure \ref{tree} shows the beginning of the ``family tree'' of the descendant chains. It lists almost all descendant chains which can be obtained from the Greek cross by means of one or two transformations of $\psimain$. The omissions are all congruent, modulo isometry, to specimens which have been included.

Notice also that Theorem \ref{mt} is to be understood purely in a combinatorial sense; i.e., the tilings referred to in it are supposed to respect the convention of Section \ref{tile} regarding polyomino tilings considered in a combinatorial context. (This restriction will be lifted in Section \ref{poly}, where we are going to be working in a geometric context instead.)

Theorem \ref{mt} is, in essence, a description of all $4$-neighbour double tiles in the polyomino setting. It splits into two separate parts in a natural manner. We proceed now to discuss both of them, in turn.

The first part states that, if a polyomino works, then it must be captured by the description. However, by Lemma \ref{dc}, all working polyominoes correspond to double chains over $\mathcal{A}_\square$. Hence, the gist of the first part is that, if a double chain over $\mathcal{A}_\square$ does not conform to the description, then it must contain a combinatorial self-intersection. We call this the ``negative'' part of Theorem \ref{mt}, and we prove it in Section \ref{neg}.

Conversely, the second -- or ``positive'' -- part of Theorem \ref{mt} states that the double chains which conform to the description do in fact correspond to working polyominoes. This time around, the core claim is that all descendant chains are free of combinatorial self-intersections. The rest follows by Lemma \ref{defl}. We give a proof of this core claim in Section \ref{pos}.

The two parts differ greatly in terms of the ideas required. The negative part is mostly combinatorial, and so in Section \ref{neg} we will be reasoning about manipulations with words. By contrast, the positive part is mostly geometric; accordingly, in Section \ref{pos} our reasoning will revolve instead around manipulations with figures in the plane.

On the face of it, we should expect the negative part to be easier to prove because it only demands of us to find one single combinatorial self-intersection within a given double chain. Conversely, the positive part ought to be tougher because for it we must sift through all subwords, and we must confirm that not even a single one of them is a combinatorial self-intersection. (So, from the point of view of formal logic, the negative part is an $\exists$-statement while the positive part is a $\forall$-statement.) These intuitions seem to be borne out by the contents of Sections \ref{neg} and \ref{pos}.

We introduced Theorem \ref{mt} as the ``polyomino version'' of our main result. There is also a ``polygon version'' which extends the description from polyominoes to general polygons. Only one word needs to be changed in the theorem's statement. In Section \ref{poly}, we do change it, and then we establish the resulting generalisation.

\section{The Negative Part of the Description} \label{neg}

In this section, we prove the negative part of our main theorem.

Section \ref{ni} rules out the non-interleaved double chains. Then we can safely relax back into the theory of the interleaved double chains that we developed over the course of Section \ref{lift}.

In the short Section \ref{stable}, we prove one simple but useful technical lemma. It will allow us to show that certain vast swathes of double chains must contain combinatorial self-intersections, and so can be excluded from consideration.

Section \ref{repack} handles the transition from $\psiall$ to $\psimain$. Effectively, in it we show that if an interleaved double chain is free of combinatorial self-intersections, then we can repackage all $g$-transforms in its descent into $g^\star$-transforms.

Section \ref{loop} deals away with the double chains descended from loops. Following this, only the double chains descended from roots will remain.

Finally, in Section \ref{root} we look more closely into what is going on at the root and in its immediate vicinity. Once we have pinned down the beginning of the descent, the rest of the proof will be straightforward.

For the most part, we will be reasoning in the setting of an arbitrary alphabet with reversal $\mathcal{A}$ instead of just $\mathcal{A}_\square$. Of course, the latter would have sufficed for our short-term goal of establishing the negative part of Theorem \ref{mt}. However, we must also keep in mind our long-term goal of extending Theorem \ref{mt} from polyominoes to general polygons. By and large, the lemmas and proofs of Sections \ref{ni}--\ref{root} will not be affected too much by this choice. The most significant exception is that, in a couple of instances, the phrase ``contains a combinatorial self-intersection'' will need to be reworded along the lines of ``cannot be realised without a geometric self-intersection''.

\subsection{Non-Interleaved Double Chains} \label{ni}

First we show that a non-interleaved double chain cannot possibly describe a double tile. Recall, from Section \ref{dcw}, that if a double chain corresponds to a solution to the double tile problem, then its type $(a_1, a_2, a_3, a_4)$ must satisfy $a_1 + a_3 \neq 0$ and $a_2 + a_4 \neq 0$.

\begin{lemma} \label{nil} Let $U$ be a non-interleaved double chain of type $(a_1, a_2, a_3, a_4)$ with $a_1 + a_3 \neq 0$ and $a_2 + a_4 \neq 0$. Then every realisation of $U$ contains a geometric self-intersection. \end{lemma} 

\begin{myproof} Let $U = U_1 : U_2 : \cdots : U_8$. We split the proof into two cases which behave very differently:

\smallskip

\emph{Case 1}. $a_2 = a_4$. Then the equations of (\textbf{NI}) imply that $U_1 \reveq U_5$, $U_3 \reveq U_7$, and $U_2 = U_6 \reveq U_4 = U_8$. Hence, \[U_{1 : 8} = \tmpcube[\mathsf{A} \to U_1, \mathsf{B} \to U_2, \mathsf{C} \to U_3].\]

Since $U_2$ is nonempty, the desired result now follows by the strong form of Lemma \ref{abcl}.

\smallskip

\emph{Case 2}. $a_2 \neq a_4$. We claim that, in this case, $U$ contains a combinatorial self-intersection.

Suppose, without loss of generality, that $a_2 > a_4$. Then the equation $U_2U_3U_4 = U_8^{-1}U_7^{-1}U_6^{-1}$ of (\textbf{NI}) shows that $U_8^{-1}$ is a prefix of $U_2$ and $U_4^{-1}$ is a prefix of $U_6$. Let \begin{align*} U_2 &= U_8^{-1}X \tag{\textbf{E}$_1$}\\ U_6 &= U_4^{-1}Y, \tag{\textbf{E}$_2$} \end{align*} and notice that both of $X$ and $Y$ are nonempty.

We substitute the right-hand sides of (\textbf{E}$_1$) and (\textbf{E}$_2$) back into the original equation of~(\textbf{NI}), simplify, and arrive at \[XU_3 = U_7^{-1}Y^{-1}. \tag{\textbf{E}$_3$}\]

Now (\textbf{E}$_3$) and the equation $U_3 = U_7^{-1}$ of (\textbf{NI}) together imply that $X$ and $Y^{-1}$ are anagrams. Hence, the word $XY$ is balanced. We proceed to show that the reverse of this word occurs as a combinatorial self-intersection within $U$.

Indeed, by \begin{gather*} U_2U_3 \stackrel{(\textbf{E}_1)}{=\joinrel=} U_8^{-1}XU_3 \stackrel{(\textbf{E}_3)}{=\joinrel=} U_8^{-1}U_7^{-1}Y^{-1}\\ U_4U_5U_6 \stackrel{(\textbf{NI})}{=\joinrel=} U_2^{-1}U_1^{-1}U_8^{-1} \stackrel{(\textbf{E}_1)}{=\joinrel=} X^{-1}U_8U_1^{-1}U_8^{-1} \end{gather*} we get that $Y^{-1}$ is a suffix of $U_{2 : 3}$ and $X^{-1}$ is a prefix of $U_{4 : 6}$. Thus $Y^{-1}X^{-1}$ becomes a subword of $U_{2 : 6}$, and the proof is complete. \end{myproof}

\subsection{Stable Self-Intersections} \label{stable}

Since the non-interleaved double chains are ruled out by Lemma \ref{nil}, from now on we are going to focus on the interleaved double chains. Throughout the rest of this section, let $U = U_1 : U_2 : \cdots : U_8$ be an interleaved double chain of type $(a_1, a_2, a_3, a_4)$ over the alphabet with reversal~$\mathcal{A}$. Just as before, we assume that $a_1 + a_3 \neq 0$ and $a_2 + a_4 \neq 0$; otherwise, $U$ cannot possibly describe a double tile.

We define a \emph{stable self-intersection} in $U$ to be a nonempty balanced subword of some slice of $U$ of the form $U_{i : (i + 3)}$. Notice that we allow $i \le 0$ and $i + 3 \ge 9$, and we also allow the subword to coincide with the slice itself. Of course, if $U$ contains a stable self-intersection, then it must also contain a combinatorial self-intersection.

The point of stable self-intersections is that they are ``heritable'', in the following sense:

\begin{lemma} \label{sl} Suppose that $U$ contains a stable self-intersection. Then so do all $\psimain$-de\-scen\-dants of $U$ as well. \end{lemma} 

\begin{myproof} Let $\varphi \in \psimain$ and $V = \varphi(U) = V_1 : V_2 : \cdots : V_8$. It suffices to show that either $U_{i : (i + 3)}$ or its reverse occurs as a subword within some slice of $V$ of the form $V_{j : (j + 3)}$.

Suppose first that $\varphi = f_1$; the other three $f$-transforms behave similarly. For $i \in \{1, 8\}$, we can set $j = 1$ because \[V_{1 : 4} \stackrel{f_1}{=\joinrel=} U_6^{-1}U_5^{-1}U_4^{-1}U_2U_3U_4 \stackrel{(\textbf{I})}{=\joinrel=} U_6^{-1}U_{8 : 4};\] for $i \in \{2, 3\}$, we can set $j = 2$ because \[V_{2 : 5} \stackrel{f_1}{=\joinrel=} U_2U_3U_4U_2^{-1}U_1^{-1}U_8^{-1} \stackrel{(\textbf{I})}{=\joinrel=} U_{2 : 6}U_8^{-1};\] and the remaining four cases for $i$ are analogous.

Suppose, now, that $\varphi = \gso$; of course, $\gse$ behaves similarly. When $i \in \{1, 2\}$, clearly $U_{i : (i + 3)}^{-1}$ is already a subword of $V_3 = U_{1 : 5}^{-1}$; and the remaining six cases for $i$ are analogous. \end{myproof}

\subsection{Repackaging the \texorpdfstring{$g$}{g}-Transforms} \label{repack}

Our task here will be to ``repackage'' all $g$-transforms in the descent of $U$ into $g^\star$-transforms.

\begin{lemma} \label{rgl} Suppose that $U$ is free of combinatorial self-intersections. Then it is a $\psimain$-descendant either of some root or of some loop. \end{lemma} 

The proof will occupy the rest of the present sub-section.

Just as in Section \ref{rl}, let $B$ be a base for $U$ (either a root or a loop) and let $\bfphi$ be a descent for $U$ from $B$. We can assume without loss of generality that $\bfphi$ has been obtained by means of the algorithm of Section \ref{rl}, with one modification: Each time during the reduction process when we are offered a choice between an $f_1$-reduction and an $f_3$-reduction, immediately we carry out an $\fodd$-reduction; and similarly also for a choice between an $f_2$-reduction and an $f_4$-reduction. So, strictly speaking, $\bfphi$ will really be a $(\psiall \cup \{\fodd, \feven\})$-descent for $U$ from $B$.

We call an instance of $\godd$ in $\bfphi$ \emph{regular} when it is immediately preceded by an instance of $\feven$; and similarly also for $\geven$. Our task now becomes to show that in fact all instances of $\godd$ and $\geven$ in $\bfphi$ are regular. Then we will be able to unify each $\godd$-instance with the preceding $\feven$-instance, so as to obtain an instance of $\gso$; and similarly for $\geven$ as well.

Suppose, for the sake of contradiction, that there exists an irregular instance of $\godd$ or $\geven$ in $\bfphi$. Out of all such irregular instances, consider the last one. In the ``evolutionary history'' of $U$ from $B$ according to $\bfphi$, let the output of this particular $g$-instance be the interleaved double chain $T$. Since there are no further irregular $g$-instances on the road from $T$ to $U$, we conclude that $U$ is a $\psimain$-descendant of $T$.

With this, our task boils down to showing that there exists a stable self-intersection in $T$. It will then follow by Lemma \ref{sl} that $U$ contains a combinatorial self-intersection, and we will get our desired contradiction.

Observe that $\bfphi$ cannot possibly begin with a $g$-transform. Indeed, if $B$ is a root, then one of $\godd(B)$ and $\geven(B)$ is not well-defined; while at the other one the modified algorithm of Section \ref{rl} would have prescribed either an $\fodd$-reduction or an $\feven$-reduction. Otherwise, if $B$ is a loop, then both of $\godd(B)$ and $\geven(B)$ are loops, too, and at both of them the algorithm would have prescribed a stop.

This means that there must exist an $f$-instance in $\bfphi$ earlier than our marked irregular $g$-instance. (Notice that, here, an ``$f$-instance'' means an instance either of one of the four $f$-transforms, or of one of $\fodd$ and $\feven$.) Out of all such $f$-instances, consider the last one. In the evolutionary history of $U$ from $B$ according to $\bfphi$, let the output of this particular $f$-instance be the interleaved double chain $S$. Since there are no further $f$-instances on the road from $S$ to $T$, we conclude that $T$ is a $\{\godd, \geven\}$-descendant of $S$.

Let $(s_1, s_2, s_3, s_4)$ be the type of $S = S_1 : S_2 : \cdots : S_8$. Since $S$ is the output of an $f$-instance in $\bfphi$, it must admit a true $f$-reduction. Suppose, for concreteness, that $S$ admits a true $f_1$-reduction; i.e., $s_1 \ge s_2 + s_4 \neq 0$. The other three cases are analogous.

Consider, next, the run of $g$-instances in $\bfphi$ which transforms $S$ into $T$. The instances of $\godd$ and $\geven$ in this run must alternate, because each $g$-transform is its own inverse and the modified algorithm of Section \ref{rl} guarantees that each lift in the evolutionary history of $U$ (from $B$ according to $\bfphi$) is true. Furthermore, since $s_1 + s_3 \ge s_1 \ge s_2 + s_4$, this run must begin with an instance of $\geven$.

Hence, we define the sequence of interleaved double chains $T^{(1)}$, $T^{(2)}$, $T^{(3)}$, $\ldots$ by $T^{(1)} = S$~and \[\begin{aligned} T^{(i + 1)} &= \geven(T^{(i)}) & &\text{for all odd $i$}\\ T^{(i + 1)} &= \godd(T^{(i)}) & &\text{for all even $i$}. \end{aligned}\]

Our analysis so far shows that $T$ must occur somewhere in this sequence. Let $T = T^{(k)}$, with $k \ge 2$.

The rest of the proof will be easier to express in terms of imaginary words and imaginary concatenation. Let $W$ be the formal $\fodd$-reduction of $S$; i.e., let $W = W_1 : W_2 : \cdots : W_8$ be the interleaved double pseudo-chain defined by \[\begin{aligned} W_i &= \widehat{S_{i - 1}}S_{i + 4}^{-1}\widehat{S_{i + 1}} & &\text{for all odd $i$}\\ W_i &= S_i & &\text{for all even $i$}. \end{aligned}\]

Then, formally, $S = \fodd(W)$. (Of course, on the right-hand side we must reinterpret the construction rules for $\fodd$ in terms of imaginary concatenation.)

Notice that most of the parts of $W$ will in fact be ordinary words. For $W_2$, $W_4$, $W_6$, $W_8$ this is obvious. Furthermore, for $W_1$ and $W_5$ it follows by the same argument as in Sections \ref{ft} and~\ref{gt}, because $s_1 \ge s_2 + s_4$. This leaves only the parts $W_3$ and $W_7$ of $W$. Both of them will be ordinary words, too, if $s_3 \ge s_2 + s_4$. Otherwise, both of them will be negative words.

Let \[W_\divideontimes = W_{1 : 8} = W_1S_7^{-1}W_5S_3^{-1}.\]

Clearly, this is an ordinary word as well. Furthermore, it is nonempty; otherwise, $s_1 = s_2 + s_4$ and $s_3 = 0$, making $S$ a loop and forcing the modified algorithm of Section \ref{rl} to stop at $S$. Finally, $W_\divideontimes$ must be balanced, too, because it is the concatenation of a double pseudo-chain.

To finish the proof, we will show that $T$ contains a stable self-intersection obtained from $W_\divideontimes$ by means either of a reversal or of some cyclic shift.

First of all, let us express $T = T^{(k)}$ in terms of $W$. We claim that, if $k$ is odd, then \[\begin{aligned} T^{(k)}_i &= W_{(i - k + 4) : (i + k + 4)}^{-1} & &\text{for all odd $i$}\\ T^{(k)}_i &= W_{(i - k + 1) : (i + k - 1)} & &\text{for all even $i$}; \end{aligned}\] and, if $k$ is even, then \[\begin{aligned} T^{(k)}_i &= W_{(i - k + 5) : (i + k + 3)} & &\text{for all odd $i$}\\ T^{(k)}_i &= W_{(i - k) : (i + k)}^{-1} & &\text{for all even $i$}. \end{aligned}\]

These formulas are straightforward to establish by induction on $k$.

We will also make use of one simple observation: Suppose that $i \le j$ with $i$, $j \not \equiv 3 \mymod 4$. Then the slice $W_{i : j}$ of $W$ must be an ordinary word. Indeed, each instance of $W_3$ in this slice will be flanked by two instances of $W_2$ and $W_4$, with the concatenation $W_2W_3W_4$ of all three working out to the ordinary word $S_7^{-1}$; and similar reasoning applies also to each instance of $W_7$.

Thus equipped, we proceed to consider the following cases:

\smallskip

\emph{Case 1}. $k \ge 4$. Then $W_{0 : 2k}^{-1}$ is a part of $T^{(k)}$. But $W_\divideontimes = W_{1 : 8}$ is a subword of $W_{0 : 2k}$ by our observation about the slices of $W$.

\smallskip

\emph{Case 2.} $k = 3$. Then \begin{align*} T^{(3)}_{3 : 5} &= W_{4 : 10}^{-1}W_{2 : 6}W_{6 : 12}^{-1} = (W_{4 : 5}W_{6 : 10})^{-1}W_{2 : 6}(W_{6 : 10}W_{11 : 12})^{-1}\\ &\stackrel{(\textbf{I})}{=\joinrel=} W_{6 : 10}^{-1}W_{0 : 1}W_{2 : 6}W_{7 : 8}W_{6 : 10}^{-1} = W_{6 : 10}^{-1}W_{0 : 8}W_{6 : 10}^{-1}. \end{align*}

Since both of $W_{0 : 8}$ and $W_{6 : 10}$ are ordinary words by our observation about the slices of $W$, we conclude that $W_\divideontimes = W_{1 : 8}$ is a subword of $T^{(3)}_{3 : 5}$.

\smallskip

\emph{Case 3}. $k = 2$. Now $T = \geven(S)$, and so the corresponding instance of $\geven$ in $\bfphi$ is irregular. Since $S$ already admits a true $f_1$-reduction, we conclude that it cannot admit a true $f_3$-reduction. I.e., $s_3 < s_2 + s_4$.

By means of manipulations similar to the ones we carried out in Case 2, we find that \begin{align*} T^{(2)}_{2 : 4} &= W_{0 : 4}^{-1}W_{6 : 8}W_{2 : 6}^{-1} = (W_{0 : 1}W_{2 : 4})^{-1}W_{6 : 8}(W_{2 : 4}W_{5 : 6})^{-1}\\ &\stackrel{(\textbf{I})}{=\joinrel=} W_{2 : 4}^{-1}W_{4 : 5}W_{6 : 8}W_{1 : 2}W_{2 : 4}^{-1} = W_{2 : 4}^{-1}W_{4 : 2}W_{2 : 4}^{-1}. \end{align*}

Once again, both of $W_{2 : 4}$ and $W_{4 : 2}$ are ordinary words by our observation about the slices of~$W$. Hence, $W_{4 : 2}$ is a subword of $T^{(2)}_{2 : 4}$. On the other hand, explicitly, $W_{4 : 2} = S_4W_5S_3^{-1}W_1S_2$.

Since $s_3 < s_2 + s_4$, there exist two nonnegative integers $s'$ and $s''$ such that $s' \le s_2$, $s'' \le s_4$, and $s' + s'' = s_3$. Let, then, $S_7^{-1} = XY$ with $|X| = s'$ and $|Y| = s''$. Since $s' \le s_2$ and $X$ is a prefix of $S_7^{-1}$, from the equation $S_2S_3 = S_7^{-1}S_6^{-1}$ of (\textbf{I}) we get $X$ is a prefix of $S_2$, too. Similarly, since $s'' \le s_4$ and $Y$ is a suffix of $S_7^{-1}$, from the equation $S_3S_4 = S_8^{-1}S_7^{-1}$ of (\textbf{I}) we get that $Y$ is also a suffix of $S_4$.

Therefore, the cyclic shift $YW_5S_3^{-1}W_1X$ of $W_\divideontimes$ is a subword of $W_{4 : 2}$ as well as a subword of~$T^{(2)}_{2 : 4}$. The proof of Lemma \ref{rgl} is complete.

\subsection{Loops} \label{loop}

Our next task will be to show that the chains which describe double tiles must be descended from roots, not loops. First, though, we sort out one technical point. Consider a loop $L$ one of whose parts is empty. By the equations of (\textbf{I}), modulo a suitable cyclic shift, $L$ must be of the form \[\varepsilon : X : YX : Y : \varepsilon : X^{-1} : Y^{-1}X^{-1} : Y^{-1}\] for some words $X$ and $Y$ over $\mathcal{A}$. But then $f_3^{-1}(L)$ becomes a root.

Accordingly, we define a loop to be \emph{genuine} when all of its parts are nonempty. We have just seen that, if $U$ is a $\psimain$-descendant of a non-genuine loop, then it is also a $\psimain$-descendant of some root. Hence, from now on we can focus on the genuine loops only.

\begin{lemma} \label{ll} Suppose that $U$ is a $\psimain$-descendant of a genuine loop. Then every realisation of $U$ contains a geometric self-intersection. \end{lemma} 

The proof will take up the rest of the present sub-section.

Let $\bfphi$ be a $\psimain$-descent for $U$ from the genuine loop $L = L_1 : L_2 : \cdots : L_8$ of type $(\ell_1, \ell_2, \ell_3, \ell_4)$. Suppose, without loss of generality, that $\ell_1 = \max\{\ell_1, \ell_2, \ell_3, \ell_4\}$. Since $\ell_1 + \ell_3 = \ell_2 + \ell_4$, also $\ell_3 = \min\{\ell_1, \ell_2, \ell_3, \ell_4\}$. Then let $L_1 = ABC$ with $|A| = \ell_2 - \ell_3$, $|B| = \ell_3$, and $|C| = \ell_4 - \ell_3$.

From the equation $L_1L_2 = L_6^{-1}L_5^{-1}$ of (\textbf{I}), we get that $L_6^{-1}$ coincides with the prefix of $L_1$ of length $\ell_2$; i.e., $L_6 = B^{-1}A^{-1}$. Continuing to reason similarly, eventually we find expressions in terms of $A$, $B$, $C$ for all eight parts of $L$. Explicitly, \begin{align*} L_1 &= ABC & L_2 &= B^{-1}A & L_3 &= B & L_4 &= C^{-1}B^{-1}\\ L_5 &= A^{-1}BC^{-1} & L_6 &= B^{-1}A^{-1} & L_7 &= B & L_8 &= CB^{-1}. \end{align*}

Consider, then, the template \[\tmploop = \mathsf{A}\mathsf{B}\mathsf{C} : \mathsf{B}^{-1}\mathsf{A} : \mathsf{B} : \mathsf{C}^{-1}\mathsf{B}^{-1} : \mathsf{A}^{-1}\mathsf{B}\mathsf{C}^{-1} : \mathsf{B}^{-1}\mathsf{A}^{-1} : \mathsf{B} : \mathsf{C}\mathsf{B}^{-1}.\]

Our observations may now be summarised as \[L = \tmploop[\mathsf{A} \to A, \mathsf{B} \to B, \mathsf{C} \to C].\]

Or, in other words, all loops are instantiations of one and the same fixed template.

Since $\tmploop$ contains a combinatorial self-intersection, so do all genuine loops as well. However, the same is not true in general of their $\psimain$-descendants. For example, $\tmploop$ is itself a loop, and its $\psimain$-descendant $f_3(\tmploop)$ is free of combinatorial self-intersections. Hence, for the proof of Lemma \ref{ll}, subtler reasoning will be required.

To begin with, let $U^\mathsf{ABC} = \bfphi(\tmploop)$. Since the transformations of $\psimain$ can be expressed entirely in terms of reversal and concatenation, all of them commute with substitutions. Thus $U = U^\mathsf{ABC}[\mathsf{A} \to A, \mathsf{B} \to B, \mathsf{C} \to C]$.

We proceed now to work out the general structure of a $\psimain$-descendant of $\tmploop$.

Let $W = w_1w_2 \ldots w_k$ be a word over the alphabet $\{\mathsf{A}, \mathsf{C}, \mathsf{A}^{-1}, \mathsf{C}^{-1}\}$. There are eight distinct ways to intersperse alternating instances of $\mathsf{B}$ and $\mathsf{B}^{-1}$ between the letters of $W$. Explicitly, if $k$ is odd, then one such interspersal is $\mathsf{B}w_1\mathsf{B}^{-1}w_2\mathsf{B}w_3\mathsf{B}^{-1}w_4 \ldots \mathsf{B}w_k\mathsf{B}^{-1}$. Out of it, we get seven more by possibly omitting the initial instance of $\mathsf{B}$; possibly omitting the final instance of $\mathsf{B}^{-1}$; and possibly applying the swapping substitution $[\mathsf{B} \leftrightarrow \mathsf{B}^{-1}]$. The analysis for an even $k$ is analogous. (Some of the interspersals might coincide if $W$ is too short; however, all eight are in fact distinct when $k \ge 2$.)

It is a fun little exercise in notational design to come up with sensible notation for these eight options. We offer the following solution: Each interspersal is denoted by enclosing $W$ within some pair of brackets. Both at the beginning and at the end of $W$, a floor bracket indicates that the nearest interspersed symbol is $\mathsf{B}^{-1}$, and a ceiling bracket indicates that it is $\mathsf{B}$. Furthermore, both at the beginning and at the end of $W$, a thin bracket indicates that the nearest symbol of the interspersal belongs to the original word $W$; while a thick bracket indicates that it is one of $\mathsf{B}$ and $\mathsf{B}^{-1}$. For example, the expression $\lceil W \rrfloor$ is well-defined if and only if $k$ is even, and it denotes the interspersal $w_1\mathsf{B}w_2\mathsf{B}^{-1}w_3 \ldots \mathsf{B}w_k\mathsf{B}^{-1}$.

(These definitions, as written, work smoothly if $k \ge 2$ but must be applied with care when $k \in \{0, 1\}$. Specifically, we must ensure that concatenation rules such as $\lfloor W_1 \rceil \llfloor W_2 \rrceil = \lfloor W_1W_2 \rrceil$ hold for all $W_1$ and $W_2$ whose lengths are of the correct parities for the left-hand side to be well-defined. One minor subtlety is that then necessarily $\lfloor \varepsilon \rfloor = \widehat{\mathsf{B}}$ and $\lceil \varepsilon \rceil = \widehat{\mathsf{B}^{-1}}$. However, these unusual edge cases will not arise in the proof.)

Let $\tmploop^\mathsf{AC} = \tmploop[\mathsf{B} \to \varepsilon]$ and $U^\mathsf{AC} = \bfphi(\tmploop^\mathsf{AC})$. We claim that there exists some $\ell$ such that the $\ell$-th cyclic shifts $S = S_1 : S_2 : \cdots : S_8$ and $T = T_1 : T_2 : \cdots : T_8$ of $U^\mathsf{AC}$ and $U^\mathsf{ABC}$, respectively, are related by either \[T = \lfloor S_1 \rfloor : \llceil S_2 \rceil : \llfloor S_3 \rrfloor : \lceil S_4 \rrceil : \lfloor S_5 \rfloor : \llceil S_6 \rceil : \llfloor S_7 \rrfloor : \lceil S_8 \rrceil \tag{\textbf{L}$_1$}\] or \[T = \lceil S_1 \rceil : \llfloor S_2 \rfloor : \llceil S_3 \rrceil : \lfloor S_4 \rrfloor : \lceil S_5 \rceil : \llfloor S_6 \rfloor : \llceil S_7 \rrceil : \lfloor S_8 \rrfloor\makebox[0pt][l]{.} \tag{\textbf{L}$_2$}\]

This is straightforward to demonstrate by induction on $\bfphi$. In fact, we can be significantly more precise: Setting $\ell = 0$ works when the number of instances of $\gso$ in $\bfphi$ is even; and setting $\ell = 2$ works otherwise. Furthermore, (\textbf{L}$_1$) holds when the total number of instances of the two $g^\star$-transforms in $\bfphi$ is odd; and (\textbf{L}$_2$) holds otherwise.

(It is worth noting that $f_1^{-1}(\tmploop^\mathsf{AC})$ is a root, and so $U^\mathsf{AC}$ is always a $\psimain$-descendant of this root. Furthermore, it follows easily from (\textbf{L}$_1$) and (\textbf{L}$_2$) that $U^\mathsf{ABC}$ contains a combinatorial self-intersection if and only if $U^\mathsf{AC}$ does. This gives us some clarity as to which $\psimain$-descendants of the genuine loops are certain to contain combinatorial self-intersections.)

We are ready now to finish up the proof. Suppose, for concreteness, that $S$ and $T$ satisfy~(\textbf{L}$_1$). The case of (\textbf{L}$_2$) is analogous.

Notice that $T_{1 : 2} = \lfloor S_1 \rfloor \llceil S_2 \rceil = \lfloor S_{1 : 2} \rceil$ and $T_{3 : 4} = \llfloor S_3 \rrfloor \lceil S_4 \rrceil = \mathsf{B}^{-1} \lceil S_{3 : 4} \rfloor \mathsf{B}$. Similar reasoning applies also to $T_{5 : 6}$ and $T_{7 : 8}$. Furthermore, $S_{1 : 2} \reveq S_{5 : 6}$ and $S_{3 : 4} \reveq S_{7 : 8}$ by the equations of (\textbf{I}). We conclude that \[T_{1 : 8} = \tmpcube[\mathsf{A} \to \lfloor S_{1 : 2} \rceil, \mathsf{B} \to \mathsf{B}^{-1}, \mathsf{C} \to \lceil S_{3 : 4} \rfloor].\]

Recall, though, that $T$ is the $\ell$-th cyclic shift of $U^\mathsf{ABC}$ and $U = U^\mathsf{ABC}[\mathsf{A} \to A, \mathsf{B} \to B, \mathsf{C} \to C]$. Furthermore, since $L$ is a genuine loop, the word $B = L_3 = L_7$ must be nonempty. The desired geometric self-intersection is now guaranteed to exist by the strong form of Lemma \ref{abcl}. The proof of Lemma \ref{ll} is complete.

\subsection{Roots} \label{root}

By now, we have narrowed things down to $U$ being a $\psimain$-descendant of some root $R$. Here, we look into the beginning of the descent of $U$, at and near $R$. Suppose, for concreteness, that $R$ is of the form \[X : \varepsilon : Y : \varepsilon : X^{-1} : \varepsilon : Y^{-1} : \varepsilon,\] where $X$ and $Y$ are two words over $\mathcal{A}$. The opposite case, when $R^\shift{1}$ is of this form, is analogous.

\begin{lemma} \label{begin} Suppose that $U$ is free of combinatorial self-intersections. Then, without loss of generality, the $\psimain$-descent of $U$ from $R$ begins with an instance of $f_2$ immediately followed by an instance of $f_4$. \end{lemma}

\begin{myproof} Let $\bfphi$ be a $\psimain$-descent for $U$ from $R$. Recall that the type parameters of $U$ satisfy $a_1 + a_3 \neq 0$ and $a_2 + a_4 \neq 0$. So, in particular, $\bfphi$ cannot be empty, and it cannot begin with an instance of $f_1$ or $f_3$.

Notice, next, that $\gse(R) = \feven(R^\shift{4})$. Hence, if $\bfphi$ begins with an instance of $\gse$, we can replace $R$ with $R^\shift{4}$ as well as the opening $\gse$-instance of $\bfphi$ with the pair $\{f_2, f_4\}$. The desired conclusion will then be true of the modified descent.

On the other hand, \[\gso(R) = Y^{-1}X^{-1}Y : X^{-1}Y^{-1} : XY^{-1}X^{-1} : Y^{-1}X : YXY^{-1} : XY : X^{-1}YX : YX^{-1}\] contains the stable self-intersection $YX^{-1}Y^{-1}X$ in the concatenation of its first three parts. By Lemma \ref{sl}, it follows that $\bfphi$ cannot begin with an instance of $\gso$.

We conclude that $\bfphi$ must begin with an instance of either $f_2$ or $f_4$. Since these two options are analogous, let us assume for concreteness that $\bfphi$ begins with an $f_2$-instance. Let also $k$ be the length of the initial run of $f_2$-instances in $\bfphi$; i.e., $\bfphi$ begins with $k$ consecutive instances of $f_2$, and after that either it ends or else it continues with an instance of some other element of~$\psimain$.

By induction on $k$, it is straightforward to check that \[f_2^k(R) = X : (YX)^k : Y : \varepsilon : X^{-1} : (Y^{-1}X^{-1})^k : Y^{-1} : \varepsilon.\]

Denote this interleaved double chain by $S$. Since it contains the combinatorial self-intersection $XYX^{-1}Y^{-1}$ within $S_{2 : 6}$, we conclude that the initial run of $f_2$-instances in $\bfphi$ cannot coincide with $\bfphi$, and so it must be followed in $\bfphi$ by an instance of some other element $\varphi$ of $\psimain$. (Notice also that this combinatorial self-intersection is not stable as it spans too many parts of $S$.)

Let $T = \varphi(S) = T_1 : T_2 : \cdots : T_8$. Then $T$ must be free of stable self-intersections, because otherwise Lemma \ref{sl} would force a contradiction. We proceed to show that this observation eliminates all but one possibilities for $\varphi$.

Since $f_1$ and $f_3$ behave analogously, we consider only the first one of them. Then \[T = (XY)^kX : (YX)^k : Y : \varepsilon : (X^{-1}Y^{-1})^kX^{-1} : (Y^{-1}X^{-1})^k : Y^{-1} : \varepsilon,\] and we find the stable self-intersection $XYX^{-1}Y^{-1}$ within $T_{2 : 5}$. (Remarkably, this is the only time when we need to apply Lemma \ref{sl} to a stable self-intersection spanning more than three consecutive parts of its corresponding double chain.)

On the other hand, if $\varphi = \gse$, then $T_4 = S_{2 : 6}^{-1}$. So the combinatorial self-intersection we found in $S_{2 : 6}$ will reoccur, reversed, within $T_4$, and this time around it will be stable.

Finally, if $\varphi = \gso$, then \begin{gather*} T = (Y^{-1}X^{-1})^{k + 1}Y : (X^{-1}Y^{-1})^{k + 1} : X(Y^{-1}X^{-1})^{k + 1} : Y^{-1}X : {}\\ {} : (YX)^{k + 1}Y^{-1} : (XY)^{k + 1} : X^{-1}(YX)^{k + 1} : YX^{-1}, \end{gather*} and we discover the stable self-intersection $X^{-1}Y^{-1}XY$ within $T_{3 : 5}$.

Or, to summarise, we have demonstrated by elimination that $\varphi = f_4$. Since $f_2$ and $f_4$ commute, we can now rearrange $\bfphi$ so that it begins as desired. \end{myproof}

\begin{lemma} \label{ril} In the setting of Lemma \ref{begin}, the root $R$ must be free of combinatorial self-intersections as well. \end{lemma} 

\begin{myproof} Suppose, for the sake of contradiction, that $R$ contains a combinatorial self-in\-ter\-sec\-tion. Since the even-numbered parts of $R$ are all empty, we can assume without loss of generality that this combinatorial self-intersection is a subword of $R_{1 : 5}$. Denote $W = \feven(R)$. Then the explicit form of $W$ spelled out below shows that $R_{1 : 5} = XYX^{-1}$ occurs as a subword within $W_{2 : 4}$. Hence, $W$ contains a stable self-intersection, and by Lemma \ref{sl} we arrive at a contradiction. \end{myproof}

We are ready now to finish up the proof of the negative part of Theorem \ref{mt}.

Let $\Phi$ be a double tile in the polyomino setting. Then it is described by some double chain $U$ over $\mathcal{A}_\square$. Clearly, all results of the present section apply to $U$. Thus we can assume, without loss of generality, that $U$ is a $\psimain$-descendant of $\feven(R)$ for some root $R$ as in Lemma \ref{begin}.

Recall the definition of the double chain $\maltese$ from Section \ref{smt}, and observe that \begin{align*} \feven(R) &= X : YX : Y : X^{-1}Y : X^{-1} : Y^{-1}X^{-1} : Y^{-1} : XY^{-1}\\ &= \maltese[\mathsf{R} \to X, \mathsf{U} \to Y]. \end{align*}

Let $\bfphi$ be a $\psimain$-descent for $U$ from $\feven(R)$. Define also $U_\square = \bfphi(\maltese)$; so, in particular, $U_\square$ is a $\psimain$-descendant of the Greek cross by definition. Since the construction rules for the transformations of $\psimain$ can be expressed entirely in terms of reversal and concatenation, all of these transformations commute with substitutions. (We made use of this observation also in Section \ref{loop}.) Therefore, \[U = U_\square[\mathsf{R} \to X, \mathsf{U} \to Y]. \tag{\textbf{R}}\]

By Lemma \ref{ril}, we get that $R$ is free of combinatorial self-intersections, and so it describes some polyomino $\Phi_\text{\rm Root} = [\psi_\square(R)]$. Furthermore, both of $X$ and $Y$ must be nonempty. By (\textbf{R}), it follows that $U_\square$ is free of combinatorial self-intersections, too. (Otherwise, we obtain a combinatorial self-intersection in $U$.) We conclude that $U_\square$ also describes some polyomino $\Phi_\square = [\psi_\square(U_\square)]$.

Then (\textbf{R}) shows that the double tile $\Phi$ is a deformation of the polyomino $\Phi_\square$ which uses the basic building block $\Phi_\text{\rm Root}$.

Or, in summary, we have established that every double tile $\Phi$ in the polyomino setting is a deformation of some polyomino $\Phi_\square$ described by a $\psimain$-descendant $U_\square$ of the Greek cross. This is exactly what we set out to do, and our proof of the negative part of Theorem \ref{mt} is complete.

\section{Basic Properties of the Descendant Chains} \label{prop}

Here, we work out the basic properties of the descendant chains. This material will serve as the foundation for Sections \ref{pos} and \ref{pro}--\ref{clover}.

Section \ref{unique} shows that the correspondence between descents and descendant chains is essentially one-to-one. Next, Section \ref{symm} shows that every descendant chain is centrally symmetric. This is quite helpful; it will often allow us to ``halve'' the amount of work we need to do. In Sections \ref{pv} and \ref{nv}, we introduce two systems of vectors associated with the descendant chains which will help us get a handle on their geometric properties. Continuing this train of thought, in Section \ref{orient} we analyse the orientations of the descendant chains. Finally, in Section \ref{gl} we revisit the construction rules for our lifting transformations from a geometric point of view.

Recall that, by definition, a descendant chain is a syntactic object; specifically, a word chain over $\mathcal{A}_\square$. However, our focus is now going to shift on to the geometric side of things. So, for convenience, henceforth we will use the term ``descendant chain'' also for the corresponding chains of polylines. (Either floating or literal.)

Notice that, since words over $\mathcal{A}_\square$ are in such a close correspondence with the polylines they describe, a lot of the associated apparatus carries over immediately. For example, the type of a descendant chain of polylines is still given by the lengths of its parts.

Throughout this section, let $U = U_1 : U_2 : \cdots : U_8$ be a descendant chain of polylines. Formally, this means that $U$ is either the floating-polyline $\psi_\square$-image of a descendant chain of words over $\mathcal{A}_\square$, or else a literal-polyline representative of this $\psi_\square$-image. Let also $\varphi \in \psimain$ and $V = \varphi(U) = V_1 : V_2 : \cdots : V_8$. Formally, this means that $\varphi$ maps the syntactic descendant chain which describes $U$ onto the syntactic descendant chain which describes $V$. In a setting where we understand $U$ to be a chain of literal polylines, we denote its division points by $P_1$, $P_2$, $\ldots$, $P_8$ so that $P_i$ is the initial point of $U_i$ for all $i$. Similarly for $V$ as well, with the division points being $Q_1$, $Q_2$, $\ldots$, $Q_8$.

For convenience, in the context of descendant chains we are going to use the term ``descent'' as shorthand for ``$\psimain$-descent from the Greek cross''.

\subsection{Uniqueness} \label{unique}

To begin with, we look into the correspondence between descents and descendant chains.

\begin{proposition} \label{ud} The descent of a descendant chain is unique modulo $f$-commutativity. \end{proposition} 

\begin{myproof} Since $|U|$ increases with each lift, the only descent which produces the Greek cross is the empty one. Suppose, from now on, that $U$ is not the Greek cross.

We will show that, by looking solely at the type $(a_1, a_2, a_3, a_4)$ of $U$, we can determine the final lift in its descent modulo $f$-commutativity. Clearly, this would suffice.

We say that $U$ is $f_i$-\emph{feasible} when its type parameters satisfy the constraint $a_i > a_{i - 1} + a_{i + 1}$. We also say that $U$ is $\gso$-\emph{feasible} when its type parameters satisfy the constraint $a_1 + a_3 > a_2 + a_4$ as well as the constraints $a_i < a_{i - 1} + a_{i + 1}$ for all $i$. The definition of $U$ being $\gse$-\emph{feasible} is analogous; the only difference is that, in the first constraint, the direction of the inequality is reversed.

Now, if $U$ admits a descent ending with $\varphi \in \psimain$, then it must be $\varphi$-feasible. However, there are only two pairs of lifting transformations whose feasibility constraints are compatible: $\{f_1, f_3\}$ and $\{f_2, f_4\}$.

Suppose, for concreteness, that $U$ is both $f_1$-feasible and $f_3$-feasible. Observe that $f_1$-fea\-si\-bil\-i\-ty is always preserved by an $f_3$-reduction. Furthermore, the descent of an $f_1$-feasible descendant chain cannot be empty, and it can only end with an instance of one of $f_1$ and $f_3$.

Together, these two observations imply that every descent for $U$ must end with a (possibly empty) run of instances of $f_3$, immediately preceded by an instance of $f_1$. Hence, modulo $f$-commutativity, every descent for $U$ must end with an instance of $f_1$. \end{myproof}

One corollary of the proof is that just the type of $U$ already suffices to pin down the descent of $U$, once again modulo $f$-commutativity. We return to questions of uniqueness in Section \ref{pro}, with Proposition \ref{un}.

\subsection{Symmetry} \label{symm}

Recall, from Section \ref{bc}, that we consider a chain to be symmetric when the corresponding isometry preserves it modulo reversals and cyclic shifts.

\begin{lemma} \label{symml} Every descendant chain $U$ is centrally symmetric. Specifically, a suitable central symmetry swaps $U$ and $U^\shift{4}$. \end{lemma} 

\begin{myproof} By induction on the descent of $U$. The base case is clear. For the induction step, simply observe that, modulo translation, $\varphi$ commutes both with central symmetry and with $4$-step cyclic shift, for all $\varphi \in \psimain$. \end{myproof}

We define the \emph{center} of $U$ to be the center of the central symmetry which swaps $U$ and $U^\shift{4}$. The same central symmetry will also swap $U_i$ and $U_{i + 4}$, for all $i$.

One natural question arises here, and we go on a brief detour in order to discuss it.

Let $W = W_1 : W_2 : \cdots : W_8$ be any double chain of polylines associated with a $4$-neigh\-bour double tile. Then $W$ is a ``deformation'' of a descendant chain by Theorem \ref{mt}. (In the polyomino setting. For the setting of general polygons, we must invoke Theorem \ref{pt} instead.)

Observe that, in general, deformations do not preserve central symmetry. Indeed, $4$-neigh\-bour double tiles are not in general centrally symmetric. However, deformations do preserve the central symmetry of the set of the associated chain's division points. Hence, the set of the division points of $W$ will always be centrally symmetric.

The question, then, is as follows: Does there exist a simpler proof of this fact -- specifically, one which does not rely on Theorem \ref{mt} as a prerequisite?

This question is made more pertinent by the following observation: Recall, from Section \ref{dcw}, that $|W_i| = |W_{i + 4}|$ for all $i$. On the other hand, the above argument shows that also $\sigma(W_i) + \sigma(W_{i + 4}) = \mathbf{0}$ for all $i$. Superficially, these two systems of identities look quite similar. However, the former system is a straightforward corollary of (\textbf{I}); whereas, by contrast, no analogous derivation is immediately obvious for the latter system.

In Section \ref{orient}, we will see that the octagon formed by the division points of a descendant chain is not only centrally symmetric but also convex. From the point of view of this octagon, each deformation has the same effect as some affine transformation. Thus, by the same argument as above, the division points of $W$ must form a convex octagon, too. Our question regarding central symmetry can now be posed for convexity as well.

\subsection{Partial Vectors} \label{pv}

We define the \emph{partial vectors} of $U$ to be the spans of the parts of $U$. We denote them by $\bfu_i = \sigma(U_i) = P_i \tto P_{i + 1}$. Then $\bfu_{i + 4} = -\bfu_i$ for all $i$, by Lemma \ref{symml}.

Denote the partial vectors of $V$ by $\bfv_1$, $\bfv_2$, $\ldots$, $\bfv_8$, with $\bfv_i = \sigma(V_i) = Q_i \tto Q_{i + 1}$ for all $i$. We proceed to work out the recurrence relations which convert the partial vectors of $U$ into the partial vectors of $V$, for each $\varphi \in \psimain$. The calculations are routine, and so we present only the end results.

For $\varphi = f_i$, we get that \[\begin{aligned} \bfv_j &= \bfu_j & &\text{for all $j \not \equiv i \mymod 4$}\\ \bfv_j &= \bfu_{j - 1} + \bfu_j + \bfu_{j + 1} & &\text{for all $j \equiv i \mymod 4$}. \end{aligned}\]

The two $g^\star$-lifts behave identically to one another -- not just analogously, but literally in the exact same way. Furthermore, their behaviour is uniform over all eight partial vectors. When $\varphi \in \{\gso, \gse\}$, we get that \[\bfv_i = -(\bfu_{i - 1} + \bfu_i + \bfu_{i + 1}) \qquad \text{for all $i$}.\]

\subsection{Neighbourhood Vectors} \label{nv}

We define the \emph{neighbourhood vectors} of $U$ by $\bfs_i = \bfu_i + \bfu_{i + 1}$ for all $i$. Or, equivalently, $\bfs_i = \sigma(U_{i : (i + 1)}) = P_i \tto P_{i + 2}$. Then $\bfs_{i + 4} = -\bfs_i$ for all $i$ because $U_{(i + 4) : (i + 5)} \sim U_{i : (i + 1)}^{-1}$. (Or, alternatively, by Lemma \ref{symml}.)

Observe also that \[U_iU_{i + 1} + \bfs_{i + 2} = U_{i + 5}^{-1}U_{i + 4}^{-1} \qquad \text{for all $i$}.\]

This is the geometric analogue of (\textbf{I}). (Of course, in the original (\textbf{I}) the symbol $U$ denotes a word chain, whereas here it denotes a chain of literal polylines.)

Our choice of term for the neighbourhood vectors is due to the following reason: Suppose that $U$ is associated with a $4$-neigh\-bour double tile $\Phi$, and consider the two tilings of $\Phi$. The displacements from each tile to its four neighbours in these two tilings will then coincide with the eight neighbourhood vectors of $U$. Specifically, in one tiling the displacements will be $\bfs_1$, $\bfs_3$, $\bfs_5$, $\bfs_7$, with the tiling itself using the lattice $\mathcal{L}(\bfs_1, \bfs_3)$; while in the other tiling they will be $\bfs_2$, $\bfs_4$, $\bfs_6$, $\bfs_8$, with the tiling itself using the lattice $\mathcal{L}(\bfs_2, \bfs_4)$.

Denote the neighbourhood vectors of $V$ by $\bft_1$, $\bft_2$, $\ldots$, $\bft_8$, with $\bft_i = \bfv_i + \bfv_{i + 1} = \sigma(V_{i : (i + 1)}) = Q_i \tto Q_{i + 2}$. Just as in the previous sub-section, for each $\varphi \in \psimain$ we proceed to work out the recurrence relations which convert the neighbourhood vectors of $U$ into the neighbourhood vectors of $V$.

For $\varphi = f_i$, we get that \[\begin{aligned} \bft_j &= \bfs_j + \bfs_{j + 1} - \bfs_{j + 2} & &\text{for all $j \equiv i \mymod 4$}\\ \bft_j &= \bfs_j + \bfs_{j - 1} - \bfs_{j - 2} & &\text{for all $j \equiv i + 3 \mymod 4$}\\ \bft_j &= \bfs_j & &\text{otherwise}. \end{aligned}\]

Once again, the two $g^\star$-lifts behave identically to one another. Furthermore, their behaviour over the neighbourhood vectors is exactly the same as their behaviour over the partial vectors. When $\varphi \in \{\gso, \gse\}$, we get that \[\bft_i = -(\bfs_{i - 1} + \bfs_i + \bfs_{i + 1}) \qquad \text{for all $i$}.\]

\subsection{Orientations} \label{orient}

For all $i$ and $j$, let $u_{ij} = \bfu_i \times \bfu_j$ and $s_{ij} = \bfs_i \times \bfs_j$. Define also $v_{ij}$ and $t_{ij}$ similarly for $V$.

\begin{lemma} \label{puij} Every descendant chain $U$ satisfies $u_{ij} > 0$ for all $1 \le i < j \le 4$. \end{lemma} 

Notice that Lemma \ref{puij} in fact determines the sign of $u_{ij}$ for all $i$ and $j$. Geometrically, it tells us that the partial vectors $\bfu_1$, $\bfu_2$, $\ldots$, $\bfu_8$ occur in this cyclic order as we make one full counterclockwise revolution around the origin. Or, equivalently, it tells us that the division points $P_1$, $P_2$, $\ldots$, $P_8$ of $U$ are, in this order, the vertices of a convex octagon oriented counterclockwise. (From Section \ref{symm}, we already know that this octagon is centrally symmetric as well.)

\begin{myproof} By induction on the descent of $U$. The base case is clear. For the induction step, assume that the claim holds for $U$. We must check that it holds for $V$, too.

Suppose first that $\varphi = f_1$. Then \begin{gather*} v_{12} = (\bfu_8 + \bfu_1 + \bfu_2) \times \bfu_2 = u_{12} + u_{24}\\ v_{13} = (\bfu_8 + \bfu_1 + \bfu_2) \times \bfu_3 = u_{13} + u_{23} + u_{34}\\ v_{14} = (\bfu_8 + \bfu_1 + \bfu_2) \times \bfu_4 = u_{14} + u_{24}\\ v_{23} = u_{23} \qquad v_{24} = u_{24} \qquad v_{34} = u_{34}. \end{gather*}

The other three $f$-lifts behave similarly.

Suppose, now, that $\varphi \in \{\gso, \gse\}$. Since the same formula $\bfv_i = -(\bfu_{i - 1} + \bfu_i + \bfu_{i + 1})$ applies to all indices $i$, it suffices to verify the positivity of just $v_{12}$ and $v_{13}$. We find that \begin{align*} v_{12} &= (\bfu_8 + \bfu_1 + \bfu_2) \times (\bfu_1 + \bfu_2 + \bfu_3)\\ &= u_{13} + u_{14} + u_{23} + u_{24} + u_{34}\\ v_{13} &= (\bfu_8 + \bfu_1 + \bfu_2) \times (\bfu_2 + \bfu_3 + \bfu_4)\\ &= u_{12} + u_{13} + u_{14} + u_{23} + 2u_{24} + u_{34}, \end{align*} and the induction step is complete. \end{myproof}

\begin{lemma} \label{psij} Every descendant chain $U$ satisfies $s_{ij} > 0$ for all $1 \le i < j \le 4$. \end{lemma} 

Once again, Lemma \ref{psij} determines the sign of $s_{ij}$ for all $i$ and $j$; and, geometrically, it tells us that the vectors $\bfs_1$, $\bfs_2$, $\ldots$, $\bfs_8$ occur in this cyclic order as we make one full counterclockwise revolution around the origin.

\begin{myproof} Since the same formula $\bfs_i = \bfu_i + \bfu_{i + 1}$ applies to all indices $i$, it suffices to verify the positivity of just $s_{12}$ and $s_{13}$. We find that \begin{align*} s_{12} &= (\bfu_1 + \bfu_2) \times (\bfu_2 + \bfu_3) = u_{12} + u_{13} + u_{23}\\ s_{13} &= (\bfu_1 + \bfu_2) \times (\bfu_3 + \bfu_4) = u_{13} + u_{14} + u_{23} + u_{24}, \end{align*} and now the result follows by Lemma \ref{puij}. \end{myproof}

It is interesting to note that Lemma \ref{psij} does not seem to admit a direct proof by induction on the descent of $U$. In contrast with Lemma \ref{puij}, the induction step does not quite go through for the $f$-lifts.

\begin{lemma} \label{odc} Every descendant chain $U$ satisfies $\area(U) > 0$. Thus, in particular, if $U$ is free of self-intersections, it must be oriented counterclockwise. \end{lemma} 

\begin{myproof} By Lemma \ref{area}, we get that $\area(U) = s_{13} = s_{24}$. However, both of $s_{13}$ and $s_{24}$ are positive by Lemma \ref{psij}. \end{myproof}

This observation will be important for the last step of the proof of Lemma \ref{isg} in Section \ref{good}.

\subsection{The Geometry of the Lifts} \label{gl}

We go on to derive the geometric forms of the construction rules for our lifting transformations. Throughout this sub-section, we assume that both of $U$ and $V$ are centered at the origin. (The center of a descendant chain is defined in Section \ref{symm}.)

For each $\varphi \in \psimain$, already we know that the parts of $V$ are translation copies of certain slices of $U$ and their reversals. We proceed now to work out the corresponding displacements explicitly. Specifically, here by ``explicitly'' we mean ``in terms of the neighbourhood vectors of~$U$''. Once again, the calculations are routine, and we present only the end results.

For the $f$-lifts, let \[\begin{aligned} \bfeta_i &= \half(\bfs_{i + 1} + \bfs_{i + 2})\\ \bfomega_i &= \half(\bfs_{i + 1} - \bfs_{i + 2}) \end{aligned} \qquad \text{for all $i = 1$, $2$, $3$, $4$}.\]

Then with $\varphi = f_i$ we get that \[\begin{gathered} V_i = U_{(i + 3) : (i + 5)}^{-1} - \bfeta_i\\ V_j = U_j + \bfomega_i \qquad \text{for all $j \in \{i + 1, i + 2, i + 3\}$}\\ V_{i + 4} = U_{(i - 1) : (i + 1)}^{-1} + \bfeta_i\\ V_j = U_j - \bfomega_i \qquad \text{for all $j \in \{i + 5, i + 6, i + 7\}$}. \end{gathered} \tag{\textbf{F}$_i$}\]

Continuing on to the $g^\star$-lifts, with $\varphi = \gso$ we get that \[\begin{aligned} V_i &= U_{(i - 2) : (i + 2)}^{-1} + \half(\bfs_i + \bfs_{i + 1} + \bfs_{i + 2} + \bfs_{i + 3}) & &\text{for all odd $i$}\\ V_i &= U_{(i + 3) : (i + 5)} + \half(\bfs_i + \bfs_{i + 1} + \bfs_{i + 2} + \bfs_{i + 3}) & &\text{for all even $i$}. \end{aligned} \tag{\textbf{G}$_\text{\rm Odd}$}\]

The construction rules for $\gse$ are analogous. The only difference is that the words ``odd'' and ``even'' are swapped. We denote the corresponding system of equations by (\textbf{G}$_\text{\rm Even}$).

\section{The Positive Part of the Description} \label{pos}

Our goal in this section will be to show that every descendant chain is free of self-intersections. The positive part of our main theorem will then follow right away.

The proof will be by induction on the descent of the chain. However, our induction hypothesis will involve not only the chain itself, but also some of its closest descendants. Furthermore, in addition to all of these chains being free of self-intersections, our induction hypothesis will include also certain geometric relationships between them.

In Section \ref{good}, we introduce the key notion of ``goodness''. Sections \ref{dogood}--\ref{pgood} explore this notion in further detail and generalise it. Then in Section \ref{full} we state and discuss our strengthened induction hypothesis. Finally, in Sections \ref{stepf} and \ref{stepg} we carry out the induction step.

For most of the proof, we will be reasoning about the geometric properties of the $f$-lifts and their compositions. The $g^\star$-lifts will only come into play at the very end, in Section \ref{stepg}.

Our definitions and lemmas will be phrased in terms of an arbitrary descendant chain $U$. However, they will in fact apply more generally to all interleaved double chains $U$ as in Lemmas \ref{symml} and \ref{puij}--\ref{odc}.

Throughout the rest of the present section, let $U = U_1 : U_2 : \cdots : U_8$ be a descendant chain of literal polylines, centered at the origin, with division points $P_1$, $P_2$, $\ldots$, $P_8$ and neighbourhood vectors $\bfs_1$, $\bfs_2$, $\ldots$, $\bfs_8$. Let also $V = \varphi(U) = V_1 : V_2 : \cdots : V_8$, with $\varphi \in \psimain$, be a lift of $U$, once again centered at the origin, with division points $Q_1$, $Q_2$, $\ldots$, $Q_8$ and neighbourhood vectors $\bft_1$, $\bft_2$, $\ldots$, $\bft_8$. Or, in short, both of $U$ and $V$ are as in Section \ref{gl}.

For convenience, we introduce the following convention: Whenever a translation copy of $U$ is denoted with a superscript, its parts and division points will be denoted with the same superscript; and similarly for $V$ as well. For example, if $V^\diamond$ is a translation copy of $V$, then its parts will be $V^\diamond_1$, $V^\diamond_2$, $\ldots$, $V^\diamond_8$ and its division points will be $Q^\diamond_1$, $Q^\diamond_2$, $\ldots$, $Q^\diamond_8$.

\subsection{Goodness} \label{good}

For the sake of clarity, in this sub-section we focus on the $f$-lift $f_1$. The analogous results for the other three $f$-lifts will be implied.

Let $V = f_1(U)$. Recall the definition of the vectors $\bfeta_1 = \half(\bfs_2 + \bfs_3)$ and $\bfomega_1 = \half(\bfs_2 - \bfs_3)$ from Section \ref{gl}. Then let also \begin{align*} U^- &= U - \bfomega_1\\ U^+ &= U + \bfomega_1. \end{align*}

Notice that there is substantial overlap between $U^-$ and $U^+$, on the one hand, and $V$, on the other. Indeed, by (\textbf{F}$_1$) we get that $V_i = U^+_i$ for all $i \in \{2, 3, 4\}$ and $V_i = U^-_i$ for all $i \in \{6, 7, 8\}$.

Furthermore, $V_1 = U^-_{1 : 2} \cup U^+_{8 : 1}$ and $V_5 = U^+_{5 : 6} \cup U^-_{4 : 5}$. The former identity follows because, by (\textbf{F}$_1$), an $\bfeta_1$ translation maps $U^-_{1 : 2}$, $U^+_{8 : 1}$, $V_1$ onto $U_{5 : 6}^{-1}$, $U_{4 : 5}^{-1}$, $U_{4 : 6}^{-1}$, respectively; and the latter identity is derived analogously.

We conclude that $V = U^-_{4 : 2} \cup U^+_{8 : 6}$ and both slices on the right-hand side are in fact sub-polylines of $V$. So, in particular, $U^-$ and $U^+$ together cover $V$.

This leaves only the parts $U^-_3$ and $U^+_7$ of $U^-$ and $U^+$ unaccounted for. Experimentally, in many small cases they are contained within $[V]$. (E.g., for all $U$ and $V = f_1(U)$ in Figure \ref{tree}.) Furthermore, in many small cases they touch $V$ only at their endpoints. These observations can be expressed succinctly in the notation of Section \ref{line}. Thus we arrive at the next definition:

\begin{figure}[ht] \centering \includegraphics{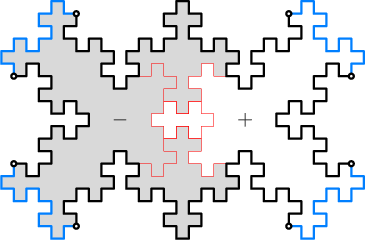} \caption{} \label{fgood} \end{figure}

We say that $U$ is $f_1$-\emph{good} when the following conditions are satisfied: (i) $U$ is free of self-intersections; (ii) $f_1(U)$ is free of self-intersections; (iii) $\lopen U^-_3 \ropen \subseteq \langle V \rangle$; (iv) $\lopen U^+_7 \ropen \subseteq \langle V \rangle$.

One example is shown in Figure \ref{fgood}. In it, $U$ is the descendant chain of Figure \ref{chain-ge-go}. The polylines $U^-_3$ and $U^+_7$ are outlined in red, and the translation copy $U^-$ of $U$ has been highlighted.

Notice that conditions (iii) and (iv) are equivalent, by central symmetry. So, in order to verify goodness, it suffices to establish just one of them.

It is clear what the analogous definitions would be for the other $f$-lifts. For $f_i$, we let $U^- = U - \bfomega_i$ and $U^+ = U + \bfomega_i$. Then we consider their parts $U^-_{i + 2}$ and $U^+_{i - 2}$.

There are some direct corollaries of goodness which will be important to us. Suppose that $U$ is $f_1$-good. Then, by the preceding discussion, $U^- \subseteq [V]$ and $U^+ \subseteq [V]$. Hence, also $\langle U^- \rangle \subseteq \langle V \rangle$ and $\langle U^+ \rangle \subseteq \langle V \rangle$. We assume all of these properties to be implied whenever we stipulate that some descendant chain is $f_1$-good.

How do we establish goodness? Before we answer this question, we make some preliminary observations.

Suppose that $U$ is free of self-intersections and consider two translation copies $U'$ and $U''$ of $U$, related by a translation $\bftau = U' \tto U''$ with $\bftau \in \mathcal{L}(\bfs_1, \bfs_3)$. Then $U'$ and $U''$ can be embedded into the tiling of the plane with copies of $[U]$ which uses the lattice $\mathcal{L}(\bfs_1, \bfs_3)$.

The combinatorial structure of the tiling now allows us to deduce exactly how $U'$ and $U''$ intersect. When $\bftau = \bfs_1$, we get that $U' \cap U'' = U'_{3 : 4} \reveq U''_{7 : 8}$; and the remaining cases with $\bftau \in \{\pm \bfs_1, \pm \bfs_3\}$ are similar. When $\bftau = \bfs_1 + \bfs_3$, we get that $U' \cap U'' = P'_5 = P''_1$; and the remaining cases with $\bftau \in \{-\bfs_1, \bfs_1\} + \{-\bfs_3, \bfs_3\}$ are similar. Finally, when $\bftau = \alpha \bfs_1 + \beta \bfs_3$ with $\alpha^2 + \beta^2 \ge 4$, we get that $U'$ and $U''$ are disjoint.

The scenario where $\bftau$ is drawn out of the lattice $\mathcal{L}(\bfs_2, \bfs_4)$ can be analysed analogously. We will make use of these observations throughout the present section.

We are ready now to introduce our main tool for establishing goodness:

\begin{lemma} \label{isg} Suppose that $U$ is free of self-intersections. Suppose also that $U^-$ is disjoint from $U^+_{2 : 4}$. (By central symmetry, this condition is equivalent to $U^+$ being disjoint from $U^-_{6 : 8}$.) Then $U$ is $f_1$-good. \end{lemma} 

Notice that the converse implication holds as well. I.e., if $U$ is $f_1$-good, then $U^-$ is disjoint from $U^+_{2 : 4}$ and $U^+$ is disjoint from $U^-_{6 : 8}$.

\begin{myproof} Let \begin{align*} U^\text{I} &= U - \bfeta_1\\ U^\text{II} &= U + \bfeta_1. \end{align*}

By (\textbf{F}$_1$), we get that $V_1 \reveq U^\text{I}_{4 : 6}$ and $V_5 \reveq U^\text{II}_{8 : 2}$.

We proceed to check conditions (ii) and (iii) of the definition of goodness, one by one.

First we show that $V$ is free of self-intersections. Partition $V$ into $\lopen V_1 \ropen$, $V_{2 : 4}$, $\lopen V_5 \ropen$, $V_{6 : 8}$. Each one of these slices of $V$ is a translation copy, modulo reversal, of some slice of $U$. So none of them can contain self-intersections when considered in isolation. We are left to check that no unwanted intersections can occur between them, either.

Since $U^\text{I} + \bfs_3 = U^-$, these two copies of $U$ can be embedded into the tiling of $[U]$ which uses the lattice $\mathcal{L}(\bfs_1, \bfs_3)$. We get that $\lopen U^\text{I}_{4 : 6} \ropen$ is disjoint from $U^-_{4 : 8}$. Similarly, $\lopen U^\text{I}_{4 : 6} \ropen$ is disjoint from $U^+_{2 : 6}$, too. However, $U^\text{I}_{4 : 6} \reveq V_1$ whereas $U^-_{4 : 8}$ and $U^+_{2 : 6}$ together cover $V_{2 : 8}$. Hence, $\lopen V_1 \ropen$ is disjoint from $V_{2 : 8}$.

By central symmetry, also $\lopen V_5 \ropen$ is disjoint from $V_{6 : 4}$.

What remains is to show that $V_{2 : 4}$ and $V_{6 : 8}$ are disjoint. But this follows because $V_{2 : 4} = U^+_{2 : 4}$, $V_{6 : 8} \subseteq U^-$, and the fact that $U^-$ is disjoint from $U^+_{2 : 4}$ has been given to us for free.

We move on to condition (iii). First we will verify that $\lopen U^-_3 \ropen$ is disjoint from $V$.

This is clear for $V_{6 : 8} = U^-_{6 : 8}$, as $U^-$ is free of self-intersections. The embedding of $U^-$ and $U^\text{I}$ into the tiling of $[U]$ which uses the lattice $\mathcal{L}(\bfs_1, \bfs_3)$ tells us that $\lopen U^-_3 \ropen$ is disjoint also from $U^\text{I}_{4 : 6} \reveq V_1$, and a similar argument establishes its disjointness from $V_5$. Finally, the entirety of $U^-$ being disjoint from $V_{2 : 4} = U^+_{2 : 4}$ has been, as above, granted to us for free.

So $\lopen U^-_3 \ropen$ and $V$ are disjoint. It follows that either $\lopen U^-_3 \ropen \subseteq \langle V \rangle$, or else $\lopen U^-_3 \ropen$ and $[V]$ are disjoint. We must rule out the latter scenario.

Consider a small vicinity $\mathcal{D}$ of $P^-_3 = P^\text{I}_5$. Notice that this point lies on $\lopen U^\text{I}_{4 : 6} \ropen \reveq \lopen V_1 \ropen$, and consider the short initial section of $\lopen U^-_3 \ropen$ as well as the short interior section of $\lopen V_1 \ropen$ contained within $\mathcal{D}$.

Both of $U$ and $V$ are oriented counterclockwise by Lemma \ref{odc}. Since $V$ is oriented counterclockwise, inside of $\mathcal{D}$ we get that $\langle V \rangle$ lies ``on the left'' of $V_1$. On the other hand, since $U$ is oriented counterclockwise, too, the embedding of $U^-$ and $U^\text{I}$ into the tiling of $[U]$ which uses the lattice $\mathcal{L}(\bfs_1, \bfs_3)$ tells us that inside of $\mathcal{D}$ a short initial section of $\lopen U^-_3 \ropen$ also lies ``on the left'' of~$V_1$, by virtue of $U^\text{I}_{4 : 6} \reveq V_1$. We conclude that $\lopen U^-_3 \ropen$ and $\langle V \rangle$ intersect, and condition (iii) follows. \end{myproof}

\subsection{What Good Does Goodness Do?} \label{dogood}

We say that the vector $\bfmu$ is a \emph{separator} for $U$ when $U$ and $U + \bfmu$ are disjoint. We care about separators because they can be helpful for ruling out self-intersections in the lifts of $U$.

There are some separators that we get pretty much for free. As we observed in the previous sub-section, from the combinatorial structure of the tiling of $[U]$ using the lattice $\mathcal{L}(\bfs_1, \bfs_3)$ we can deduce that $\alpha\bfs_1 + \beta\bfs_3$ is a separator for $U$ whenever the integers $\alpha$ and $\beta$ satisfy $\alpha^2 + \beta^2 \ge 4$. Similarly also for the lattice $\mathcal{L}(\bfs_2, \bfs_4)$.

The main significance of goodness is that it allows us to find many other separators for~$U$. Indeed, suppose that $U$ is $f_1$-good with $V = f_1(U)$. Let also $\mathfrak{T}$ be the tiling of the plane with copies of $[V]$ which uses, say, the lattice $\mathcal{L}(\bft_1, \bft_3)$. The analysis for the lattice $\mathcal{L}(\bft_2, \bft_4)$ is analogous.

Let $\bfrho = \alpha\bft_1 + \beta\bft_3$ be a nonzero vector of $\mathcal{L}(\bft_1, \bft_3)$, and set $V^\ddag = V + \bfrho$. We can assume, without loss of generality, that both of $[V]$ and $[V^\ddag]$ are tiles in $\mathfrak{T}$.

Define $U^-$ and $U^+$ as in the previous sub-section, and consider also $U' = U^- + \bfrho$ and $U'' = U^+ + \bfrho$. Then both of $U^-$ and $U^+$ are contained within $[V]$; and, by translation, also both of $U'$ and $U''$ are contained within $[V^\ddag]$.

\begin{figure}[ht] \centering \includegraphics{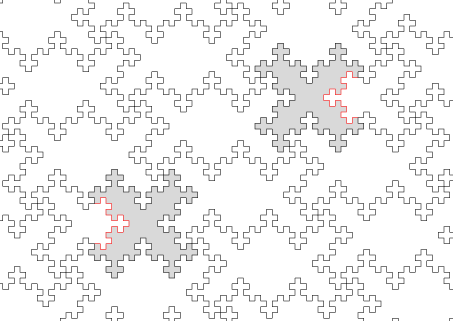} \caption{} \label{sep} \end{figure}

If $[V]$ and $[V^\ddag]$ are disjoint, then certainly both of $U^-$ and $U^+$ will be disjoint from both of $U'$ and $U''$, too. But, by the preceding discussion, $[V]$ and $[V^\ddag]$ are in fact disjoint whenever $\alpha^2 + \beta^2 \ge 4$. For example, Figure \ref{sep} shows $\alpha = 2$ and $\beta = 1$, with $U$ and $V$ as in Figure \ref{fgood}.

Even if $[V]$ and $[V^\ddag]$ are not disjoint, they can only touch along their boundaries. However, in the previous sub-section we developed a detailed understanding of how $U^-$ and $U^+$ overlap with the boundary $V$ of $[V]$. Naturally, this understanding translates to $U'$, $U''$, and $V^\ddag$ as well. So we might be able to show that some of $U^-$ and $U^+$ are disjoint from some of $U'$ and $U''$ also in some cases when $[V]$ and $[V^\ddag]$ intersect.

We arrive at the following lemma:

\begin{lemma} \label{sepg} Suppose that $U$ is $f_i$-good for some $f$-lift $f_i$, with $\bfrho_1$, $\bfrho_2$, $\bfrho_3$, $\bfrho_4$ being the neighbourhood vectors of $f_i(U)$. Let $\bfrho$ be any nonzero vector in one of the lattices $\mathcal{L}(\bfrho_1, \bfrho_3)$ and $\mathcal{L}(\bfrho_2, \bfrho_4)$. Let also $\bfdelta$ be any vector out of $\{\mathbf{0}, \pm 2\bfomega_i\}$. Then $\bfrho + \bfdelta$ is a separator for $U$, with a finite number of exceptions. Furthermore, no exceptions occur when $\bfrho$ is of the form either $\alpha\bfrho_1 + \beta\bfrho_3$ or $\alpha\bfrho_2 + \beta\bfrho_4$ with $\alpha^2 + \beta^2 \ge 4$. \end{lemma} 

We proceed now to list the exceptions. We are only going to spell out the details of the analysis for $f_1$ and the lattice $\mathcal{L}(\bft_1, \bft_3)$. The remaining cases can all be handled analogously, and we assume the corresponding results to be implied.

Since the negation of an exception must be an exception as well, it suffices to consider just one vector out of each pair of mutual negations.

Suppose first that $\bfrho = \bft_1$. Then $[V] \cap [V^\ddag] = V_{3 : 4} \reveq V^\ddag_{7 : 8}$. From the previous sub-section, we know that $V_{3 : 4}$ and $U^-$ are disjoint but $V_{3 : 4} \subseteq U^+$. Similarly, $V^\ddag_{7 : 8} \subseteq U'$ but $V^\ddag_{7 : 8}$ and $U''$ are disjoint. So, in this case, $U^+$ and $U'$ do intersect -- but the remaining three pairs out of $\{U^-, U^+\} \times \{U', U''\}$ are disjoint.

Similar arguments resolve also the two cases when $\bfrho = \bft_1 \pm \bft_3$.

Suppose, finally, that $\bfrho = \bft_3$. In this case, there are no pairs out of $\{U^-, U^+\} \times \{U', U''\}$ for which disjointness can be guaranteed. One counterexample is given by $U = f_1(\maltese)$ and $V = f_1 \circ f_1(\maltese)$, when all four pairs intersect.

(Notice, though, that with $\bfrho = \bft_3$ not all of the intersections are forced. For example, if $U = \maltese$ and $V = f_1(\maltese)$, then $U^+$ and $U'$ become disjoint. Hence, in this case, it would be more appropriate to talk about ``possible exceptions'' rather than just ``exceptions''.)

\begin{table}[ht] \centering \begin{tabular}{|c|c|}
$\bfrho$ & $\bfdelta$\\
\hline
$\{\bft_1, \bft_1 \pm \bft_3, -\bft_4, \pm \bft_2 - \bft_4\}$ & $-2\bfomega_1$\\
\hline
$\{\bft_2, \bft_3\}$ & $\{\mathbf{0}, \pm 2\bfomega_1\}$\\
\hline
\end{tabular} \caption{} \label{exft} \end{table}

The conclusions from our analysis can be summarised as follows:

\begin{lemma} \label{exf} The possible exceptions of Lemma \ref{sepg} for the $f_1$-lift are the ones listed in Table~\ref{exft}, together with their negations. \end{lemma} 

\subsection{Generalised Goodness} \label{ggood}

Here, we generalise the notion of goodness introduced in Section \ref{good}.

Let $\frf = (f_{I(1)}, f_{I(2)}, \ldots, f_{I(k)})$ be a $k$-element sequence of $f$-lifts. Define $W^{(1)} = U$ and $W^{(i + 1)} = f_{I(i)}(W^{(i)})$ for all $1 \le i \le k$, as well as $W = W^{(k + 1)}$.

We say that $U$ is $\frf$-\emph{good} when $W^{(i)}$ is $f_{I(i)}$-good for all $1 \le i \le k$.

So, in particular, with $k = 1$ we recover our original notion of goodness. We already know that in this case certain translation copies of $U$ fit inside of $[W]$. We proceed to extend this observation to generalised goodness. Suppose, from now on, that $U$ is in fact $\frf$-good.

For each $1 \le i \le k + 1$, let $\bftau^{(i)}_1$, $\bftau^{(i)}_2$, $\ldots$, $\bftau^{(i)}_8$ be the neighbourhood vectors of $W^{(i)}$. Then, for each $1 \le j \le 4$, let also \[\bfomega^{(i)}_j = \half(\bftau^{(i)}_{j + 1} - \bftau^{(i)}_{j + 2}).\]

(These vectors are defined relative to $W^{(i)}$ in the same way as the vectors $\bfomega_1$, $\bfomega_2$, $\bfomega_3$, $\bfomega_4$ are defined relative to $U$.)

For each $1 \le i \le k + 1$, we are going to define a set of vectors $\Delta_i$ such that, for all $\bfdelta \in \Delta_i$, the translation copy $U + \bfdelta$ of $U$ fits inside of $[W^{(i)}]$. The definition will be by induction on $i$. For the base case, $i = 1$, we let $\Delta_1 = \{\mathbf{0}\}$.

Suppose, for the induction step, that $\Delta_i$ has been defined already. Then $U + \bfdelta$ fits inside of $[W^{(i)}]$ for all $\bfdelta \in \Delta_i$. On the other hand, since $W^{(i)}$ is $f_{I(i)}$-good, also both of $W^{(i)} - \bfomega^{(i)}_{I(i)}$ and $W^{(i)} + \bfomega^{(i)}_{I(i)}$ fit inside of $[W^{(i + 1)}]$. Hence, we define \[\Delta_{i + 1} = \Delta_i + \{-\bfomega^{(i)}_{I(i)}, \bfomega^{(i)}_{I(i)}\},\] and we conclude that $U + \bfdelta$ fits inside of $[W^{(i + 1)}]$ for all $\bfdelta \in \Delta_{i + 1}$.

Denote $\Delta = \Delta_{k + 1}$. We will sometimes need to reference individual elements of $\Delta$. So, for convenience, in a context where $\frf$ has been fixed, let us agree to enumerate $\Delta$ in the following manner: Given a vector $\bfdelta \in \Delta$, first write $\bfdelta = \pm \bfomega^{(1)}_{I(1)} \pm \bfomega^{(2)}_{I(2)} \pm \cdots \pm \bfomega^{(k)}_{I(k)}$. There are $k$ many $\pm$ signs on the right-hand side. Turn each minus sign into a zero and each plus sign into a one, and then read off the resulting binary number $r$. (Of course, $0 \le r \le 2^k - 1$.) We assign to $\bfdelta$ the index $r$, and we denote it by $\bfdelta_r$.

For example, when $\frf = (f_1, f_3)$, we get that $\bfomega^{(1)}_1 = \bfomega_1$ and $\bfomega^{(2)}_3 = \bfomega_3$. So the elements of $\Delta$ become $\bfdelta_0 = -\bfomega_1 - \bfomega_3$, $\bfdelta_1 = -\bfomega_1 + \bfomega_3$, $\bfdelta_2 = \bfomega_1 - \bfomega_3$, and $\bfdelta_3 = \bfomega_1 + \bfomega_3$. (This special case will be important in the next sub-section.)

Notice that this enumeration is not always injective; the same element of $\Delta$ might happen to receive several distinct indices. For example, when $\frf = (f_1, f_1, f_1)$, we get that $\bfomega^{(1)}_1 = \bfomega^{(2)}_1 = \bfomega^{(3)}_1 = \bfomega_1$, and so $\bfdelta_1 = \bfdelta_2 = \bfdelta_4 = -\bfomega_1$. This ambiguity is not going to matter much to us.

By the previous sub-section, goodness allows us to conclude that certain translation copies of $U$ are disjoint. So does generalised goodness, in the exact same manner.

Indeed, let $\bftau_i = \bftau^{(k + 1)}_i$ for all $i$. Then $[W]$ tiles the plane with both of the lattices $\mathcal{L}(\bftau_1, \bftau_3)$ and $\mathcal{L}(\bftau_2, \bftau_4)$. Let $[W']$ and $[W'']$ be two distinct tiles in one of these two tilings, and let also $U'$ and $U''$ be two translation copies of $U$ which fit inside of $[W']$ and $[W'']$, respectively. If $[W']$ and $[W'']$ are disjoint, then of course so will be $U'$ and $U''$ as well. Furthermore, even if $[W']$ and $[W'']$ touch, we might still be able to prove that $U'$ and $U''$ are disjoint anyway, provided that we have some additional information about how $U'$ and $U''$ overlap with the boundaries $W'$ and $W''$ of their corresponding tiles.

Thus we arrive at the following generalisation of Lemma \ref{sepg}:

\begin{lemma} \label{sepgg} Suppose that $U$ is $\frf$-good. Let $\bftau$ be any nonzero element of one of the lattices $\mathcal{L}(\bftau_1, \bftau_3)$ and $\mathcal{L}(\bftau_2, \bftau_4)$. Let also $\bfdelta$ be any element of $\Delta - \Delta$. Then $\bftau + \bfdelta$ is a separator for $U$, with a finite number of exceptions. Furthermore, no exceptions occur when $\bftau$ is of the form either $\alpha \bftau_1 + \beta \bftau_3$ or $\alpha \bftau_2 + \beta \bftau_4$ with $\alpha^2 + \beta^2 \ge 4$. \end{lemma}

Given a concrete $\frf$, we might be able to pin down all of the possible exceptions explicitly, by an analysis similar to the one which gave us Lemma \ref{exf}. In the next sub-section, we will indeed do so for two important special cases, both with $k = 2$. However, in all cases with $k \ge 3$, the specificity of Lemma \ref{sepgg} will turn out to be sufficient for our purposes. In fact, we are never going to need to apply Lemma \ref{sepgg} with $k \ge 4$; and an application with $k = 3$ will be required in just one single instance, in Section \ref{stepg}.

\subsection{\texorpdfstring{$\fodd$}{fOdd} and \texorpdfstring{$\feven$}{fEven}-Goodness} \label{pgood} 

There are two special cases of generalised goodness which will be particularly important to us, and now we consider them in detail.

We say that $U$ is $\fodd$-\emph{good} when it is both $(f_1, f_3)$-good and $(f_3, f_1)$-good. Similarly, we say that $U$ is $\feven$-\emph{good} when it is both $(f_2, f_4)$-good and $(f_4, f_2)$-good.

Since these two special cases are analogous to one another, throughout the rest of this sub-section we will consider $\fodd$ only, and the corresponding results for $\feven$ will be implied.

Unpacking the definitions, we see that $\fodd$-goodness involves a total of four descendant chains: $U$, $f_1(U)$, $f_3(U)$, and $\fodd(U)$. Furthermore, it entails both $f_1$-goodness and $f_3$-goodness.

Notice that, in the setting of the previous sub-section, the special cases $\frf = (f_1, f_3)$ and $\frf = (f_3, f_1)$ yield the same $W = \fodd(U)$. So they also yield the same neighbourhood vectors $\bftau_1$, $\bftau_2$, $\ldots$, $\bftau_8$. Finally, they yield the same $\Delta$, too, which we denote by $\Delta_\text{\rm Odd}$.

(Though $\frf = (f_1, f_3)$ and $\frf = (f_3, f_1)$ do not yield the same enumeration of $\Delta_\text{\rm Odd}$; the labels $\bfdelta_1$ and $\bfdelta_2$ are swapped. We will use the enumeration induced by $\frf = (f_1, f_3)$, and for $\Delta_\text{Even}$ we will use the one induced by $\frf = (f_2, f_4)$.)

Recall, from the previous sub-section, our explicit calculation of $\Delta_\text{\rm Odd}$ with $\frf = (f_1, f_3)$. It will be helpful to introduce special notation for the translation copies of $U$ displaced by the elements of this set, similarly to how we defined $U^-$ and $U^+$ in Section \ref{good}. Let \begin{align*} U^{--} &= U - \bfomega_1 - \bfomega_3 & U^{-+} &= U - \bfomega_1 + \bfomega_3\\ U^{+-} &= U + \bfomega_1 - \bfomega_3 & U^{++} &= U + \bfomega_1 + \bfomega_3. \end{align*}

\begin{figure}[ht] \centering \includegraphics{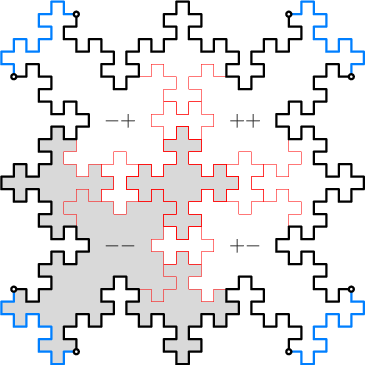} \caption{} \label{ffgood} \end{figure}

By (\textbf{F}$_1$) and (\textbf{F}$_3$), we get that these translation copies of $U$ together cover $W$. Specifically, $W_1 = U^{--}_{1 : 2} \cup U^{+-}_{8 : 1}$; $W_2 = U^{+-}_2$; $W_3 = U^{+-}_{3 : 4} \cup U^{++}_{2 : 3}$; $W_4 = U^{++}_4$; and similarly for the remaining parts of~$W$.

One example is shown in Figure \ref{ffgood}, once again with $U$ being the double chain of Figure~\ref{chain-ge-go}. The parts of $U^{--}$, $U^{-+}$, $U^{+-}$, $U^{++}$ contained within $\langle W \rangle$ are outlined in red, and the translation copy $U^{--}$ of $U$ has been highlighted.

Let $V = f_1(U)$. Define $U^-$ and $U^+$ as in Section \ref{good}, as well as $\Lambda = f_3(U)$ and \begin{align*} V^- &= V - \bfomega_3 & V^+ &= V + \bfomega_3\\ \Lambda^- &= \Lambda - \bfomega_1 & \Lambda^+ &= \Lambda + \bfomega_1. \end{align*}

Notice that a $-\bfomega_3$ translation maps $U^-$, $U^+$, $V$ onto $U^{--}$, $U^{+-}$, $V^-$, respectively; and, similarly, an $\bfomega_3$ translation maps the same ordered triple onto $U^{-+}$, $U^{++}$,~$V^+$.

\begin{lemma} \label{f1f3} Suppose that $U$ is both $f_1$-good and $f_3$-good. Then $U$ is $(f_1, f_3)$-good if and only if it is $(f_3, f_1)$-good. \end{lemma}

Or, in other words, if $U$ is both $f_1$-good and $f_3$-good, then $(f_1, f_3)$-goodness and $(f_3, f_1)$-goodness are equivalent for $U$, and so either one of them entails $\fodd$-goodness.

\begin{myproof} Suppose that $U$ is $(f_1, f_3)$-good. (In addition to being $f_1$-good and $f_3$-good.) We will show that it is $(f_3, f_1)$-good, too. The opposite implication can be handled analogously.

Unpacking the definitions, we must check that $\lopen \Lambda^-_3 \ropen \subseteq \langle W \rangle$. However, $\lopen \Lambda^-_3 \ropen \subseteq \lopen U^{--}_{3 : 4} \ropen \cup \lopen U^{-+}_{2 : 3} \ropen$ by (\textbf{F}$_3$). So it suffices to verify that $\lopen U^{--}_{3 : 5} \ropen \subseteq \langle W \rangle$ and $\lopen U^{-+}_{1 : 3} \ropen \subseteq \langle W \rangle$. We consider the former inclusion in detail, and the latter one can be handled analogously. Furthermore, for the analysis we partition $\lopen U^{--}_{3 : 5} \ropen$ into $\lopen U^{--}_3 \ropen$ and $(U^{--}_{4 : 5} \ropen$.

First, $\lopen U^-_3 \ropen \subseteq \langle V \rangle$ because $U$ is $f_1$-good. Translating by $-\bfomega_3$, we get that $\lopen U^{--}_3 \ropen \subseteq \langle V^- \rangle \subseteq \langle W \rangle$, where the latter inclusion follows because $V$ is $f_3$-good.

Second, $(U^-_{4 : 5} \ropen \subseteq \lopen V_5 \ropen$ by (\textbf{F}$_1$). Translating by $-\bfomega_3$, we get that $(U^{--}_{4 : 5} \ropen \subseteq \lopen V^-_5 \ropen \subseteq \langle W \rangle$, where once again the latter inclusion follows because $V$ is $f_3$-good. \end{myproof}

Next we introduce one useful tool for establishing $\fodd$-goodness. It is the $\fodd$-analogue of Lemma \ref{isg}.

\begin{lemma} \label{isfog} Suppose that $U$ is both $f_1$-good and $f_3$-good. Suppose also that $U^{--}$ is disjoint from $U^{++}_4$, and $U^{-+}$ is disjoint from $U^{+-}_2$. (By central symmetry, these conditions are equivalent to $U^{++}$ being disjoint from $U^{--}_8$, and $U^{+-}$ is disjoint from $U^{-+}_6$, respectively.) Then $U$ is $\fodd$-good. \end{lemma} 

\begin{myproof} By Lemma \ref{f1f3}, it suffices to check that $V$ is $f_3$-good. By Lemma \ref{isg}, to this end it would be enough to verify that $V^-$ is disjoint from $V^+_{4 : 6}$.

Since $V \subseteq U^- \cup U^+$, translating by $-\bfomega_3$ we get that $V^- \subseteq U^{--} \cup U^{+-}$. So it suffices to show that both of $U^{--}$ and $U^{+-}$ are disjoint from $V^+_{4 : 6}$. We consider $U^{--}$ in detail, and $U^{+-}$ can be handled analogously. Furthermore, for the analysis we partition $V^+_{4 : 6}$ into $V^+_4$, $V^+_5$, $V^+_6$.

First, the fact that $U^{--}$ is disjoint from $U^{++}_4 = V^+_4$ has been given to us free of charge.

Second, $U - \bfomega_3$ is disjoint from $U_{4 : 6} + \bfomega_3$ because $U$ is $f_3$-good. Translating by $-\bfomega_1$, we get that $U^{--}$ is disjoint from $U^{-+}_{4 : 6} \supseteq U^{-+}_6 = V^+_6$.

Third, $V_5 \subseteq U + \bfeta_1$ by (\textbf{F}$_1$). Translating by $\bfomega_3$, we get that \[V^+_5 \subseteq U^\diamond = U + \bfeta_1 + \bfomega_3.\]

The displacement between $U^{--}$ and $U^\diamond$ works out to $\bfs_1 + \bfs_2 + \bfs_4$, and we are only left to show that the latter vector is a separator for $U$. Since $U$ is $f_3$-good, this follows by Lemmas \ref{sepg} and \ref{exf} for the $f_3$-lift, with $\bfrho = \bfrho_2$ and $\bfdelta = \mathbf{0}$. \end{myproof}

\begin{table}[ht] \centering \begin{tabular}{|c|c|}
$\bftau$ & $\bfdelta$\\
\hline
$\{\bftau_1, -\bftau_4\}$ & $\{\bfdelta_0, \bfdelta_1\} - \{\bfdelta_2, \bfdelta_3\}$\\
\hline
$\{\bftau_2, \bftau_3\}$ & $\{\bfdelta_0, \bfdelta_2\} - \{\bfdelta_1, \bfdelta_3\}$\\
\hline
$\{\bftau_1 + \bftau_3, \bftau_2 - \bftau_4\}$ & $\bfdelta_0 - \bfdelta_3$\\
\hline
$\{-\bftau_1 + \bftau_3, \bftau_2 + \bftau_4\}$ & $\bfdelta_2 - \bfdelta_1$\\
\hline
\end{tabular} \caption{} \label{exfot} \end{table}

Finally, we work out the separators guaranteed by $\fodd$-goodness. Or, equivalently, we work out the exceptions of Lemma \ref{sepgg} in the special cases when $\frf = (f_1, f_3)$ and $\frf = (f_3, f_1)$. These two special cases will yield the same exceptions because, as we remarked in the beginning, they share the same $W$ and the same $\Delta$.

\begin{lemma} \label{exfo} The possible exceptions of Lemma \ref{sepgg} with $\mathfrak{f} = (f_1, f_3)$ are the ones listed in Table \ref{exfot}, together with their negations. \end{lemma} 

We omit the proof because it is analogous to the analysis of the possible exceptions for Lemma~\ref{sepg} in Section \ref{dogood}. The required information about how the translation copies of an $\fodd$-good $U$ overlap with $W$ is as follows: $U^{--}_{6 : 2} \subseteq \lopen W_{7 : 1} \ropen$; $\lopen U^{--}_{3 : 5} \ropen \subseteq \langle W \rangle$; and similarly also for $U^{-+}$, $U^{+-}$, $U^{++}$. The former inclusion follows by (\textbf{F}$_1$) and (\textbf{F}$_3$). The latter one is shown as in the proof of Lemma \ref{f1f3}.

\subsection{Full \texorpdfstring{$f$}{f}-Goodness} \label{full}

We say that $U$ is \emph{fully} $f$-\emph{good} when it is both $\fodd$-good and $\feven$-good.

Notice that full $f$-goodness entails also $f_i$-goodness for all $f$-lifts $f_i$. Furthermore, unpacking the definitions, we see that full $f$-goodness involves a total of seven descendant chains: $U$, $f_i(U)$ for all $f$-lifts $f_i$, and both of $\fodd(U)$ and $\feven(U)$.

For the rest of the proof, we will show, by induction on the descent of $U$, that every descendant chain $U$ is fully $f$-good. Then it will follow, in particular, that every descendant chain $U$ is free of self-intersections.

The base case of the induction is the Greek cross, and it is straightforward to check that it is in fact fully $f$-good. We handle the induction step for the $f$-lifts in Section \ref{stepf}, and the induction step for the $g^\star$-lifts is in Section \ref{stepg}.

\subsection{Induction Step for the \texorpdfstring{$f$}{f}-Lifts} \label{stepf}

Here, we carry out the induction step for the $f$-lifts. Clearly, all of them behave analogously to one another. So we are only going to consider $f_1$ in detail, and the corresponding results for the remaining three $f$-lifts will be implied.

Suppose, then, that $U$ is fully $f$-good. We must show that $f_1(U)$ is fully $f$-good as well.

Notice that, because we have broken symmetry by choosing $f_1$ specifically, $\fodd$ and $\feven$ are not going to behave analogously relative to $f_1(U)$. So we will need to establish the $\fodd$-goodness and $\feven$-goodness of $f_1(U)$ separately. Indeed, the proofs will turn out to be very different.

Throughout this sub-section, let $V = f_1(U)$. Define also $U^-$ and $U^+$ as in Section~\ref{good}.

\begin{lemma} \label{f1f1} Suppose that $U$ is $f_1$-good. Then $f_1(U)$ is $f_1$-good as well. \end{lemma}

\begin{myproof} Let $V^- = V - \bfomega_1$ and $V^+ = V + \bfomega_1$. By Lemma \ref{isg}, it suffices to show that $V^-$ is disjoint from $V^+_{2 : 4}$. (See Figure \ref{f1f1f}.)

\begin{figure}[ht] \centering \includegraphics{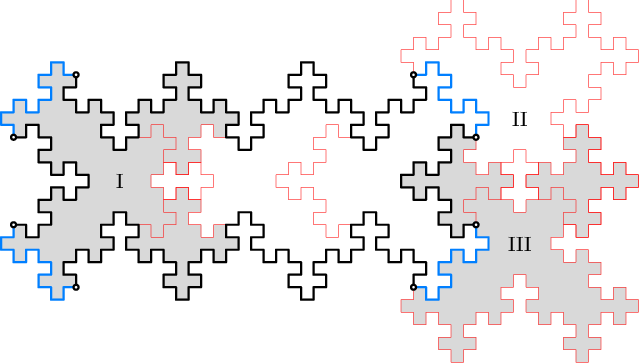} \caption{} \label{f1f1f} \end{figure} 

Since $V \subseteq U^- \cup U^+$, translating by $-\bfomega_1$ we get that $U$ and \[U^\text{I} = U - 2\bfomega_1 = U - \bfs_2 + \bfs_3\] together cover~$V^-$. So it is enough to check that both of $U$ and $U^\text{I}$ are disjoint from $V^+_{2 : 4}$.

For $U$, this is straightforward. Since $U$ is $f_1$-good, $U^-$ is disjoint from $U^+_{2 : 4} = V_{2 : 4}$. Translating by $\bfomega_1$, we get that $U$ is disjoint from $V^+_{2 : 4}$, as desired.

For $U^\text{I}$, observe that $V_{2 : 4} = U_{2 : 4} + \bfomega_1$ by (\textbf{F}$_1$), and so $V^+_{2 : 4} = U_{2 : 4} + 2\bfomega_1$. On the other hand, $U_{2 : 4} = U_{2 : 3} \cup U_{3 : 4} = (U_{6 : 7}^{-1} - \bfs_4) \cup (U_{7 : 8}^{-1} + \bfs_1) \subseteq (U - \bfs_4) \cup (U + \bfs_1)$. Translating by $2\bfomega_1$, we get that \begin{align*} U^\text{II} &= U + \bfs_1 + \bfs_2 - \bfs_3\\ U^\text{III} &= U + \bfs_2 - \bfs_3 - \bfs_4 \end{align*} together cover $V^+_{2 : 4}$.

What remains is to show that $U^\text{I}$ is disjoint from $U^\text{II}$ and $U^\text{III}$. The corresponding displacements are $U^\text{I} \tto U^\text{II} = \bfs_1 + 2\bfs_2 - 2\bfs_3$ and $U^\text{I} \tto U^\text{III} = 2\bfs_2 - 2\bfs_3 - \bfs_4$. Since $U$ is $f_1$-good, both of these vectors are indeed separators for $U$ by Lemmas \ref{sepg} and \ref{exf} for the $f_1$-lift, with $\bfrho \in \{\bfrho_1, -\bfrho_4\}$ and $\bfdelta = 2\bfomega_1$. \end{myproof}

\begin{lemma} \label{f1f2} Suppose that $U$ is both $f_1$-good and $f_2$-good. Then $f_1(U)$ is $f_2$-good as well. \end{lemma}

Similarly, if $U$ is both $f_1$-good and $f_4$-good, then $f_1(U)$ is $f_4$-good. Notice that this is not a cyclically shifted variant of Lemma \ref{f1f2}; rather, a cyclic shift would tell us instead that, if $U$ is both $f_1$-good and $f_4$-good, then $f_4(U)$ is $f_1$-good. But the proof is analogous anyway.

\begin{myproof} Setting $\frf = (f_1, f_2)$, define $\bfomega^{(2)}_j$ for $1 \le j \le 4$ as in Section \ref{ggood}. Let also $V^- = V - \bfomega^{(2)}_2$ and $V^+ = V + \bfomega^{(2)}_2$. By Lemma \ref{isg}, it suffices to show that $V^-$ is disjoint from $V^+_{3 : 5}$. (See Figure~\ref{f1f2f}.)

\begin{figure}[ht] \centering \includegraphics{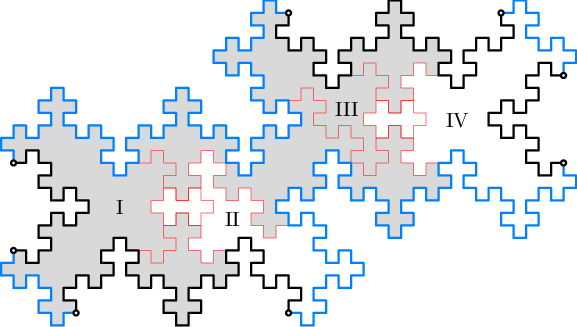} \caption{} \label{f1f2f} \end{figure} 

For convenience, denote \begin{align*} U^\text{I} &= U^- - \bfomega^{(2)}_2 & U^\text{II} &= U^+ - \bfomega^{(2)}_2\\ U^\text{III} &= U^- + \bfomega^{(2)}_2 & U^\text{IV} &= U^+ + \bfomega^{(2)}_2. \end{align*}

Since $V \subseteq U^- \cup U^+$, translating by $-\bfomega^{(2)}_2$ we get that $V^- \subseteq U^\text{I} \cup U^\text{II}$. Similarly, since $V_{3 : 5} = V_{3 : 4} \cup V_5 \subseteq U^+_{3 : 4} \cup (U^+_{5 : 6} \cup U^-_{4 : 5})$, translating by $\bfomega^{(2)}_2$ we get that $V^+_{3 : 5} \subseteq U^\text{III}_{4 : 5} \cup U^\text{IV}_{3 : 6}$.

So it is enough to verify that $U^\text{I}$ and $U^\text{II}$ are disjoint from $U^\text{III}_{4 : 5}$ and $U^\text{IV}_{3 : 6}$. These are four separate claims, and we consider them one by one.

First, $U^\text{I}$ and $U^\text{III}_{4 : 5}$. Since $U^\text{I} \tto U^\text{III} = 2\bfomega^{(2)}_2 = \bfs_2 - \bfs_4$, we get that $U^\text{I}$ and $U^\text{III}$ can be embedded into the tiling of $[U]$ which uses the lattice $\mathcal{L}(\bfs_2, \bfs_4)$. Hence, their intersection is the point $P^\text{I}_4 = P^\text{III}_8$, which does not belong to $U^\text{III}_{4 : 5}$.

Second, $U^\text{II}$ and $U^\text{IV}_{3 : 6}$. By the same reasoning as in the previous paragraph, the intersection of $U^\text{II}$ and $U^\text{IV}$ is the point $P^\text{II}_4 = P^\text{IV}_8$, which does not belong to $U^\text{IV}_{3 : 6}$.

Third, $U^\text{II}$ and $U^\text{III}_{4 : 5}$. Since $U^\text{II} = U - \bfomega_2$ and $U^\text{III} = U + \bfomega_2$, the $f_2$-goodness of $U$ implies that $U^\text{II}$ is already disjoint from $U^\text{III}_{3 : 5} \supseteq U^\text{III}_{4 : 5}$.

Fourth, $U^\text{I}$ and $U^\text{IV}_{3 : 6}$. We claim that, in fact, $U^\text{I}$ and $U^\text{IV}$ are disjoint altogether. Indeed, the displacement between them works out to $U^\text{I} \tto U^\text{IV} = 2\bfs_2 - \bfs_3 - \bfs_4$. However, since $U$ is $f_1$-good, this vector is a separator for $U$ by Lemmas \ref{sepg} and \ref{exf} for the $f_1$-lift, with $\bfrho = \bfrho_2 - \bfrho_4$ and $\bfdelta = \mathbf{0}$. \end{myproof}

\begin{lemma} \label{f1fo} Suppose that $U$ is $\fodd$-good. Then $f_1(U)$ is $\fodd$-good as well. \end{lemma}

\begin{myproof} Unlike the proofs of the other lemmas in this sub-section, this one does not require the use of separators.

Unpacking the definitions, we are given that $U$ is both $f_1$-good and $f_3$-good; $f_1(U)$ is $f_3$-good; and $f_3(U)$ is $f_1$-good. From these premises, we must derive that: (i) $f_1(U)$ is $f_1$-good; (ii) $f_1(U)$ is $f_3$-good; (iii) $f_1 \circ f_1(U)$ is $f_3$-good; and (iv) $f_3 \circ f_1(U)$ is $f_1$-good.

Of these, (i) follows by Lemma \ref{f1f1} and (ii) is given to us free of charge. Furthermore, since we have ticked off (i) and (ii) already, we can drop (iii) by Lemma \ref{f1f3}. This only leaves (iv).

Since $f_1$ and $f_3$ commute, we can rewrite (iv) to say that $f_1 \circ f_3(U)$ is $f_1$-good. However, in this form it follows immediately by Lemma \ref{f1f1} applied to $f_3(U)$. \end{myproof}

\begin{lemma} \label{f1fe} Suppose that $U$ is fully $f$-good. Then $f_1(U)$ is $\feven$-good. \end{lemma}

\begin{myproof} Define the vectors $\bfomega^{(2)}_1$, $\bfomega^{(2)}_2$, $\bfomega^{(2)}_3$, $\bfomega^{(2)}_4$ as in the proof of Lemma \ref{f1f2}. Then also set \begin{align*} V^{--} &= V - \bfomega^{(2)}_2 - \bfomega^{(2)}_4 & V^{-+} &= V - \bfomega^{(2)}_2 + \bfomega^{(2)}_4\\ V^{+-} &= V + \bfomega^{(2)}_2 - \bfomega^{(2)}_4 & V^{++} &= V + \bfomega^{(2)}_2 + \bfomega^{(2)}_4. \end{align*}

\begin{figure}[ht] \centering \includegraphics{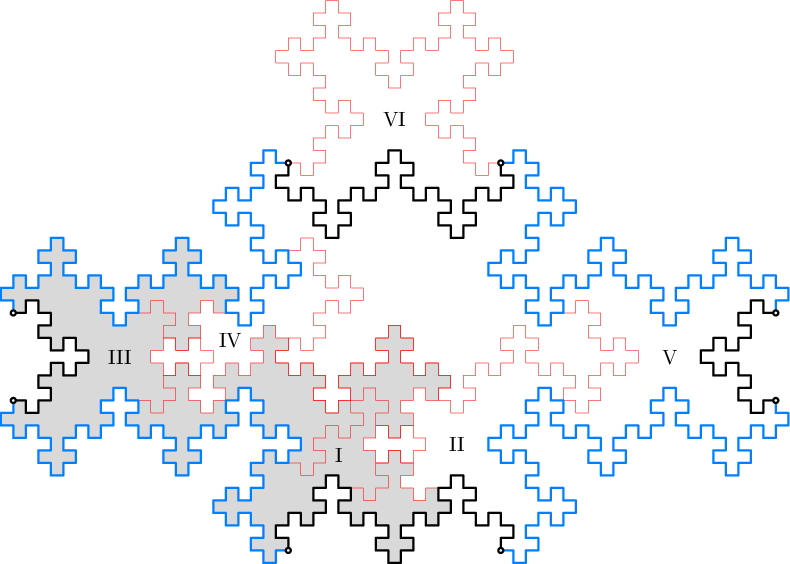} \caption{} \label{f1fef} \end{figure} 

By Lemma \ref{f1f2}, we get that $f_1(U)$ is both $f_2$-good and $f_4$-good. So, by Lemma \ref{isfog}, it suffices to show that $V^{--}$ is disjoint from $V^{++}_5$ and $V^{-+}$ is disjoint from $V^{+-}_3$. (See Figure \ref{f1fef}.)

Let \begin{align*} U^\text{I} &= U^- - \bfomega^{(2)}_2 - \bfomega^{(2)}_4 & U^\text{II} &= U^+ - \bfomega^{(2)}_2 - \bfomega^{(2)}_4\\ U^\text{III} &= U^- - \bfomega^{(2)}_2 + \bfomega^{(2)}_4 & U^\text{IV} &= U^+ - \bfomega^{(2)}_2 + \bfomega^{(2)}_4. \end{align*}

Since $V \subseteq U^- \cup U^+$, translating suitably we get that $V^{--} \subseteq U^\text{I} \cup U^\text{II}$ and $V^{-+} \subseteq U^\text{III} \cup U^\text{IV}$.

Furthermore, $V_3 \subseteq U^+$ and $V_5 \subseteq U + \bfeta_1$ by (\textbf{F}$_1$). Translating suitably, we get that \begin{align*} V^{+-}_3 &\subseteq U^\text{V} = U^+ + \bfomega^{(2)}_2 - \bfomega^{(2)}_4\\ V^{++}_5 &\subseteq U^\text{VI} = U + \bfeta_1 + \bfomega^{(2)}_2 + \bfomega^{(2)}_4. \end{align*}

With this, everything boils down to showing that both of $U^\text{I}$ and $U^\text{II}$ are disjoint from $U^\text{VI}$; and both of $U^\text{III}$ and $U^\text{IV}$ are disjoint from $U^\text{V}$.

\begin{table}[ht] \centering \begin{tabular}{|c|c|c|c|}
$\bfalpha$ & $\frf$ & $\bftau$ & $\bfdelta$\\
\hline
$(-1, 2, 1, -1)$ & $(f_2, f_4)$ & $\bftau_2$ & $\bfdelta_1 - \bfdelta_0$\\
\hline
$(-1, 1, 2, -1)$ & $(f_2, f_4)$ & $\bftau_3$ & $\bfdelta_2 - \bfdelta_0$\\
\hline
$(1, 2, -2, -1)$ & $(f_1, f_3)$ & $2\bftau_1$ & $\bfdelta_0 - \bfdelta_1$\\
\hline
$(1, 1, -1, -1)$ & $(f_1, f_3)$ & $2\bftau_1$ & $\bfdelta_0 - \bfdelta_3$\\
\hline
\end{tabular} \caption{} \label{f1fet} \end{table}

Set $\ols = (\bfs_1, \bfs_2, \bfs_3, \bfs_4)$. Then, using the shorthand notation of Section \ref{vec}, we can record every linear combination of the neighbourhood vectors of $U$ in the form $\bfalpha \cdot \ols$. For example, the displacement between $U^\text{I}$ and $U^\text{VI}$ works out to $U^\text{I} \tto U^\text{VI} = -\bfs_1 + 2\bfs_2 + \bfs_3 - \bfs_4$. So, in this case, $\bfalpha = (-1, 2, 1, -1)$.

The four displacements we need to analyse for the proof are $U^\text{I} \tto U^\text{VI}$, $U^\text{II} \tto U^\text{VI}$, $U^\text{III} \tto U^\text{V}$, and $U^\text{IV} \tto U^\text{V}$. The corresponding coefficient sequences are listed in the first column of Table~\ref{f1fet}, in this order from top to bottom. Each one of these displacements is indeed a separator for $U$ by Lemmas \ref{sepgg} and \ref{exfo}. The values of $\frf$, $\bftau$, and $\bfdelta$ which must be plugged in are summarised in the remaining columns of Table \ref{f1fet}. \end{myproof}

Together, Lemmas \ref{f1fo} and \ref{f1fe} take care of the induction step for the $f$-lifts.

\subsection{Induction Step for the \texorpdfstring{$g^\star$}{g*}-Lifts} \label{stepg}

We continue with the induction step for the $g^\star$-lifts. Just as in the previous sub-section, it is clear that the two of them behave analogously to one another. So we are only going to consider $\gso$ in detail, and the corresponding results for $\gse$ will be implied.

Suppose, now, that $U$ is fully $f$-good. We must show that $\gso(U)$ is fully $f$-good as well.

Let $V = \gso(U)$. By (\textbf{G}$_\text{\rm Odd}$), we get that \begin{align*} V_1 \subseteq U^\text{I} &= U + \half(\bfs_1 + \bfs_2 + \bfs_3 + \bfs_4)\\ V_2 \subseteq U^\text{II} &= U + \half(-\bfs_1 + \bfs_2 + \bfs_3 + \bfs_4)\\ &\ldots\\ V_8 \subseteq U^\text{VIII} &= U + \half(\bfs_1 + \bfs_2 + \bfs_3 - \bfs_4). \end{align*}

\begin{figure}[ht] \null \hfill \begin{subfigure}[c]{180pt} \centering \includegraphics{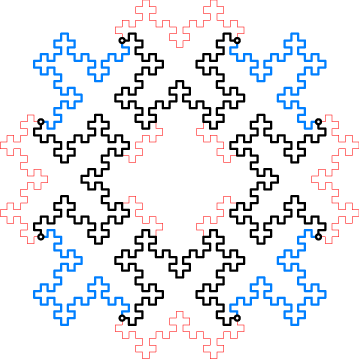} \caption{} \label{wreath-go} \end{subfigure} \hfill \begin{subfigure}[c]{180pt} \centering \includegraphics{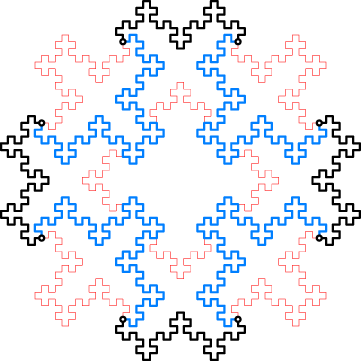} \caption{} \label{wreath-ge} \end{subfigure} \hfill \null \caption{} \label{wreaths} \end{figure}

Notice that this is true of $\gse$ as well, with the exact same $U^\text{I}$, $U^\text{II}$, $\ldots$, $U^\text{VIII}$. For example, Figures \ref{wreath-go} and \ref{wreath-ge} show the $\gso$ and $\gse$-lifts of the double chain $U$ in Figure \ref{chain-ge-go}, respectively, with the ``surplus'' parts of $U^\text{I}$, $U^\text{II}$, $\ldots$, $U^\text{VIII}$ outlined in red.

Just as before, our reasoning will ultimately bottom out at showing that certain vectors are separators for $U$. Fortunately, by this point we have become very efficient at doing this.

Indeed, from the previous sub-section we know that if $U$ is fully $f$-good, then so are all of its $f$-lifts as well. By induction on $\frf$, it follows that a fully $f$-good $U$ is also $\frf$-good for all finite sequences of $f$-lifts $\frf$. Thus, for a fully $f$-good $U$, we can safely apply Lemma \ref{sepgg} with $\frf$ arbitrary. We will make extensive use of this observation throughout the present sub-section.

In order to avoid doing the same work multiple times, we proceed now to factor out some symmetries. Let $\bfalpha = (\alpha_1, \alpha_2, \alpha_3, \alpha_4)$ be a coefficient sequence and set $\ols = (\bfs_1, \bfs_2, \bfs_3, \bfs_4)$, as in the proof of Lemma \ref{f1fe}. Suppose also that the full $f$-goodness of $U$ and Lemma \ref{sepgg} together imply that $\bfalpha \cdot \ols$ is a separator for $U$.

Clearly, if $U$ is fully $f$-good, then so is $U^\shift{1}$ as well. (Strictly speaking, $U^\shift{1}$ is not a descendant chain; but recall our remark at the beginning of this section about our definitions and lemmas being more generally applicable.) By Lemma \ref{sepgg} with $U^\shift{1}$ in place of $U$, we get that $\alpha_1\bfs_2 + \alpha_2\bfs_3 + \alpha_3\bfs_4 + \alpha_4\bfs_5$ must be a separator for $U$. However, the latter vector is really just $(-\alpha_4, \alpha_1, \alpha_2, \alpha_3) \cdot \ols$.

So we define a \emph{skew rotation} of $\bfalpha$ to be any coefficient sequence obtained from $\bfalpha$ by iterating the transformation $(\alpha_1, \alpha_2, \alpha_3, \alpha_4) \to (-\alpha_4, \alpha_1, \alpha_2, \alpha_3)$. There are a total of eight distinct skew rotations of $\bfalpha$ when $\bfalpha \neq \mathbf{0}$.

From the point of view of our toolbox, the skew rotations of $\bfalpha$ are mutually interchangeable; if we can prove (using the full $f$-goodness of $U$ and Lemma \ref{sepgg}) that one of them is the coefficient sequence of a separator for $U$, then we can prove this for all of them.

Hence, in order to process a large batch of would-be separators ``in bulk'', we may choose one representative skew rotation out of each coefficient sequence -- hopefully in such a way that the number of distinct representatives is much smaller than the number of would-be separators in the batch. This is exactly what we are going to do in the proofs of Lemmas \ref{gint}--\ref{gfpar}.

For convenience, in these proofs we are going to treat Lemmas \ref{exf} and \ref{exfo} as ``supplements'' to Lemma \ref{sepgg}, and we will assume that they have been invoked implicitly whenever needed. On the other hand, in all instances where $\frf \not \in \{f_1, f_2, f_3, f_4, (f_1, f_3), (f_2, f_4)\}$, just the original Lemma \ref{sepgg} will already be sufficiently strong for our purposes.

\begin{lemma} \label{gint} Suppose that $U$ is fully $f$-good. Then $\gso(U)$ is free of self-intersections. \end{lemma} 

\begin{myproof} Each part of $V$ is a translation copy, modulo reversal, of some slice of $U$. So none of them can contain self-intersections when considered in isolation.

Observe next that $U^\text{I} - \bfs_1 = U^\text{II}$. So $U^\text{I}$ and $U^\text{II}$ can be embedded into the tiling of $[U]$ which uses the lattice $\mathcal{L}(\bfs_1, \bfs_3)$. Since $V_1 \reveq U^\text{I}_{7 : 3}$ and $V_2 = U^\text{II}_{5 : 7}$, the combinatorial structure of this tiling tells us that $V_1$ and $V_2$ meet only at their common endpoint $Q_2 = P^\text{I}_7 = P^\text{II}_5$. Similar reasoning applies also to $V_i$ and $V_{i + 1}$ for all $i$.

We are left to verify that all other pairs of distinct parts of $V$ are disjoint. We will show that $U^\text{I}$ is disjoint from all of $U^\text{III}$, $U^\text{IV}$, $\ldots$, $U^\text{VII}$. The remaining cases can be handled analogously.

\begin{table}[ht] \centering \begin{tabular}{|c|c|c|c|}
$\bfalpha$ & $\frf$ & $\bftau$ & $\bfdelta$\\
\hline
$(0, 0, 1, 1)$ & $f_1$ & $\bftau_2 + \bftau_4$ & $\mathbf{0}$\\
\hline
$(0, 1, 1, 1)$ & $f_1$ & $2\bftau_2 + \bftau_4$ & $\mathbf{0}$\\
\hline
$(1, 1, 1, 1)$ & $(f_1, f_3)$ & $2\bftau_2$ & $\bfdelta_0 - \bfdelta_3$\\
\hline
\end{tabular} \caption{} \label{gintt} \end{table}

There are five pairwise displacements to consider. However, we can cover their coefficient sequences with only three representative skew rotations. These are listed in Table \ref{gintt}, together with the values of $\frf$, $\bftau$, and $\bfdelta$ which must be plugged into Lemma \ref{sepgg} in order to obtain them. \end{myproof}

For the proofs of Lemmas \ref{gf} and \ref{gfpar}, let \[\bfomega^\text{High}_i = \half(\bft_{i + 1} - \bft_{i + 2}) \qquad \text{for all $i = 1$, $2$, $3$, $4$}.\]

(These vectors are defined relative to $V$ in the same way as the vectors $\bfomega_1$, $\bfomega_2$, $\bfomega_3$, $\bfomega_4$ are defined relative to $U$.)

\begin{lemma} \label{gf} Suppose that $U$ is fully $f$-good. Then $\gso(U)$ is $f_i$-good for all $f$-lifts $f_i$. \end{lemma} 

\begin{myproof} We will show that $V$ is $f_1$-good. The other three $f$-lifts can be handled analogously.

\begin{figure}[ht] \centering \includegraphics{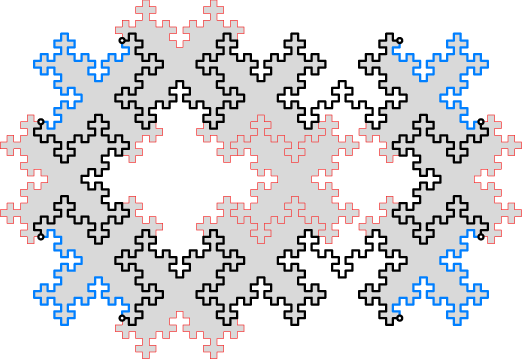} \caption{} \label{gof1} \end{figure}

By Lemma \ref{gint}, we get that $V$ is free of self-intersections. Hence, by Lemma \ref{isg}, what remains is to verify that the $\bfomega^\text{High}_1$ translation copies of $U^\text{I}$, $U^\text{II}$, $\ldots$, $U^\text{VIII}$ are disjoint from the $-\bfomega^\text{High}_1$ translation copies of $U^\text{VI}$, $U^\text{VII}$, $U^\text{VIII}$. (See Figure \ref{gof1}.)

\begin{table}[ht] \centering \begin{tabular}{|c|c|c|c|}
$\bfalpha$ & $\frf$ & $\bftau$ & $\bfdelta$\\
\hline
$(0, 0, 1, 2)$ & $f_1$ & $\bftau_2 + 2\bftau_4$ & $\bfdelta_1 - \bfdelta_0$\\
\hline
$(0, 0, 2, 1)$ & $f_1$ & $2\bftau_2 + \bftau_4$ & $\bfdelta_0 - \bfdelta_1$\\
\hline
$(0, 0, 2, 2)$ & $f_1$ & $2\bftau_2 + 2\bftau_4$ & $\mathbf{0}$\\
\hline
$(0, 1, 2, 1)$ & $f_1$ & $3\bftau_2 + \bftau_4$ & $\bfdelta_0 - \bfdelta_1$\\
\hline
$(0, 1, 2, 2)$ & $f_1$ & $3\bftau_2 + 2\bftau_4$ & $\mathbf{0}$\\
\hline
$(0, 2, 2, 1)$ & $f_2$ & $2\bftau_2 + 3\bftau_4$ & $\mathbf{0}$\\
\hline
$(1, 0, 1, 1)$ & $f_3$ & $\bftau_3$ & $\mathbf{0}$\\
\hline
$(1, 1, 0, 1)$ & $f_3$ & $\bftau_2$ & $\mathbf{0}$\\
\hline
$(1, 1, 2, 1)$ & $(f_2, f_3)$ & $\bftau_2 + 2\bftau_4$ & $\mathbf{0}$\\
\hline
$(1, 2, 1, 1)$ & $(f_4, f_3)$ & $2\bftau_1 + \bftau_3$ & $\mathbf{0}$\\
\hline
$(1, 1, 2, 2)$ & $(f_2, f_3)$ & $\bftau_2 + 3\bftau_4$ & $\mathbf{0}$\\
\hline
$(2, 2, 1, 1)$ & $(f_4, f_3)$ & $3\bftau_1 + \bftau_3$ & $\mathbf{0}$\\
\hline
$(1, 2, 2, 1)$ & $(f_2, f_3)$ & $2\bftau_2 + 2\bftau_4$ & $\bfdelta_0 - \bfdelta_3$\\
\hline
\end{tabular} \caption{} \label{gft} \end{table}

There are $24$ pairwise displacements to consider. However, we can cover their coefficient sequences with only $16$ representative skew rotations. Three of them are already in Table \ref{gintt}. The remaining $13$, together with their derivations, are listed in Table \ref{gft}. \end{myproof}

\begin{lemma} \label{gfpar} Suppose that $U$ is fully $f$-good. Then $\gso(U)$ is both $\fodd$-good and $\feven$-good. \end{lemma} 

\begin{myproof} We will show that $V$ is $\fodd$-good. Once again, $\feven$-goodness can be handled analogously.

\begin{figure}[ht] \centering \includegraphics{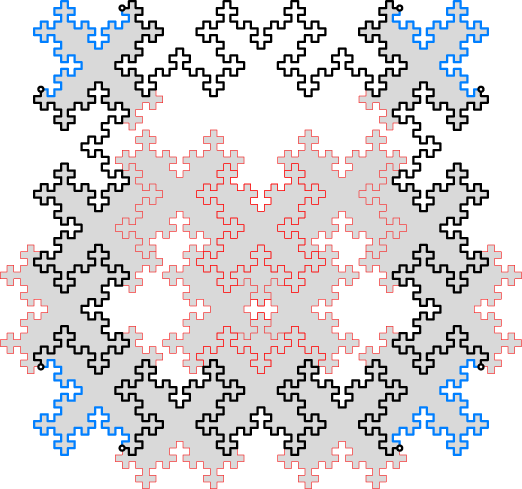} \caption{} \label{gofo} \end{figure}

By Lemma \ref{gf}, we get that $V$ is both $f_1$-good and $f_3$-good. Hence, by Lemma \ref{isfog}, what remains is to verify that the $\bfomega^\text{High}_1 + \bfomega^\text{High}_3$ translation copies of $U^\text{I}$, $U^\text{II}$, $\ldots$, $U^\text{VIII}$ are disjoint from the $-\bfomega^\text{High}_1 - \bfomega^\text{High}_3$ translation copy of $U^\text{VIII}$; and the $-\bfomega^\text{High}_1 + \bfomega^\text{High}_3$ translation copies of $U^\text{I}$, $U^\text{II}$, $\ldots$, $U^\text{VIII}$ are disjoint from the $\bfomega^\text{High}_1 - \bfomega^\text{High}_3$ translation copy of $U^\text{II}$. (See Figure \ref{gofo}.)

\begin{table}[ht] \centering \begin{tabular}{|c|c|c|c|}
$\bfalpha$ & $\frf$ & $\bftau$ & $\bfdelta$\\
\hline
$(1, 1, 1, 2)$ & $(f_1, f_3)$ & $2\bftau_2 + \bftau_4$ & $\bfdelta_0 - \bfdelta_1$\\
\hline
$(2, 1, 1, 1)$ & $(f_1, f_3)$ & $\bftau_1 + 2\bftau_3$ & $\bfdelta_0 - \bfdelta_1$\\
\hline
$(1, 2, 2, 2)$ & $(f_2, f_3)$ & $2\bftau_2 + 3\bftau_4$ & $\bfdelta_0 - \bfdelta_3$\\
\hline
$(2, 2, 2, 1)$ & $(f_4, f_3)$ & $3\bftau_1 + 2\bftau_3$ & $\bfdelta_0 - \bfdelta_3$\\
\hline
$(2, 2, 2, 2)$ & $(f_1, f_2, f_3)$ & $\bftau_2 + 3\bftau_4$ & $\bfdelta_5 - \bfdelta_0$\\
\hline
\end{tabular} \caption{} \label{gfpart} \end{table}

There are $16$ pairwise displacements to consider. However, we can cover their coefficient sequences with only eight representative skew rotations. Three of them, namely $(1, 1, 1, 1)$, $(1, 1, 2, 2)$, and $(2, 2, 1, 1)$, are already in Tables \ref{gintt} and \ref{gft}. The remaining five, together with their derivations, are listed in Table \ref{gfpart}. \end{myproof}

Lemma \ref{gfpar} completes the induction step for the $g^\star$-lifts, and with this also the proof of the positive part of Theorem \ref{mt}.

\section{The Generalisation to Polygons} \label{poly}

Here, we generalise our main result from polyominoes to polygons.

We state the generalisation in Section \ref{sg}. The positive part of it we get at no extra cost. For the negative part, first we introduce the notion of alphabetisability in Section \ref{alph}. The gist is that an alphabetisable double chain of polylines maps onto a double chain of words -- and we already know how to analyse such chains. Then, in Section \ref{adc}, we complete the proof by showing that all double chains of polylines which correspond to $4$-neighbour double tiles are in fact alphabetisable.

\subsection{Statement of the Generalisation} \label{sg}

The generalisation of our main result from polyominoes to polygons goes as follows:

\begin{theorem} \label{pt} Every polygon which admits two different $4$-neighbour lattice tilings of the plane is a deformation of some polyomino described by an $\{f_1, f_2, f_3, f_4, \gso, \gse\}$-descendant of the Greek cross. Conversely, each such descendant of the Greek cross describes a polyomino, and every deformation of such a polyomino admits two different $4$-neighbour lattice tilings of the plane. \end{theorem} 

This differs from the statement of Theorem \ref{mt} by a single word -- the first instance of ``polyomino'' has been changed to ``polygon''.

Once again, the result splits into a positive part and a negative part. The positive part we get for free; indeed, the only difference with the positive part of Theorem \ref{mt} is that, this time around, the word ``deformation'' is to be understood in the sense of ``general polygonal deformation''. Of course, the material of Section \ref{def} applies equally well to deformations in the combinatorial setting of polyominoes and deformations in the geometric setting of general polygons.

On the other hand, the negative part of the generalisation differs from the negative part of Theorem \ref{mt} in one crucial way: The natural association between polyominoes and words over finite alphabets does not carry over to general polygons in an obvious manner. Hence, the syntactic machinery we developed in Section \ref{neg} is not obviously applicable. We bridge this gap over the course of the rest of the present section.

Notice that the convention of Section \ref{tile} regarding polyomino tilings considered in a combinatorial context does not apply to the generalisation. Here, we are considering all polygons (polyominoes included) in a purely geometric context instead.

\subsection{Alphabetisability} \label{alph}

Our proof of the negative part of Theorem \ref{mt} in Section \ref{neg} was carried out entirely in terms of manipulations with words. The reason we could afford to do this is because polyominoes are naturally encodable by means of words over $\mathcal{A}_\square$.

For the negative part of Theorem \ref{pt}, our strategy will be to extend this kind of encodability to general polygons. Or, rather, to general polygonal $4$-neighbour double tiles.

We say that an interleaved double chain of polylines $F$ is \emph{alphabetisable} when we can find an alphabet with reversal $\mathcal{A}$, a valuation $\psi$, and an interleaved double chain $U$ of words over $\mathcal{A}$ such that $F = \psi(U)$. The definition of a non-interleaved double chain of polylines $F$ being alphabetisable is analogous.

Let $\Phi$ be a $4$-neighbour double tile. Suppose that $F$ is the double chain of polylines associated with $\Phi$ as in Section \ref{dcp}, and that it admits the alphabetisation $U$.

Recall that, in Section \ref{neg}, we did not make any assumptions about the underlying alphabet. (Even though our end goal was a result about polyominoes, with $\mathcal{A} = \mathcal{A}_\square$.) Thus all of Lemmas \ref{nil}--\ref{ril} apply to $U$.

So does the analysis at the end of Section \ref{root}, with one caveat: When $\mathcal{A} \neq \mathcal{A}_\square$, the fact that $R$ is free of combinatorial self-intersections does not in general imply that $\psi(R)$ is free of geometric self-intersections. To patch things up, we extend the definition of a stable self-intersection as follows: The phrase ``a nonempty balanced subword'' is replaced with ``a nonempty subword whose $\psi$-image is either closed or it contains a geometric self-intersection''. Clearly, with the revised definition we can deduce that $\psi(R)$ must indeed be free of geometric self-intersections. Otherwise, reasoning as in the proof of Lemma \ref{ril}, we get a contradiction with $F = \psi(U)$ being free of geometric self-intersections.

Or, in summary, all of the material in Section \ref{neg} will go through for $U$. Therefore, the negative part of Theorem \ref{pt} will be true of $\Phi$, as desired.

With this, everything boils down to showing that if a double chain of polylines is the boundary of a polygon, then it must be alphabetisable. The proof is in Section \ref{adc}.

We close this sub-section by demonstrating that the notion of alphabetisability is indeed meaningful; i.e., that not all double chains of polylines are alphabetisable. For the sake of clarity, first we will exhibit a degenerate counterexample where some parts of the chain are empty. Then we will show how it can be augmented.

Let $\kappa(a)$ be the floating polyline represented by $(0, 0) \sto (a, 0)$. Given two positive real numbers $a$ and $b$, consider the interleaved double chain of floating polylines \[F(a, b) = \kappa(a) : \kappa(b) : \varepsilon : \varepsilon : \kappa(a)^{-1} : \kappa(b)^{-1} : \varepsilon : \varepsilon.\]

Suppose that $F(a, b)$ is alphabetisable by means of the word chain $U$. If $a > b$, then $U$ admits a true $f_1$-reduction and $f_1^{-1}(U)$ alphabetises $F(a - b, b)$. Similarly, if $a < b$, then $U$ admits a true $f_2$-reduction and $f_2^{-1}(U)$ alphabetises $F(a, b - a)$. By iterating these two reductions, essentially we are carrying out the Euclidean algorithm over $a$ and $b$.

Suppose now that $a/b$ is irrational. Then $U$ will admit an infinite series of true $f$-reductions. This is a contradiction. Hence, if $a/b$ is irrational, then $F(a, b)$ cannot be alphabetised.

For a counterexample where all parts of the chain are nonempty, simply take a suitable descendant of $F(a, b)$. One such would be \[f_4 \circ f_3(F(a, b)) = \kappa(a) : \kappa(b) : \kappa(b) : \kappa(a)^{-1}\kappa(b) : \kappa(a)^{-1} : \kappa(b)^{-1} : \kappa(b)^{-1} : \kappa(a)\kappa(b)^{-1}.\]

\subsection{The Alphabetisability of Double Chains} \label{adc}

We go on now to an investigation of the alphabetisability of double chains.

Let $F$ be a double chain of literal polylines. We define the \emph{vertices} and \emph{segments} of $F$ to be the vertices and segments of the parts of $F$. So, in particular, each division point of $F$ is also a vertex of $F$. We consider each empty part of $F$ to consist of a single empty segment; but we assume that the parts of $F$ are not all empty.

Let $\Gamma$ be a circle of the same length as $F$. Fix also a mapping $\xi : \Gamma \to \mathbb{R}^2$ such that, as $P$ makes one full counterclockwise revolution around $\Gamma$, its image $Q = \xi(P)$ traces out $F$; and at all times the length of the circular arc traversed by $P$ equals the length of the polyline traced out by~$Q$. So, in geometric terms, $\xi$ is essentially an arclength parametrisation of $F$. Notice that $\xi$ might not be a bijection, because $F$ is allowed to contain self-intersections.

We can now subdivide $\Gamma$ into closed arcs which, via $\xi$, biject onto the segments of $F$. (Of course, the arcs will overlap at their endpoints.) Denote these arcs by $\gamma_1$, $\gamma_2$, $\ldots$, $\gamma_k$ and the segments of $F$ by $s_1$, $s_2$, $\ldots$, $s_k$, so that $\gamma_i$ bijects onto $s_i$ for all $i$.

Since $F$ is a double chain, as in Section \ref{dcp} we can regroup it into $A_1 : B_1 : A_2 : B_2$ with $A_1 \revsim A_2$ and $B_1 \revsim B_2$; and we can also regroup some cyclic shift of $F$ into $C_1 : D_1 : C_2 : D_2$ with $C_1 \revsim C_2$ and $D_1 \revsim D_2$. Let $A_1$, $B_1$, $A_2$, $B_2$ map onto the closed arcs $A^\circ_1$, $B^\circ_1$, $A^\circ_2$, $B^\circ_2$ of $\Gamma$, respectively, and define $C^\circ_1$, $D^\circ_1$, $C^\circ_2$, $D^\circ_2$ similarly.

Since $|A^\circ_1| = |A^\circ_2|$ and $|B^\circ_1| = |B^\circ_2|$, there exist two perpendicular lines $\ell_A$ and $\ell_B$ through the center of $\Gamma$ such that reflection across $\ell_A$ swaps $A^\circ_1$ and $A^\circ_2$; while reflection across $\ell_B$ swaps $B^\circ_1$ and $B^\circ_2$. Define $\ell_C$ and $\ell_D$ similarly.

Construct the graph $G_\text{I}$ on vertex set $\Gamma$ as follows: Suppose, for concreteness, that $|A^\circ_1| = |A^\circ_2| \neq 0$. Whenever two points of $\lopen A^\circ_1 \ropen$ and $\lopen A^\circ_2 \ropen$ are reflections of one another across $\ell_A$, we join them by an edge; and similarly whenever two points of $B^\circ_1$ and $B^\circ_2$ are reflections of one another across $\ell_B$. Then $G_\text{I}$ will be a perfect matching over $\Gamma$. The construction of $G_\text{II}$ is analogous, but based on $C^\circ_1$, $D^\circ_1$, $C^\circ_2$, $D^\circ_2$ instead.

Given two nonempty oriented segments $s'$ and $s''$, we write $s' \eqdir s''$ when $s'$ and $s''$ are parallel and point in the same direction. We also write $s' \opdir s''$ when $s'$ and $s''$ are parallel but point in opposite directions. We use the same notation for nonzero vectors, too.

Suppose that an edge of $G_\text{I}$ joins the points $P \in \lopen \gamma_i \ropen$ and $Q \in \lopen \gamma_j \ropen$. Then $s_i \opdir s_j$. This follows because a small vicinity of $\xi(P)$ contained within $\lopen s_i \ropen$ and a small vicinity of $\xi(Q)$ contained within $\lopen s_j \ropen$ must be related by the same translation as one of the pairs $A_1 \revsim A_2$ and $B_1 \revsim B_2$. (Though we cannot claim that $s_i \revsim s_j$ in general; e.g., consider what is going on in the counterexamples of the previous sub-section.) The same observation applies to $G_\text{II}$ as well.

Let $G = G_\text{I} \cup G_\text{II}$; i.e., $G$ is the graph on vertex set $\Gamma$ whose edge set is the union of the edge sets of $G_\text{I}$ and $G_\text{II}$. Since both of $G_\text{I}$ and $G_\text{II}$ are perfect matchings, each connected component of $G$ is either a cycle or a doubly infinite path. (Notice that we allow multiple edges and two-edge ``cycles'' in $G$.) Furthermore, each such cycle or doubly infinite path must alternate between edges of $G_\text{I}$ and edges of $G_\text{II}$.

Call a connected component of $G$ \emph{special} when one of its vertices $\xi$-maps onto a vertex of~$F$. Otherwise, call it \emph{regular}. We also say that a vertex or edge of $G$ is \emph{special} when it belongs to a special connected component; and similarly for it being \emph{regular}. Observe that there are only finitely many special connected components. Since each connected component of $G$ is countable, we conclude that the regular ones must be the majority.

We proceed now to look into the structure of an arbitrary connected component of $G$. Let $P_0$ be any vertex of $G$ and let $Q_0$ be the point diametrically opposite $P_0$ on $\Gamma$. Denote also $R_0 = \{P_0, Q_0\}$. The connected components of $P_0$ and $Q_0$ will turn out to be closely interlinked, and we are going to analyse them together.

Let $P_1$ be the $G_\text{I}$-neighbour of $P_0$ and let $Q_1$ be the $G_\text{I}$-neighbour of $Q_0$. Either both of these neighbours are obtained by reflection across $\ell_A$, or else both of them are obtained by reflection across $\ell_B$. However, the reflections of $R_0$ across $\ell_A$ and $\ell_B$ coincide because $\ell_A \perp \ell_B$. Denote the resulting pair of points by $R_1$. Then we can tell that $R_1 = \{P_1, Q_1\}$, even though we do not know which reflection exactly has produced $P_1$ and $Q_1$.

Similarly, let $P_{-1}$ be the $G_\text{II}$-neighbour of $P_0$ and let $Q_{-1}$ be the $G_\text{II}$-neighbour of $Q_0$. Once again, the reflections of $R_0$ across $\ell_C$ and $\ell_D$ coincide, and we denote the resulting pair of points by $R_{-1}$, with $R_{-1} = \{P_{-1}, Q_{-1}\}$.

Continuing on in the same way, for all even positive integers $i$ let $P_{i + 1}$ and $Q_{i + 1}$ be the $G_\text{I}$-neighbours of $P_i$ and $Q_i$; and also let $R_{i + 1} = \{P_{i + 1}, Q_{i + 1}\}$ be the reflection of $R_i$ across both of $\ell_A$ and $\ell_B$. Similarly, for all odd positive integers $i$ let $P_{i + 1}$ and $Q_{i + 1}$ be the $G_\text{II}$-neighbours of $P_i$ and $Q_i$; and also let $R_{i + 1} = \{P_{i + 1}, Q_{i + 1}\}$ be the reflection of $R_i$ across both of $\ell_C$ and $\ell_D$. For all negative integers $i$, we define $P_{i - 1}$, $Q_{i - 1}$, $R_{i - 1}$ in terms of $P_i$, $Q_i$, $R_i$ analogously, but with the words ``even'' and ``odd'' swapped. Then the connected components of $P$ and $Q$ will be $\cdots$---$P_{-1}$---$P_0$---$P_1$---$\cdots$ and $\cdots$---$Q_{-1}$---$Q_0$---$Q_1$---$\cdots$, respectively.

Let $\theta$ be the counterclockwise angle from $\ell_A$ to $\ell_C$. Recall, from elementary Euclidean geometry, that two consecutive reflections across two intersecting lines amount to a rotation by twice the angle between these lines. Hence, for all integers $i$, a rotation by $(-1)^i \cdot 2\theta$ around the center of $\Gamma$ maps $R_i$ onto $R_{i + 2}$. This observation will be crucial to our argument. Let $R_\text{Even}$ be the union of all $R_i$ with $i$ even, and define $R_\text{Odd}$ analogously.

Over the course of the rest of the proof, we consider two independent cases based on the rationality of $\theta$. (Of course, by $\theta$ being rational we mean that it is a rational multiple of $360^\circ$.) The two cases are covered by Lemmas \ref{it} and \ref{rt}, respectively.

We define a polyline or a chain of polylines to be \emph{flat} when all of its nonempty segments are pairwise parallel.

\begin{lemma} \label{it} Suppose that $\theta$ is irrational. Then $F$ is flat. \end{lemma} 

\begin{myproof} Choose $P_0$ and $Q_0$ so that both of them are regular. This is indeed possible because the union of all special connected components is countable.

For all $i$, let $P_i \in \lopen \gamma_{p(i)} \ropen$. Then $s_{p(i)} \opdir s_{p(i + 1)}$ for all $i$. It follows that there exists some nonzero vector $\bfp$ such that $s_{p(i)} \eqdir \bfp$ for all even $i$ and $s_{p(i)} \opdir \bfp$ for all odd $i$. Define the mapping $q : \mathbb{Z} \to \{1, 2, \ldots, k\}$ and the nonzero vector $\bfq$ similarly.

Since $R_i$ and $R_{i + 2}$ are related by a rotation of $2\theta$ for all even integers $i$, and $\theta$ is irrational, we get that $R_\text{Even}$ is dense in $\Gamma$. So each index $j$ with $s_j$ nonempty is of the form either $p(i)$ or $q(i)$ for some even integer $i$. Hence, all nonempty segments $s$ of $F$ satisfy either $s \eqdir \bfp$ or $s \eqdir \bfq$. By the exact same reasoning, but based on $R_\text{Odd}$ instead, also all nonempty segments $s$ of $F$ satisfy either $s \opdir \bfp$ or $s \opdir \bfq$.

However, these two conclusions can coexist only when $\bfp \parallel \bfq$. Thus the nonempty segments of $F$ must be pairwise parallel, as desired. \end{myproof}

Clearly, a flat chain of polylines can never be the boundary of a polygon. This means that we can exclude the irrational $\theta$ from consideration.

\begin{lemma} \label{rt} Suppose that $\theta$ is rational. Then $F$ is alphabetisable. \end{lemma} 

\begin{myproof} Since $R_i$ and $R_{i + 2}$ are related by a rotation of $2\theta$ for all even integers $i$, and $\theta$ is rational, we get that $R_\text{Even}$ is finite for all initial choices of $P_0$ and $Q_0$. Similarly, so is $R_\text{Odd}$ as well. However, the union $R_\text{Even} \cup R_\text{Odd}$ covers the connected components of $P_0$ and $Q_0$. We conclude that, if $\theta$ is rational, then every connected component of $G$ is finite.

So, in particular, every special connected component is finite. It follows that the vertices of the special connected components subdivide $\Gamma$ into a finite number of closed arcs $\delta_1$, $\delta_2$, $\ldots$,~$\delta_\ell$. This subdivision is a refinement of the subdivision of $\Gamma$ into $\gamma_1$, $\gamma_2$, $\ldots$, $\gamma_k$, but with the empty $\gamma$-arcs omitted. The subdivision of each nonempty $\gamma$-arc into $\delta$-arcs induces, via $\xi$, also a subdivision of the corresponding nonempty segment of $F$ into sub-segments. For all $i$, let $\xi$ map $\delta_i$ onto the sub-segment $t_i$.

Observe that the regular edges of $G_\text{I}$ can be partitioned into ``sheaves'' so that the edges of each sheaf form a perfect matching between $\lopen \delta_i \ropen$ and $\lopen \delta_j \ropen$ for some pair of $\delta$-arcs $\delta_i$ and $\delta_j$; and similarly for $G_\text{II}$ as well.

Construct a new graph $H$ on the vertex set $\{t_1, t_2, \ldots, t_\ell\}$ by collapsing each such sheaf into an edge of $H$ joining the sub-segments $t_i$ and $t_j$. Then the edge $t_i$---$t_j$ of $H$ guarantees that $t_i \revsim t_j$. Furthermore, each one of the connected components $H_1$, $H_2$, $\ldots$, $H_d$ of $H$ will be an even cycle which alternates between edges derived from the sheaves of $G_\text{I}$ and edges derived from the sheaves of $G_\text{II}$.

Consider the alphabet with reversal $\mathcal{A} = \{\mathsf{A}_1, \mathsf{A}_1^{-1}, \mathsf{A}_2, \mathsf{A}_2^{-1}, \ldots, \mathsf{A}_d, \mathsf{A}_d^{-1}\}$. For all $1 \le i \le d$, we label the sub-segments of $H_i$ with the letters $\mathsf{A}_i$ and $\mathsf{A}_i^{-1}$ of $\mathcal{A}$, in an alternating fashion. Let also $\psi$ be the valuation which maps the letters of $\mathcal{A}$ onto the floating polylines represented by the sub-segments of $F$ labelled with them.

We can now associate each part of $F$ with the word obtained by reading off the labels of its constituent sub-segments, in order. Let these words form the word chain $U$ over $\mathcal{A}$. Clearly, $\mathcal{A}$, $\psi$, and $U$ together give us the desired alphabetisation. \end{myproof}

Thus we can deal away with the irrational $\theta$ by Lemma \ref{it}; and Lemma \ref{rt} tells us that the rational $\theta$ always yield alphabetisable double chains of polylines. This completes the proof of Theorem \ref{pt}.

\section{Gaussian Primes} \label{prime}

Here, we discuss the question of \cite{M} about the special case of the double tile problem when the lattices of the two tilings are induced by a given Gaussian prime.

Section \ref{qs} reviews this question in its original context. Section \ref{pro} is a brief technical detour which will allow us to express our subsequent observations more neatly. Finally, in Section \ref{clover} we apply our results to the $4$-neighbour restriction of the Gaussian prime double tile problem.

\subsection{Statement of the Question} \label{qs}

Let $z = a + bi$ be a nonzero Gaussian integer. The complex multiples of $z$ are the integer linear combinations of $z$ and $iz = -b + ai$. So we can identify them with the vectors of $\mathcal{L}((a, b), (-b, a))$. We denote this lattice by $\mathcal{L}(z)$. Similarly, $\mathcal{L}(\overline{z}) = \mathcal{L}((b, a), (-a, b))$.

In \cite{M}, Mebane asks:

\begin{quotation} Let $a$, $b$ be positive integers such that $p = a^2 + b^2$ is an odd prime. Call a polyomino \emph{good} if it does not contain two squares with coordinates $(x_1, y_1)$, $(x_2, y_2)$ such that $(x_1 - x_2) + (y_1 - y_2)i$ is a Gaussian integer multiple of either $a + bi$ or $a - bi$. Is it always true that the number of good polyominoes of size $p$ is $2$, up to congruence? \end{quotation}

Here, we call such polyominoes $z$-\emph{good} or $(a, b)$-\emph{good}. For convenience, we define $z$ to be a \emph{proper} Gaussian prime when both of $a$ and $b$ are positive integers and $p = |z|^2$ is an odd prime. Suppose, from now on, that $z$ does satisfy the conditions of this definition.

Observe that a polyomino $\Phi$ is $z$-good if and only if it is a double tile which uses the lattices $\mathcal{L}(z)$ and $\mathcal{L}(\overline{z})$. Indeed, consider the implicit infinite grid of unit square cells that $\Phi$ ``inhabits'', and partition the cells of this grid into equivalence classes based on $\mathcal{L}(z)$. The divisibility condition with $z$ can be rephrased as $\Phi$ meeting each equivalence class at most once; following which, the area condition becomes equivalent to $\Phi$ meeting each equivalence class exactly once. The argument for $\mathcal{L}(\overline{z})$ is analogous.

Clearly, the rectangle of size $1 \times p$ is always $z$-good. The conjectured second $z$-good polyomino comes from a construction due to Mebane (op.\ cit.), which we outline below. Notice that, in stark contrast with our methods, it is based on an analysis of the interior of the polyomino instead of its boundary.

Consider a square of size $(a + b) \times (a + b)$ subdivided into a grid $S$ of unit square cells. Let $G$ be the graph whose vertices are the cells of $S$, and where two cells are joined by an edge if and only if the displacement $(x, y)$ between them satisfies $\{|x|, |y|\} = \{a, b\}$.

\begin{figure}[ht] \null \hfill \begin{subfigure}[c]{90pt} \centering \includegraphics{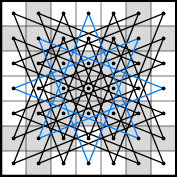} \caption{} \label{prime-a} \end{subfigure} \hfill \begin{subfigure}[c]{60pt} \centering \includegraphics{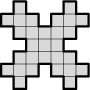} \caption{} \label{prime-b} \end{subfigure} \hfill \begin{subfigure}[c]{110pt} \centering \includegraphics{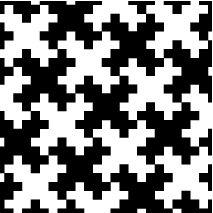} \caption{} \label{prime-c} \end{subfigure} \hfill \null \caption{} \label{pf} \end{figure} 

It is straightforward to check that $G$ consists of $(a - b)^2$ isolated vertices and some pairwise disjoint even cycles. (These cycles turn out to exhibit a rich ``fractal-like'' structure; see \cite{B}.) For example, Figure \ref{prime-a} shows the case of $z = 5 + 2i$ and $p = 29$.

Suppose that some $z$-good polyomino $\Phi$ fits inside of $S$. Then, for each edge of $G$, at most one of its endpoints can belong to $\Phi$. So $\Phi$ must contain all isolated vertices of $G$ together with one alternating half out of each cycle of $G$.

Conversely, suppose that we can choose one alternating half out of each cycle of $G$ so that these halves form a connected shape together with the isolated vertices of $G$. Then this shape will be a $z$-good polyomino. For example, the graph $G$ in Figure \ref{prime-a} yields the $(2, 5)$-good polyomino in Figure \ref{prime-b}. Its tiling which uses the lattice $\mathcal{L}(5 + 2i)$ is shown in Figure~\ref{prime-c}.

These observations do not guarantee that the construction will always work. However, experimentally, for small values of $p$ (say, for all $p < 100$) it yields exactly one $z$-good polyomino every time.

Hence, the implicit full form of the original question would be: Is it true that, for every proper Gaussian prime $z$, the complete list of $z$-good polyominoes consists of the rectangle of size $1 \times p$ together with a unique other polyomino obtained by means of the $(a + b) \times (a + b)$ construction?

\begin{figure}[ht] \null \hfill \begin{subfigure}[c]{70pt} \centering \includegraphics{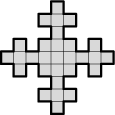} \caption{} \label{alt-a} \end{subfigure} \hfill \begin{subfigure}[c]{110pt} \centering \includegraphics{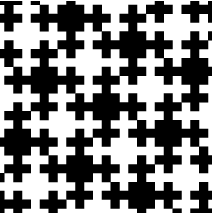} \caption{} \label{alt-b} \end{subfigure} \hfill \null \caption{} \label{alt} \end{figure}

We will show in Section \ref{clover} that some proper Gaussian primes $z$ admit additional $z$-good polyominoes. Our argument yields four counterexamples with $p < 100$ -- at $p = 29$, $73$, $89$, $97$. The first one of them is shown in Figure \ref{alt}.

Thus the original question ought to be understood somewhat more broadly: Given a proper Gaussian prime $z$, what can we say about the double tiles which use the lattices $\mathcal{L}(z)$ and $\mathcal{L}(\overline{z})$? This more general problem appears to be difficult. We will focus exclusively on its $4$-neighbour special case, and even in this restricted setting we are only going to offer some partial results.

One might wonder what happens when the condition that $a^2 + b^2$ must be an odd prime is omitted. The $(a, b)$-good polyominoes seem to be significantly less constrained then. For example (op.\ cit.), if $\gcd(a, b) = d \ge 2$, we can obtain a lot of them by subdividing a rectangle of height $d$ and area $a^2 + b^2$ into a grid of unit square cells and subsequently ``sliding'' each row of this grid slightly to the left or right. Furthermore, in Section \ref{clover} we are going to rely on the primality of $a^2 + b^2$ in an essential way in order to narrow down the $4$-neighbour $(a, b)$-good polyominoes.

On the other hand, Mebane observes (op.\ cit.) that, experimentally, for small positive integers $a$ and $b$ with $a + b$ odd, the $(a + b) \times (a + b)$ construction continues to yield a unique $(a, b)$-good polyomino even when $a^2 + b^2$ is composite.

\subsection{Proto-Descents} \label{pro}

Let $U$ be a descendant chain with neighbourhood vectors $\ols = (\bfs_1, \bfs_2, \ldots, \bfs_8)$. Let also $\varphi \in \psimain$ and $V = \varphi(U)$, with the neighbourhood vectors of $V$ being $\olt = (\bft_1, \bft_2, \ldots, \bft_8)$.

Observe that $\olt$ depends only on $\varphi$ and $\ols$; we do not need to know $U$ in order to calculate it. So we write $\varphi(\ols) = \olt$ as shorthand for ``the ordered octuples of vectors $\ols$ and $\olt$ satisfy the neighbourhood-octuple recurrence relations associated with $\varphi$ as in Section \ref{nv}''. We also let $\fodd(\ols) = f_1 \circ f_3(\ols)$ and $\feven(\ols) = f_2 \circ f_4(\ols)$.

Recall that $\gso(\ols) = \gse(\ols)$ for all $\ols$. So, in a context where we are interested more in $\ols$ than in $U$ itself, we can ``factor out'' the indistinguishability of $\gso$ and $\gse$. We write $\gsa(\ols) = \olt$ as a synonym for both of $\gso(\ols) = \olt$ and $\gse(\ols) = \olt$.

We define a \emph{proto-descent} for $U$ to be any finite sequence obtained by taking a descent for $U$ and then replacing each instance of $\gso$ or $\gse$ in it with an instance of $\gsa$. Hence, a proto-descent with $k$ instances of $\gsa$ will correspond to a total of $2^k$ descendant chains, by Proposition \ref{ud}. Furthermore, all of these descendant chains will share the same ordered octuple of neighbourhood vectors.

\begin{proposition} \label{un} The proto-descent of a descendant chain is determined uniquely, modulo $f$-commutativity, by the ordered octuple of its neighbourhood vectors. \end{proposition} 

\begin{myproof} Clearly, the neighbourhood vectors determine the partial vectors uniquely, and vice versa. So it suffices to prove the claim for the ordered octuple $\olu = (\bfu_1, \bfu_2, \ldots, \bfu_8)$ of the partial vectors of $U$.

Let $\bfphi$ be a proto-descent for $U$. In the notation of Section \ref{orient}, if $\bfphi$ ends with an instance of $f_1$, then $u_{13} > u_{23}$, $u_{34}$ and $u_{12}$, $u_{14} > u_{24}$. Conversely, if $\bfphi$ cannot be rearranged, using $f$-commutativity, so that it ends with an instance of $f_1$, then either $u_{12} \le u_{24}$ or $u_{14} \le u_{24}$. Similar observations apply also to the other three $f$-lifts. Furthermore, if $\bfphi$ ends with an instance of $\gsa$, then $u_{13}$, $u_{24} > u_{12}$, $u_{23}$, $u_{34}$, $u_{14}$.

From this point on, the proof proceeds along the same lines as the proof of Proposition \ref{ud} in Section \ref{unique}. \end{myproof}

\subsection{Clovers} \label{clover}

Let $\Phi$ be a $4$-neighbour $z$-good polyomino. By Theorem \ref{mt}, we get that $\Phi$ is a deformation of some polyomino determined by a descendant chain. Suppose that it is the $T$-deformation of the descendant polyomino $\Phi_\text{\rm Desc}$. Then $p = \area(\Phi) = \area(T) \cdot \area(\Phi_\text{\rm Desc})$. Since $p$ is prime and $\area(\Phi_\text{\rm Desc}) \ge 5$, it follows that $\area(T) = 1$. Or, in other words, $T$ is the unit square and the deformation induced by $T$ is the identity.

We conclude that every $4$-neighbour $z$-good polyomino is a descendant polyomino as well.

We define a descendant chain of literal polylines $U$ to be a \emph{clover} when its symmetry group coincides with the symmetry group of the axes-aligned unit square concentric with it. There is a simple necessary and sufficient condition for this in terms of the descent of $U$, as follows: Let \[\psisymm = \{\fodd, \feven, \gso, \gse\}.\]

Then $U$ is a clover if and only if it is a $\psisymm$-descendant of the Greek cross. The proof is straightforward.

Experimentally, for small values of $p$ (say, for all $p < 100$) the $(a + b) \times (a + b)$ construction of \cite{M} always yields a $4$-neighbour $z$-good polyomino determined by a clover. We proceed now to explore this connection.

\begin{proposition} \label{zgc} Let $z$ be a proper Gaussian prime. Then there exists a $z$-good polyomino determined by a clover. \end{proposition} 

\begin{myproof} Given a clover $U$, as in the previous sub-section let $k$ be the number of instances of $\gsa$ in its proto-descent and let $\ols$ be the ordered octuple of its neighbourhood vectors. Let also $\mystar(x, y)$ denote the ordered octuple of vectors \[((y, x), (x, y), (-x, y), (-y, x), (-y, -x), (-x, -y), (x, -y), (y, -x)).\]

We will show that, more generally, the ordered octuples $(-1)^k \cdot \ols$ obtained when $U$ ranges over all clovers coincide with the ordered octuples $\mystar(a, b)$ obtained when $a$ and $b$ range over all ordered pairs of positive integers with $a < b$, $a + b$ odd, and $\gcd(a, b) = 1$. Clearly, this would suffice.

The fact that $(-1)^k \cdot \ols$ is always of the desired form follows by straightforward induction on the $\psisymm$-descent of $U$ from the Greek cross. For the converse, we induct on $a + b$. The base case, when $a = 1$ and $b = 2$, corresponds to the Greek cross. For the induction step, observe that $\mystar(a, b) = \pm \varphi(\mystar(c, d))$ with $\varphi$, $c$, $d$ as in the first three columns of Table~\ref{cnt}; the plus and minus signs apply when $\varphi \in \{\fodd, \feven\}$ and $\varphi = \gsa$, respectively.

\begin{table}[ht] \centering \begin{tabular}{|c|c|c|c|}
$\varphi$ & $c$ & $d$ & constraint\\
\hline
$\fodd$ & $a$ & $b - 2a$ & $3a < b$\\
\hline
$\feven$ & $2a - b$ & $a$ & $b < 2a$\\
\hline
$\gsa$ & $b - 2a$ & $a$ & $2a < b < 3a$\\
\hline
\end{tabular} \caption{} \label{cnt} \end{table}

However, if $a + b \ge 5$, then exactly one of the three options yields $c$ and $d$ with $c < d$, $c + d$ odd, and $\gcd(c, d) = 1$. Explicitly, the constraints which $a$ and $b$ must satisfy in order for each option to be the winning one are listed in the fourth column of Table \ref{cnt}. \end{myproof}

Proposition \ref{zgc} shows that a non-rectangular $z$-good polyomino exists for all proper Gaussian primes $z$. One corollary of the proof is that, modulo $f$-commutativity, all clovers which determine $z$-good polyominoes share the same proto-descent. If this proto-descent contains $k$ instances of $\gsa$, then out of it we get a total of $2^k$ clovers which will yield the same number of non-rectangular $z$-good polyominoes.

The first few primes $p$ with $k \ge 1$ are $29$, $73$, $89$, $97$; for all of them, $k = 1$ as we remarked in Section \ref{qs}. The first few primes $p$ with $k \ge 2$ are much larger -- $509$, $797$, $953$; once again, in fact $k = 2$ for all of them.

We now put forward a conjecture:

\begin{conjecture} \label{cc} Let $z$ be a proper Gaussian prime. Then every $4$-neighbour $z$-good polyomino is determined by a clover. \end{conjecture} 

We already know that every $4$-neighbour $z$-good polyomino is determined by some descendant chain $U$. If each neighbourhood vector $(x, y)$ of $U$ satisfies $\{|x|, |y|\} = \{a, b\}$, then reasoning as in the proofs of Propositions \ref{un} and \ref{zgc} we can show without too much difficulty that $U$ is indeed a clover. Thus Conjecture \ref{cc} boils down to ruling out the scenario where the tiling lattices of $U$ coincide with $\mathcal{L}(z)$ and $\mathcal{L}(\overline{z})$ but its neighbourhood vectors are ``irregular''.

Furthermore, it appears that a significantly more general proposition might be true:

\begin{conjecture} \label{uc} The proto-descent of a descendant chain is determined uniquely, modulo $f$-commutativity, by the unordered pair of its tiling lattices. \end{conjecture} 

Clearly, Conjecture \ref{uc} implies Conjecture \ref{cc} immediately. Proposition \ref{un} may be viewed as a (much) weaker variant of Conjecture \ref{uc} where, in addition to the two tiling lattices, the generating vectors of these lattices have been fixed as well.

\section{Further Work} \label{further}

Here, we collect some questions suggested by the material of Sections \ref{intro}--\ref{prime}.

Sections \ref{ptmp} and \ref{agse} are about different aspects of the association between words and polylines. Then in Section \ref{curve} we go beyond polylines and we briefly touch upon continuous curves. Finally, in Section \ref{six} we discuss the general double tile problem -- without a restriction on the number of neighbours.

\subsection{Paradoxical Templates} \label{ptmp}

We know from Section \ref{abc} that there exist some templates which do not contain combinatorial self-intersections, and yet cannot be realised without geometric self-intersections. Let us call such templates \emph{paradoxical}. What can we say about them? Can we describe all of them?

We have already exhibited one example of a paradoxical template. We can obtain others from it by substitution -- but genuinely new paradoxical templates would be more interesting.

Our proof of $\tmpcube$ being paradoxical can be generalised to show that every realisation of $\mathsf{A}\mathsf{B}\mathsf{C}\mathsf{D}\mathsf{A}^{-1}\mathsf{D}^{-1}\mathsf{C}^{-1}\mathsf{B}^{-1}$ with $\sigma(\mathsf{B}) + \sigma(\mathsf{D}) = \mathbf{0}$ contains a geometric self-intersection. This template is not quite paradoxical because of the constraint imposed on its realisations. However, we can obtain new paradoxical templates from it by substituting $\mathsf{B}$ and $\mathsf{D}$ with two words $B$ and $D$ chosen so that the constraint will always be satisfied automatically. The pairs which work are exactly the ones where the concatenation $BD$ is balanced. The template $\tmpcube$ is a special case; it is recovered by the substitution $[\mathsf{B} \to \mathsf{B}, \mathsf{D} \to \mathsf{B}^{-1}]$.

\subsection{The Alphabetisability of General Systems of Equations} \label{agse}

Consider a system of equations $\mathcal{S}$ where each side of each equation is the concatenation of several variables and reverses of variables. Then $\mathcal{S}$ makes sense both over words and over floating polylines. So it is natural to ask how the syntactic and the geometric solutions of $\mathcal{S}$ relate to one another.

Just as in Section \ref{alph}, we define a floating-polyline solution of $\mathcal{S}$ to be \emph{alphabetisable} when we can find a word solution of $\mathcal{S}$ which evaluates to it. What can we say about the alphabetisability of the floating-polyline solutions of $\mathcal{S}$?

E.g., if $\mathcal{S}$ is one of (\textbf{I}) and (\textbf{NI}), then as in Section \ref{adc} we get that all of the non-alphabetisable floating-polyline solutions of $\mathcal{S}$ are flat. Let us call such $\mathcal{S}$ \emph{flattening}. Is this property common? Can it be recognised algorithmically?

Not all systems are flattening. For example, consider the single-equation system $ABCD = BADC$. It admits a lot of floating-polyline solutions which are neither flat nor alphabetisable. However, it does so for trivial reasons -- it factors into the two independent single-equation systems $AB = BA$ and $CD = DC$; accordingly, we can mix and match floating-polyline solutions between the two factors. Of course, indecomposable examples would be of greater interest. (On the other hand, each one of the two factor-systems is indeed flattening. This can be shown by an argument similar to the one in Section \ref{adc}.)

Of particular importance in this context are the systems (\textbf{I}$^6$), (\textbf{NI}$^6_1$), (\textbf{NI}$^6_2$) of Section \ref{six} which govern the general double tile problem.

\subsection{Continuous Curves} \label{curve}

Our double tiles thus far have been polyominoes and polygons. Yet broader is the setting of all regions in the plane bounded by ``well-behaved'' closed continuous curves. What can we say about this variant of the double tile problem? Once again, it would be natural to begin with the $4$-neighbour case. To what extent do our results on polyominoes and polygons carry over?

Notice that the notion of a curve being well-behaved requires formalisation, and the answers to these questions might depend on the particular formalisation we choose.

We can now talk about ``descendant chains'' with infinite descents whose parts are continuous curves instead of polylines. Here follows a brief informal sketch:

Let $\bfphi^\infty = (\ldots, \varphi_{-3}, \varphi_{-2}, \varphi_{-1})$ be a left-infinite sequence of transformations drawn out of~$\psimain$. Let also $\bfphi^{(k)}$ be the $k$-th suffix of $\bfphi^\infty$; i.e., $\bfphi^{(k)} = (\varphi_{-k}, \varphi_{-k + 1}, \ldots, \varphi_{-1})$. Then, experimentally, the descendant chains of polylines $U^{(k)} = \bfphi^{(k)}(\maltese)$ seem to approach a limiting shape when $k$ tends to infinity.

\begin{figure}[ht] \centering \includegraphics{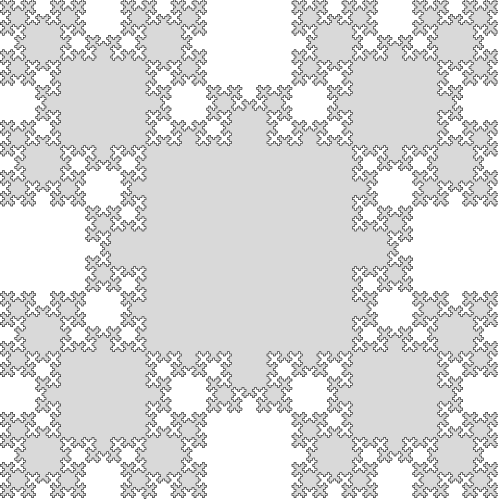} \caption{} \label{fractal} \end{figure}

Intuitively, a left-infinite $\bfphi^\infty$ is the right choice here because the transformations which occur later in the descent have a stronger effect on the shape of the descendant chain. With a left-infinite $\bfphi^\infty$, the last few terms of $\bfphi^{(k)}$ remain constant as $k$ grows; hence, the overall shape of $U^{(k)}$ remains roughly the same, too, and it is only the details that are refined. By contrast, if we let $\bfphi^{(k)}$ be the $k$-th prefix of a right-infinite $\bfphi^\infty$, then in general the last few terms of $\bfphi^{(k)}$ will vary wildly with the growth of $k$; hence, so will the shape of $U^{(k)}$ as well.

Let $U^\infty$ be a chain of continuous curves in the limiting shape of our sequence of descendant chains of polylines. Suppose that $U^\infty$ is sufficiently ``well-behaved'' and consider the region $\Phi^\infty = [U^\infty]$ enclosed by it. Then -- provided that the two tiling lattices of $U^\infty$ are distinct -- $\Phi^\infty$ will be a $4$-neighbour double tile in the continuous curve setting. For example, Figure \ref{fractal} approximates the region which corresponds to $\bfphi^\infty = (\ldots, \gso, \gso, \gso)$. We remarked in the introduction that $4$-neighbour double tiles are ``fractal-like''; notice that this region is a literal fractal, in the sense that it exhibits literal self-similarity.

What can we say about the regions determined by ``descendant chains'' with infinite descents? For example: Which ones of them are topological disks? The one approximated in Figure \ref{fractal} is not; both it and its exterior are ``pinched'' in infinitely many places. On the other hand, with $\bfphi^\infty = (\ldots, f_1, f_3, f_1, f_3, f_1, f_3)$ we do obtain a topological disk -- but it is an ordinary square.

\subsection{The 6-Neighbour Double Tile Problem} \label{six}

What can we say about double tiles in full generality -- without a restriction on the number of neighbours? Do they admit a description similar to the $4$-neighbour one of Theorems \ref{mt} and \ref{pt}?

Notice that the general double tile problem is essentially equivalent to the $6$-neighbour double tile problem because the lattice tilings in the square-grid combinatorial-isomorphism class can be viewed as a degenerate sub-species of the lattice tilings in the hexagon-grid combinatorial-isomorphism class. So we do not distinguish between the $6$-neighbour case and the general case.

Some of the machinery we developed for the $4$-neighbour case generalises to the $6$-neighbour case -- but not all of it. We proceed now to review the similarities and the differences.

To begin with, let \[\tmphex = \mathsf{A}\mathsf{B}\mathsf{C}\mathsf{A}^{-1}\mathsf{B}^{-1}\mathsf{C}^{-1}.\]

Then the $6$-neighbour analogue of Lemma \ref{abl} would go as follows: A given polygon admits a $6$-neighbour lattice tiling of the plane if and only if its boundary realises the template $\tmphex$. (Up to a suitable choice of initial point.) This is indeed true; however, the generalisation of the proof involves some subtleties.

\begin{figure}[ht] \null \hfill \begin{subfigure}[c]{50pt} \centering \includegraphics{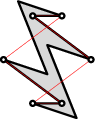} \caption{} \label{hex-a} \end{subfigure} \hfill \begin{subfigure}[c]{125pt} \centering \includegraphics{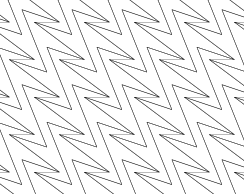} \caption{} \label{hex-b} \end{subfigure} \hfill \null \caption{} \label{hex} \end{figure}

Specifically, suppose that the division points $P_1$, $P_2$, $\ldots$, $P_6$ partition the boundary $S$ of our polygon into the chain $A_1 : B_1 : C_1 : A_2 : B_2 : C_2$ with $A_1 \revsim A_2$, $B_1 \revsim B_2$, $C_1 \revsim C_2$. Let also $T = P_1 \sto P_2 \sto \cdots \sto P_6 \sto P_1$. Our proof of Lemma \ref{abl} derives the desired tiling of $[S]$ from a simpler tiling of $[T]$. But the same strategy cannot go through as-is with $\tmphex$ because in the $6$-neighbour case $T$ can contain self-intersections. One example is shown in Figure \ref{hex}; it has $A_1 = (2, 0) \sto (0, 0)$, $B_1 = (0, 0) \sto (3, -1) \sto (1, 4) \sto (4, 3)$, $C_1 = (4, 3) \sto (0, 6)$.

Fortunately, there is a quick fix: Let $B$ be the translation copy of $B_1$ with initial point $P_4$. Consider the closed polyline $S_\divideontimes = A_1B_1C_1BB^{-1}A_2B_2C_2$. It is obtained from $S$ by the insertion of the double-back detour $BB^{-1}$. Hence, $S$ and $S_\divideontimes$ enclose essentially one and the same region of the plane. However, $S_\divideontimes$ realises the template $\tmptile$ because $A_1B_1 \revsim B^{-1}A_2$ and $C_1B \revsim B_2C_2$. Thus -- in some intuitive sense which is not too difficult to formalise -- Lemma \ref{abl} applies to $[S]$ when its boundary is rewritten as $S_\divideontimes$. We conclude that $[S]$ does admit a lattice tiling of the plane, as desired.

The notion of an interleaved double chain generalises immediately: We say that $U = U_1 : U_2 : \cdots : U_{12}$ is one when it satisfies the system of equations \[U_iU_{i + 1} \reveq U_{i + 6}U_{i + 7} \qquad \text{for all $i$}. \tag{\textbf{I}$^6$}\]

The $f$ and $g$-transforms generalise smoothly, too. The parts of $f_i(U)$, with $i \in \{1, 2, \ldots, 6\}$, are given by $V_j = U_{(j + 5) : (j + 7)}^{-1}$ when $j \equiv i \mymod 6$ and $V_j = U_j$ otherwise. The parts of $\godd(U)$ are given by $V_i = U_{i + 7}\myneg(U_i^{-1})U_{i + 5}$ for all odd $i$ and $V_i = U_i^{-1}$ for all even $i$. The construction rules for the parts of $\geven(U)$ are analogous; just as in the $4$-neighbour case, we obtain them by swapping the words ``odd'' and ``even''.

Let $(a_1, a_2, \ldots, a_6)$ be the type of $U$. We call $U$ a root when $a_{i - 1} = a_{i + 1} = 0$ for some $i$, and we call $U$ a loop when $a_1 + a_3 + a_5 = a_2 + a_4 + a_6$. Once again, every interleaved double chain is an $\{f_1, f_2, \ldots, f_6, \godd, \geven\}$-descendant either of some root or of some loop. However, the structure of the $6$-neighbour roots and loops seems to be significantly more complicated than the structure of their $4$-neighbour analogues.

\begin{figure}[ht] \null \hfill \begin{subfigure}[c]{50pt} \centering \includegraphics{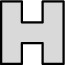} \caption{} \label{h-a} \end{subfigure} \hfill \begin{subfigure}[c]{150pt} \centering \includegraphics{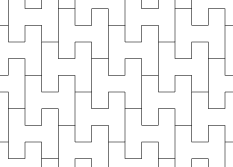} \caption{} \label{h-b} \end{subfigure} \hfill \begin{subfigure}[c]{50pt} \centering \includegraphics{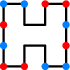} \caption{} \label{h-c} \end{subfigure} \hfill \null \caption{} \label{h} \end{figure}

Two distinct kinds of non-interleaved double chains are possible in the $6$-neighbour case. Using the notation of Section \ref{dcp}, but with $12$ division points instead of $8$, one kind corresponds to the division points occurring in the order $P_1$, $P_2$, $P_3$, $Q_1$, $Q_2$, $Q_3$, $P_4$, $P_5$, $P_6$, $Q_4$, $Q_5$, $Q_6$ and it is characterised by the system \[\begin{aligned} U_i &\reveq U_{i + 6} & &\text{for all $i \not \equiv 0 \mymod 3$}\\ U_{i : (i + 3)} &\reveq U_{(i + 6) : (i + 9)} & &\text{for all $i \equiv 0 \mymod 3$}; \end{aligned} \tag{\textbf{NI}$^6_1$}\] whereas the other kind corresponds to the division points occurring in the order $P_1$, $Q_1$, $P_2$, $P_3$, $Q_2$, $Q_3$, $P_4$, $Q_4$, $P_5$, $P_6$, $Q_5$, $Q_6$ and it is characterised by the system \[\begin{aligned} U_{1 : 2} &\reveq U_{7 : 8} & U_{2 : 4} &\reveq U_{8 : 10} & U_3 &\reveq U_9\\ U_{4 : 6} &\reveq U_{10 : 12} & U_5 &\reveq U_{11} & U_{6 : 7} &\reveq U_{12 : 1}. \end{aligned} \tag{\textbf{NI}$^6_2$}\]

Recall that in the $4$-neighbour case we were able to rule out the non-interleaved double chains early on, with Lemma \ref{nil}. By contrast, in the $6$-neighbour case there exist double tiles associated with non-interleaved double chains. For example, consider the H-shaped heptomino of Figure~\ref{h-a}; its boundary is described by the word $\mathsf{R}\mathsf{U}\mathsf{R}\mathsf{D}\mathsf{R}\mathsf{U}^3\mathsf{L}\mathsf{D}\mathsf{L}\mathsf{U}\mathsf{L}\mathsf{D}^3$. It admits exactly two tilings of the plane, both of which are $6$-neighbour lattice tilings. Figure \ref{h-b} shows one of them, and the other one is its mirror image. (Of course, by means of a suitable deformation we could easily obtain a different example where the two tilings are distinct with respect to all isometries.) The corresponding non-interleaved double chain is shown in Figure \ref{h-c}.

So, in the $6$-neighbour case, we must necessarily find some way to analyse the non-interleaved double chains as well. However, the $f$ and $g$-transforms apply only to the interleaved double chains. One way around this problem could be to construct, for each non-interleaved double chain $U$, an interleaved double pseudo-chain $U^\text{Int}$ which ``mimics'' $U$. This is done, roughly speaking, by pretending that the ``combinatorial division points'' of $U$ are in the right order for it to be interleaved and then filling in negative words whenever we must go backwards from one ``combinatorial division point'' to the next.

Explicitly, if $U$ satisfies (\textbf{NI}$^6_1$), we define $U^\text{Int}$ by \begin{align*} U^\text{Int}_1 &= U_{1 : 3} & U^\text{Int}_2 &= \widehat{U_{2 : 3}} & U^\text{Int}_3 &= U_{2 : 4}\\ U^\text{Int}_4 &= \widehat{U_{3 : 4}} & U^\text{Int}_5 &= U_{3 : 5} & U^\text{Int}_6 &= U_6; \end{align*} and the formulas for the remaining parts of $U^\text{Int}$ are obtained by shifting all indices by $6$ positions, cyclically. Similarly, if $U$ satisfies (\textbf{NI}$^6_2$), we set \begin{align*} U^\text{Int}_1 &= U_1 & U^\text{Int}_2 &= U_2 & U^\text{Int}_3 &= U_{3 : 4}\\ U^\text{Int}_4 &= \widehat{U_4} & U^\text{Int}_5 &= U_{4 : 5} & U^\text{Int}_6 &= U_6. \end{align*}

It is straightforward to see that $U^\text{Int}$ is indeed an interleaved double pseudo-chain. Furthermore, the concatenation $U^\text{Int}_{1 : 12}$ of $U^\text{Int}$ coincides with the concatenation $U_{1 : 12}$ of $U$. We can now apply $f$ and $g$-transforms to $U^\text{Int}$ purely formally, in terms of imaginary words and imaginary concatenation.

\section*{Acknowledgements}

The author is thankful to Palmer Mebane for sharing the original question out of which this work developed.

The present paper was written in the course of the author's PhD studies under the supervision of Professor Imre Leader. The author is thankful also to Prof.\ Leader for his unwavering support.

\end{document}